\documentclass[english,11pt,leqno]{amsart}

\usepackage{newtxtext,newtxmath,tikz,pdfsync}
\usepackage{anyfontsize}

\usepackage{amssymb,stmaryrd}
\usetikzlibrary{cd}

\numberwithin{equation}{section}

\usepackage{quoting}

\usepackage{titlesec,titletoc}

\titleformat{\section}
    {\bfseries\scshape}
    {\thesection.}
    {1ex}
    {\centering}
\titleformat{name=\section,numberless}
    {\bfseries\scshape}
    {}
    {0em}
    {\centering}
\titleformat{\subsection}
    {\bfseries}
    {\thesubsection.}
    {1ex}
    {\raggedright}
\titleformat{name=\subsection,numberless}
    {\bfseries}
    {}
    {0em}
    {\raggedright}

\titlecontents{section}
    [0em]                                       
    {\hspace{0em}\rmfamily \hspace{0ex}}                   
    {\hspace{0em}\thecontentslabel.\hspace{1ex}}        
    {}                              
    {\titlerule*[1em]{.}\thecontentspage}

\titlecontents{subsection}
    [0em]                                       
    {\hspace{0em}\rmfamily \hspace{0ex}}                   
    {\hspace{3ex}\thecontentslabel.\hspace{1ex}}        
    {}                              
    {\titlerule*[1em]{.}\thecontentspage}

\usepackage[utf8]{inputenc}
\usepackage{scrextend}
\usepackage{enumitem}
\setlist[enumerate,1]{label=(\arabic*), leftmargin = 7ex, itemsep=0.5ex,topsep=0.5ex,partopsep=0ex,parsep=0ex}   
\setlist[enumerate,2]{label=(\alph*), leftmargin = 4ex, itemsep=0.5ex,topsep=0.5ex,partopsep=0ex,parsep=1ex}
\setlist[enumerate,3]{label=(\roman*),leftmargin = 4ex, itemsep=0.5ex,topsep=0.5ex,partopsep=0ex,parsep=1ex}

\usepackage[utf8]{inputenc}
\usepackage{amsfonts}
\usepackage{url}
\usepackage[all,cmtip]{xy}
\usepackage{color}
\usepackage{graphicx}
\usepackage{epstopdf} 
\usepackage{time}
\usepackage{amsfonts}
\usepackage[all]{xy}
    \usetikzlibrary{decorations.pathmorphing}
    \usetikzlibrary{patterns}
\usepackage{amsthm}

\font\calli=rsfs10  at 11pt  
\font\callig=rsfs10 at 8 pt     
\font\calliBig=rsfs10 at 9 pt   

\usepackage[T1]{fontenc}
\RequirePackage[colorlinks=true, breaklinks=true, urlcolor= green, linkcolor= blue, citecolor= green, bookmarksopen=true,linktocpage=true,plainpages=false,pdfpagelabels]{hyperref}

\definecolor{better-green}{rgb}{0,.6,.1}
\definecolor{violet}{rgb}{.5,0,.5}

\newcommand{\DEF}[1]{\textbf{\emph{#1}}}    
\newcommand{\tld}[1]{\widetilde{#1}}

\newcommand{\CC}{\mathbb{C}}
\newcommand{\FF}{\mathbb{F}}
\newcommand{\KK}{\mathbb{K}}
\newcommand{\LL}{\mathbb{L}}
\newcommand{\NN}{\mathbb{N}}

\newcommand{\QQ}{\mathbb{Q}}
\newcommand{\RR}{\mathbb{R}}

\newcommand{\VV}{\mathbb{V}}
\newcommand{\ZZ}{\mathbb{Z}}

\newcommand{\A}{\mathcal{A}}
\newcommand{\B}{\mathcal{B}}
\newcommand{\C}{\mathcal{C}}
\newcommand{\D}{\mathcal{D}}
\newcommand{\E}{\mathcal{E}}
\newcommand{\F}{\mathcal{F}}

\newcommand{\J}{\mathcal{J}}
\newcommand{\M}{\mathcal{M}}
\newcommand{\N}{\mathcal{N}}
\newcommand{\R}{\mathcal{R}}
\newcommand{\X}{\mathcal{X}}

\renewcommand{\d}{\mathrm{d}}
\renewcommand{\L}{\mathrm{L}}

\newcommand{\f}{\mathfrak{f}}

\newcommand{\FT}{\hbox{\calli F}}
\newcommand{\FTs}{\hbox{\callig F}}
\newcommand{\RFT}{\hbox{\calli F}_{\mathrm{R}}}

\renewcommand{\S}{\mathcal{S}}

\newcommand{\AS}{\mathbf{AS}}
\newcommand{\CF}{\mathbf{CF}}

\newcommand{\e}{\mathrm{e}}
\renewcommand{\i}{\mathrm{i}}
\DeclareMathOperator{\RE}{Re}
\DeclareMathOperator{\IM}{Im}

\DeclareMathOperator{\an}{an}
\DeclareMathOperator{\Supp}{Supp}
\DeclareMathOperator{\ST}{S}
\DeclareMathOperator{\AF}{A}
\newcommand{\Res}[1]{\underset{#1}{\mathrm{Res}}\,}
\DeclareMathOperator{\cl}{cl}

\DeclareMathOperator{\dist}{dist}
\DeclareMathOperator{\SPAN}{span}
\DeclareMathOperator{\CONV}{conv}
\DeclareMathOperator{\VERT}{vert}

\DeclareMathOperator{\Arg}{Arg}
\DeclareMathOperator{\Log}{Log}
\DeclareMathOperator{\LE}{LE}
\DeclareMathOperator{\LC}{LC}
\DeclareMathOperator{\LM}{LM}
\DeclareMathOperator{\LF}{LF}
\DeclareMathOperator{\LT}{LT}

\newtheorem{thm}{Theorem}[section]

\newtheorem{cor}[thm]{Corollary}

\newtheorem{prop}[thm]{Proposition}
\newtheorem{lem}[thm]{Lemma}

\theoremstyle{definition}
\newtheorem{defn}[thm]{Definition}
\newtheorem{defs}[thm]{Definitions}
\newtheorem{notn}[thm]{Notation}

\newtheorem{que}[thm]{Question}

\newtheorem{rem}[thm]{Remark}
\newtheorem{exam}[thm]{Example}
\newtheorem*{exam*}{Example}
\newtheorem{rems}[thm]{Remarks}
\newtheorem{exams}[thm]{Examples}

\newtheoremstyle{gen*:sc:rm}
  {10pt}
  {10pt}
  {\rm}
  {}
  {\scshape}
  {.}
  {.5em}
  {\thmnote{#3}}  

\theoremstyle{gen*:sc:rm}
\newtheorem{plain*}{}

\newcommand*\latinnumeral[1]
  {\ifcase#1\unskip\or\unskip\or bis\fi}
\newcommand*\latintags
  {\xdef\startingequationnumber{\the\numexpr\value{equation}+1\relax}%
   \setcounter{equation}{0}%
   \aftergroup\resetarabicequations
   \def\theequation
     {\startingequationnumber~\textit{\latinnumeral{\value{equation}}}}}
\newcommand*\resetarabicequations
  {\setcounter{equation}{\numexpr\startingequationnumber\relax}}

\usepackage{xcolor}

\title[Fourier transforms of power-constructible functions]{Decay of Fourier transforms and analytic continuation of power-constructible functions}

\author{Georges Comte}
\address{Univ.   Savoie Mont Blanc, CNRS, LAMA, 73000 Chamb\'ery, France}
\email{georges.comte@univ-smb.fr}
\urladdr{https://georgescomte.perso.math.cnrs.fr/}
\author{Daniel J. Miller}
\address{Emporia State University, Campus Box 4027, 1 Kellogg Circle, Emporia, KS 66801, U.S.A.}
\email{dmille10@emporia.edu}
\urladdr{}
\author{Tamara Servi}
\address{Institut de Math\'ematiques de Jussieu -- Paris Rive Gauche, Universit\'e Paris Cit\'e and
Sorbonne Universit\'e, CNRS, IMJ-PRG, F-75013 Paris, France}
\email{tamara.servi@imj-prg.fr}

\begin{document}
\date{\today}

\begin{abstract}
For a subfield $\KK$ of $\CC$, we denote by $\C^{\KK}$ the category of algebras of functions defined on the globally subanalytic sets that are generated by all $\KK$-powers and logarithms of positively-valued globally subanalytic functions.  For any function $f\in\C^\KK(\RR)$, we study links between holomorphic extensions of $f$ and the decay of its Fourier transform $\hbox{\calliBig F}[f]$ by using tameness properties of the globally subanalytic functions from which $f$ is constructed.  We first prove a number of theorems about analytic continuation of functions in $\C^{\KK}$, including the fact that $f\in\C^{\KK}(\RR)$ extends meromorphically to $\CC$ if and only if $f$ is rational.  We then characterize the exponential rate of decay of $\hbox{\calliBig F}[f]$ by the maximal width of a horizontal strip in the plane about the real axis to which $f$ extends holomorphically.  Finally, we show that $\hbox{\calliBig F}[f]$ is integrable if $f$ is integrable and continuous.
\end{abstract}

\maketitle
\tableofcontents

\section*{Introduction}

The primary purpose of this paper is to apply techniques from tame geometry to study the decay of the Fourier transform
\begin{equation}\label{eq:intro:FTdef}
\FT[f](\xi) = \int_{\RR} f(t) \e^{-\i 2\pi \xi t} \d t
\end{equation}
as $\xi\to\pm\infty$ when $f$ is a so-called \emph{power-constructible function} in $\L^1(\RR)$, with a particular focus on its exponential decay.  Our main theorem shows that $\FT[f]$ decays exponentially if and only if $f$ is analytic and that, in this situation, the exponential rate of decay of $\FT[f]$ is characterized by the maximal open horizontal strip about the real axis in the complex plane to which $f$ may be extended analytically.  Along the way we also study various related topics, including the rate of decay of the Mellin transform, the integrability of the Fourier transform, and analytic continuation of parametric families of power-constructible functions and a number of other closely related categories of functions.  These various categories are constructed in concrete ways from the globally subanalytic functions, but before delving into the details of their definitions, let us first  motivate their use by discussing the analytic considerations underlying our main theorem.  Throughout this discussion and the rest of the paper, we use the following asymptotic terminology for any function $f:\RR\to\CC$ and sign $\sigma\in\{+,-\}$: we say that $f$ \DEF{decays at $\sigma\infty$} when $\lim_{t\to\sigma\infty} f(t) = 0$, that $f$ \DEF{decays rapidly at $\sigma\infty$} when $f(t) = O(|t|^{-r})$ as $t\to\sigma\infty$ for every $r > 0$, that \DEF{$f$ decays slowly at $\sigma\infty$} when $f$ decays at $\sigma\infty$ but not rapidly, and that \DEF{$f$ decays exponentially at $\sigma\infty$} when $f(t) = O(\e^{-r|t|})$ as $t\to\sigma\infty$ for some $r > 0$.
\\

\noindent
\textbf{Discussion.}
Consider any $f\in\L^1(\RR)$.  Then $\FT[f]$ is continuous and decays at $\pm\infty$ by the Riemann-Lebesgue lemma.  But to study the nature of the decay of $\FT[f]$ at $\pm\infty$ in more detail, we must impose assumptions upon $f$ in addition to its integrability.  One approach might be to assume that $f$ is a Schwartz function (i.e., $f^{(n)}$ decays rapidly for all $n\in\NN$), for then $\FT[f]$ is also a Schwartz function, but this does not address our main concern of whether $\FT[f]$ decays exponentially.   Another approach might be to use classical techniques from asymptotic analysis such as the principle of stationary phase, the saddle point method, or Watson's lemma. (For example, see Sidorov, Fedoryuk, and Shabunin \cite{SidFedSha}; Stein \cite{Stein}; or Wong \cite{Wong}.)  However, it seems impractical to use the saddle point method as a general tool for determining exponential decay of $\FT[f]$, and all three techniques are best suited for obtaining precise asymptotic expansions for $\FT[f]$ in monomial scales, which again does not address whether $\FT[f]$ decays exponentially.

Of greater relevance is the work of Lombardi \cite{Lom}, which uses complex analysis to study the exponential decay of oscillatory integrals of analytic functions arising from the study of certain systems of differential equations.  We do not focus on the same types of functions or integrals as Lombardi, but we do focus on exponential decay of an oscillatory integral (the Fourier transform), and we use complex analysis by basing the proof of our main theorem on two principles loosely described as follows for a given function $f\in\L^1(\RR)$:
\begin{enumerate}[label=(\Roman*)]
\item
If $\FT[f]$ decays exponentially at $\pm\infty$, then one may apply holomorphy under the integral sign to the inverse Fourier transform of $\FT[f]$ to show that $f$ extends analytically to a horizontal strip $\RR+\i(a,b)$ with $a<0<b$, where $a$ and $b$ are respectively determined from bounds on the exponential rates of decay of $\FT[f]$ in the directions $+\infty$ and $-\infty$.  (See Lemma \ref{lem:FTDecayImplyHolomExt}.)

\item
If $f$ extends analytically to a horizontal strip $\RR+\i[a,b]$ with $a < 0 < b$ and has a \emph{suitable decay at $\infty$} in this strip, then we may replace the path of integration $\RR$ in the integral \eqref{eq:intro:FTdef} with either of the parallel lines $\RR + \i a$ or $\RR + \i b$ to shows that $\FT[f]$ decays exponentially at $+\infty$ and $-\infty$ with bounds on the exponential rates of decay respectively determined by $a$ and $b$.  (See Lemma \ref{lem:HolExtImplyFTDecay}.)
\end{enumerate}
Our aim is to apply these principles as a pair to conclude that the optimal bounds on the rates of exponential decay of $\FT[f]$ at $\pm\infty$ and the maximal horizontal strip about the real axis to which $f$ extends analytically completely determine one another.  But this is only guaranteed to hold if we restrict our attention to studying special types of functions $f$ such that if $f$ is integrable on $\RR$ and extends analytically to $\RR+\i[a,b]$ with $a < 0 < b$, then this extension exhibits the suitable decay at $\infty$ alluded to in Principle II.  We look to tame geometry to supply natural categories of rings of functions with this decay property.

In tame geometry one typically studies the collection of all definable sets and functions of some chosen o-minimal structure expanding the real field.  The benefit of this is that, by definability alone, these collections are stable under a host of operations commonly used in analytic and geometric constructions, and o-minimality additionally endows these sets and functions with many tameness properties.  However, the collection of definable functions of an o-minimal structure  is generally not stable under operations defined from integrals, and more to the point, is never stable under the operation of Fourier transform.  Because of this, we depart from this typical approach of tame geometry by instead working with categories of rings of functions on the definable sets of an o-minimal structure that extend the category of rings of its definable functions.  The extension rings are rich enough to be stable under certain operations defined from integrals but are also constructed explicitly enough from the definable functions to retain control on their behavior through the tameness of the underlying o-minimal structure.  For the o-minimal structure chosen in this paper, it turns out we can still retain control over the extension rings if we also include complex powers of its definable functions, so we opt to do so.

Specifically, in this paper we start with the o-minimal structure $\RR_{\an}$, which is the expansion of the real field by all restricted analytic functions.  (See see van den Dries and Miller \cite[Example 2.5(4)]{vdDM:1994}.)  The definable sets and functions of $\RR_{\an}$ are called globally subanalytic, but we abbreviate this by saying \DEF{subanalytic} to mean globally subanalytic.  We write $\S$ for the category whose objects are the rings $\S(X)$ of all real-valued subanalytic functions on $X$, for each subanalytic set $X$, and whose morphism are pullbacks by subanalytic maps, and we call a function in $\S(X)$ for some $X$ an \DEF{$\S$-function}; we also use similar terminology and notation for categories of rings of functions extending $\S$.  Now, fix a subfield $\KK$ of $\CC$, and observe that $\cl_{\CC}(\KK) = \RR$ when $\KK\subseteq\RR$ and that $\cl_{\CC}(\KK) = \CC$ when $\KK\not\subseteq\RR$, where $\cl_{\CC}(\KK)$ is the topological closure of $\KK$ in $\CC$.  We consider a strict hierarchy
\begin{equation}\label{eq:intro:CatHier}
\S^{\KK} \subset \C^{\KK} \subset \C^{\KK,\i\S} \subset \C^{\KK,\FTs}
\end{equation}
of categories of rings on the subanalytic sets described briefly as follows: $\S^{\KK}$ is generated as $\cl_{\CC}(\KK)$-algebras from all $\KK$-powers of positively-valued subanalytic functions, $\C^{\KK}$ is generated from $\S^{\KK}$ and all logarithms of positively-valued subanalytic functions, $\C^{\KK,\i\S}$ is generated from $\C^{\KK}$ and all complex exponentials $f(x) =  \e^{\i g(x)}$ with subanalytic phase functions $g$, and $\C^{\KK,\FTs}$ consists of all parametric integrals $f(x) = \int_{\RR}g(x,t)\e^{\i t}\d t$ with $g$ in $\C^{\KK}$.  A $\C^{\KK}$-function is also called a \DEF{$\KK$-power-constructible function}, and a $\C^{\CC}$-function is more simply called a \DEF{power-constructible function}, which are the titular functions of the paper.

The main interest in $\C^{\KK}$ and $\C^{\KK,\FTs}$ stems from the fact that among all of the categories of rings of functions on the subanalytic sets that contain $\S^{\KK}$, it turns out that $\C^{\KK}$ is the smallest that is stable under parametric integration \cite{ClCoRoSe:preprint}  and that $\C^{\KK,\FTs}$ is the smallest that is stable under parametric Fourier transforms \cite{clcose:fourier_mellin_power_constructible}.  However, each of the four categories in \eqref{eq:intro:CatHier} play a more specific role in this paper: $\S^{\KK}$ serves as the generating category of rings; $\C^{\KK}$ serves as the largest category in the hierarchy \eqref{eq:intro:CatHier} that has the decay property required by Principle II; $\C^{\KK,\i\S}$ serves as a category intermediary to $\C^{\KK}$ and $\C^{\KK,\FTs}$ that shares many tameness properties with $\C^{\KK}$ (apart from said decay property); and $\C^{\KK,\FTs}$ serves as an ambient category whose asymptotic expansions in power-log scales and stability under Fourier transforms we exploit when studying the integrability of Fourier transforms of $\C^{\KK}$-functions.  One may opt to use $\S$ as the generating category by taking $\KK=\QQ$.
\\

\noindent\textbf{Outline.}
We begin our development in Section \ref{s:Background}, which starts by giving precise definitions of the categories in \eqref{eq:intro:CatHier} and their extensions to subanalytic subsets of $\CC^m$ and then gives an extensive exposition of their preparation theorems and asymptotic expansions.  Section \ref{s:ConseqPrepAsym} then deduces various consequences of this machinery that serve as technical lemmas in the rest of the paper.

The main thrust of the paper then begins in Section \ref{s:HolomMeromExt} by studying holomorphic extensions of parametric families of functions.  Simpler nonparametric versions of the section's three main theorems are stated as follows for a given category  $\D\in\{\S^{\KK}, \C^{\KK}, \C^{\KK,\i\S}\}$ and  function $f\in\D(\RR)$:
\begin{enumerate}
\item\label{intro:Main1}{\scshape  (Analytic Extensions at $\pm\infty$).}
For each sign $\sigma\in\{+,-\}$, there exists a $b>0$ such that the restriction of $f$ to $\{y\in\RR: \sigma y > b\}$ is analytic and extends to an analytic $\D$-function on $\{z\in\CC: \text{$|z| > b$ and $|\Arg(\sigma z)| < \pi$}\}$.  (See Theorem \ref{thm:HolomExtInfty}.)

\item\label{intro:Main2}{\scshape (Analytic Extensions to Horizontal Strips).}
The function $f$ is smooth if and only if it is analytic, in which case $f$ extends to an analytic $\D$-function on a horizontal strip $\RR + \i(-a,a)$ for some $a > 0$.  (See Theorem \ref{thm:CKiS:HolomExtStrip}.)

\item\label{intro:Main3}{\scshape (Meromorphic Rationality Theorem).}
If $\D\in\{\S^{\KK},\C^{\KK}\}$, then $f$ extends meromorphically to the plane $\CC$ if and only if $f$ is a rational function with coefficients in $\cl_{\CC}(\KK)$.  (See Theorem \ref{thm:CK:MeromExt}.)
\end{enumerate}
Four points are worth noting.  First, to prove \ref{intro:Main3} and its parametric version Theorem \ref{thm:CK:MeromExt}, we first establish the following Liouville property: a function $f\in\C^{\KK}(\RR)$ extends to an entire function if and only if $f$ is a polynomial with coefficients in $\cl_{\CC}(\KK)$.  Second, a consequence of the parametric version of \ref{intro:Main3} is that if $f\in\D(X\times\RR)$ and $f_x:y\mapsto f(x,y)$ extends meromorphically to $\CC$ for all $x\in X$, then the poles of $f$ in $y$ are piecewise given by subanalytic functions of $x$ on $X$.  Third, \ref{intro:Main3} is similar in spirit to the o-minimal rationality theorem of Peterzil and Starchenko \cite[Theorem 1.2]{PetStar:2001} but is not strictly comparable because when $\i\in\KK$, the category $\S^{\KK}$ contains many oscillatory functions not definable in any o-minimal structure; for example, the function $t\mapsto \sin(\log t) = (t^i-t^{-i})/(2i)$ is in $\S^{\KK}(\RR)$.  And fourth, in \ref{intro:Main1}-\ref{intro:Main3} and their parametric versions, the fact that the holomorphic extensions of the given $\D$-function is also a $\D$-function is similar to previous work on analytic continuation for functions definable in $\RR_{\an}$ and its expansion $\RR_{\an,\exp}$ by the exponential function, such as Wilkie \cite{Wil2016}, Kaiser \cite{Kai2016}, Kaiser and Speissegger \cite{KaiSpe}, and the preprint Opris \cite{Opris}; but \ref{intro:Main1}-\ref{intro:Main3} are again not strictly comparable to these o-minimal results since $\S^{\KK}$ may contain oscillatory functions.

In Section \ref{s:DecayFT} we turn to proving our main theorem characterizing the exponential decay of Fourier transforms of $\C^{\KK}$-functions, which is paraphrased as follows.
\begin{enumerate}{\setcounter{enumi}{3}
\item\label{intro:Main4}{\scshape (Decay of the Fourier Transform).}
Let $f\in\C^{\KK}(\RR)\cap\L^1(\RR)$.  Then up to modifying $f$ at a finite number of points, the four conditions
\begin{enumerate}[label=(\alph*),itemsep=-0.5ex,topsep=0ex]
\item
$\FT[f]$ decays exponentially at $\pm\infty$,

\item
$\FT[f]$ decays rapidly at $\pm\infty$,

\item
$f$ is smooth,

\item
$f$ is analytic,
\end{enumerate}
are equivalent to one another, and the optimum bounds on the exponential rate of decay of $\FT[f]$ at $\pm\infty$ are fully characterized by the maximal open horizontal strip about the real axis in the complex plane to which $f$ may be extended analytically.  (See Theorem \ref{thm:CK:FTexpDecay}.)
}\end{enumerate}
Three points are worth noting.  First, in addition to using Principles I and II to prove \ref{intro:Main4}, we apply these principles to a classical identity relating the Mellin and Fourier transforms to prove an analogous theorem relating the exponential rate of decay of the Mellin transform of a smooth $\C^{\KK}$-function $f\in\L^1(\RR_{\geq 0})$ with the maximal sector to which $f$ may be extended analytically.  Second, the residue theorem can be used to improve upon \ref{intro:Main4} when $f$ extends meromorphically by giving an exact asymptotic for $\FT[f]$ rather than just a bound, and this explains some of the utility of \ref{intro:Main3}.  Specifically, if $f\neq 0$ extends meromorphically to $\CC$, then for each sign $\sigma\in\{+,-\}$, $\FT[f](\xi) \sim C \e^{-2\pi r|\xi|}$ as $\xi\to\sigma\infty$ for some $C\neq 0$, where $r$ is the minimum distance from the real axis of the poles of $f$ in the half plane $\{z\in\CC : \sigma|\IM(z)| < 0\}$.  (See Remarks \ref{rems:DecayFT:Sing}\ref{rems:DecayFT:Sing:Poles} and Lemma \ref{lem:DecayFT:CKmeromFT}.)  And third, given a function $f\in\C^{\KK}(X\times\RR)$ with $f_x\in\L^1(\RR)\cap C^\infty(\RR)$ for all $x\in X$, we combine \ref{intro:Main4} with the parametric version of \ref{intro:Main2} to conclude that there exists a function $a\in\S_+(X)$ such that $\FT[f_x](\xi) = O(\e^{-a(x)\xi})$ as $\xi\to\pm\infty$ for all $x\in X$.  However, this parametric asymptotic bound is not optimal.  The main obstacle to optimality in the parametric case is that we do not currently have a good understanding of the \emph{virtual complex singularities} of a general function $f\in\C^{\KK}(X\times\RR)$, by which we mean the singularities of holomorphic extensions of $\{f_x\}_{x\in X}$,  except in the extreme situation when these functions extend meromorphically to all of $\CC$.  This remains a deep challenge, and in some sense we approach this problem here through the exponential rate of decay of the Fourier and Mellin transforms at infinity.

Finally, Section \ref{s:IntegFT} concludes the paper by studying the integrability of Fourier transforms.  The Fourier inversion theorem shows that when $f\in\L^1(\RR)$, a necessary condition for $\FT[f]$ to be integrable is that $f$ is continuous, up to almost everywhere equivalence.  We use properties of $\C^{\KK,\FTs}$ (particularly its stability under Fourier transforms and its asymptotic expansions) and properties of $\C^{\KK}$ to show that continuity of $f$ is actually a sufficient condition when $f$ is a $\C^{\KK}$-function.
\begin{enumerate}{\setcounter{enumi}{4}
\item\label{intro:Main5}{\scshape (Integrability of the Fourier Transform).}
If $f\in\C^{\KK}(\RR)\cap\L^1(\RR)$ is continuous, then $\FT[f]\in\L^1(\RR)$.  (See Theorem \ref{thm:IntegFT:CK}.)
}\end{enumerate}
This extends \cite[Theorem 1.3]{CluMil}, which is the special case of \ref{intro:Main5} for $\KK=\QQ$.

In addition to proving its various theorems, Sections \ref{s:HolomMeromExt}-\ref{s:IntegFT} also expend a great deal of effort giving counterexamples that demonstrate the failure of each theorem in the next larger category in the hierarchy \eqref{eq:intro:CatHier} than hypothesized in the theorem, but with a notable exception: we do not know whether \ref{intro:Main5} holds with $\C^{\KK,\i\S}$ in place of $\C^{\KK}$.  However, we do show that \ref{intro:Main5} fails with $\C^{\KK,\FTs}$ in place of $\C^{\KK}$.
\vspace*{1ex}

\noindent\textbf{Acknowledgment.}
The authors are grateful to Raf Cluckers for some very helpful conversations and preliminary work that led us to \ref{intro:Main5} and that motivated our search for a counterexample to this theorem in the larger category $\C^{\KK,\FTs}$.

\section{Background}\label{s:Background}

This section is largely a review of the supporting literature that is formulated for our current purposes.  After briefly fixing some notation to be used throughout the paper, the section focuses on the categories of functions in the hierarchy $\S^{\KK}\subset\C^{\KK}\subset\C^{\KK,\i\S}\subset\C^{\KK,\FTs}$ by discussing their precise definitions, preparation theorems,  and asymptotic expansions. 

\begin{notn}\label{notn:complex}
Write $\RE(z)$ and $\IM(z)$ for the real and imaginary parts of a complex number $z$.  For each $r>0$ and $z_0\in\CC$, write
\[
\text{$D_r(z_0) = \{z\in\CC:|z-z_0|\leq r\}$ and $D_{r}^{*}(z_0) = \{z\in\CC:0<|z-z_0|\leq r\}$.}
\]
For any nonzero $z\in\CC$, write $\Arg(z)\in(-\pi,\pi]$ and $\Log(z) = \log|z| + \i\Arg(z)$ for the principal values of the argument and logarithm of $z$.  For each $\sigma\in\{+,-\}$, define
\[
\RR_{\sigma} = \{t\in\RR : \sigma t > 0\}
\]
and
\[
\CC_{\sigma} = \{z\in\CC\setminus\{0\} : |\Arg(\sigma z)| < \pi\},
\]
and write $A_\sigma = A\cap\RR_\sigma$ for any $A\subseteq\RR$.  For any $\gamma\in\CC$, unless stated otherwise, one should assume that $z\mapsto z^\gamma$ is defined as the principal power function $z^\gamma := \e^{\gamma\Log(z)}$ on the set $\CC_+$.  For each $\sigma\in\{+,-\}$, define the closed half plane
\[
H_{\sigma} = \{z\in\CC : \sigma\IM(z) \geq 0\}.
\]
When discussing paths of integration of contour integrals, for each $a,b\in\CC$ we write $[a,b]$ for the directed line segment in the complex plane from $a$ to $b$.  For any $A\subseteq\CC^m$ and $B\subseteq\CC$, we call a function $f:A\to B$ \DEF{holomorphic} (or \DEF{analytic}) when $f$ extends to a holomorphic function $f:U\to\CC$ on some open neighborhood $U$ of $A$ in $\CC^m$.
\end{notn}

Also, throughout the paper, the word \DEF{integrable} without any further clarification means Lebesgue integrable, and $\NN = \{0,1,2,\ldots\}$.

\subsection{The Categories of Functions}\label{ss:Background:CategFcts}

Recall from the Introduction that we use the word \emph{subanalytic} to mean definable in $\RR_{\an}$ and that for each subanalytic set $X$, we write $\S(X)$ for the ring of all real-valued subanalytic functions on $X$.  Throughout Section \ref{s:Background}, we write $\Omega$ for the set of all subanalytic sets (each of which is contained in $\RR^m$ for some $m\in\NN$, by definition), and we write $\S = \{\S(X)\}_{X\in\Omega}$.  We also follow the convention that a subanalytic map is required to have a subanalytic codomain.

\begin{defs}\label{defs:Background:CategSub}
By a \DEF{category of functions on the subanalytic sets} we mean a family $\D = \{\D(X)\}_{X\in\Omega}$ satisfying the following three properties.
\begin{enumerate}
\item\label{defs:Background:CategSub:subanal}
For each $X\in\Omega$, $\D(X)$ is a set of complex-valued functions on $X$ containing $\S(X)$.

\item\label{defs:Background:CategSub:comp}
For every subanalytic map $\varphi:Y\to X$ and $f\in\D(X)$, we have $f\circ\varphi\in\D(Y)$.

\item\label{defs:Background:CategSub:union}
For all $n\in\NN$, disjoint subanalytic subsets $X$ and $Y$ of $\RR^n$, and functions $f\in\D(X)$ and $g\in\D(Y)$, we have $f\cup g\in\D(X\cup Y)$.
\end{enumerate}
Two degenerate choices of $X$ require special mention.  First, we define $\D(\emptyset) = \{\emptyset\}$ because the only function on the empty set is the empty map.  And second, we identity $\D(\RR^0)$ with a subset of $\CC$ because we identify $\RR^0$ with $\{0\}$ and identify each function $f:\RR^0\to\CC$ with the number $f(0)$.

Call a function in $\D(X)$, for some $X\in\Omega$, a \DEF{$\D$-function}.  For any subanalytic map $\varphi:Y\to X$, we define a pullback map $\varphi^*:\D(X)\to\D(Y)$ by $\varphi^*(f) = f\circ\varphi$ for each $f\in\D(X)$.  We therefore obtain a small category of sets whose objects are the members of $\D$ and whose morphisms are the pullback maps of subanalytic maps.  If each object of $\D$ is closed under the operations of addition and multiplication, then the objects of $\D$ are commutative rings and the morphism of $\D$ are ring homomorphisms, and in this situation we call $\D$ a \DEF{category of rings} of functions on the subanalytic sets.

We also need complex versions of these notions when $\D$ is a category of rings.  For each $\FF\in\{\RR,\CC\}$, define $\D_{\FF} = \{\D_{\FF}(X)\}_{X\in\Omega}$, where $\D_{\FF}(X)$ is the $\FF$-algebra of functions on $X$ generated from $\D(X)$.  Clearly, $\D_{\RR} = \D$ because $\RR=\S(\RR^0)\subseteq\D(\RR^0)$, and we have $\D\subseteq\D_{\CC}$ with equality holding if and only if $\D(\RR^0) = \CC$.  For each $m\in\NN$ and $X\subseteq\CC^m$, we call $X$ subanalytic when $X$ is a subanalytic subset of $\RR^{2m}$ under the canonical identification $\CC^m \cong \RR^{2m}$.  When $X\subseteq\CC^m$ is subanalytic and $f\in\D_{\CC}(X)$, we call $f$ a \DEF{$\D$-function} provided that $f\restriction X\cap\RR^m \in \D(X\cap\RR^m)$. 
\end{defs}

\begin{notn}\label{notn:Background:CoordProjFiber}
We use the following notation when given any $m,n\in\NN$.  We write $(x,y) = (x_1,\ldots,x_m,y_1,\ldots,y_n)$ for a tuple of variables varying over $\RR^{m+n}$, and we define the projection $\pi:\RR^{m+n}\to\RR^m$ by $\pi(x,y) = x$ for all $(x,y)\in\RR^{m+n}$.  For each set $A\subseteq\RR^{m+n}$, function $f:A\to\CC$, and point $x\in\RR^m$, define the fiber $A_x = \{y\in\RR : (x,y)\in A\}$, and define the function $f_x:A_x\to\CC$ by $f_x(y) = f(x,y)$ for each $y\in A_x$.  When $n = 1$, we simply write $(x_1,\ldots,x_m,y)$ for $(x,y)$.
\end{notn}

\begin{defs}\label{def:Background:CategStab}
Assume that $\D$ is a category of functions on the subanalytic sets.  We say that $\D$ is \DEF{stable under parametric integration} when for each $m,n\in\NN$, subanalytic set $X\subseteq\RR^m$, and function $f\in \D(X\times\RR^n)$ with $f_x\in \L^1(\RR^n)$ for each $x\in X$, the function $F:X\to\RR$ defined by
\[
F(x) = \int_{\RR^n} f(x,y)\d y, \quad\text{for each $x\in X$,}
\]
is in $\D(X)$, where $\d y = \d y_1\wedge \cdots \wedge \d y_n$.  We say that $\D$ is \DEF{stable under parametric Fourier transforms} when for each such $f$ (as hypothesized above), the function $\FT_X[f]:X\times\RR^n\to\CC$ defined by
\[
\FT_X[f](x,\xi)
= \int_{\RR^n} f(x,y) \e^{-2\pi\i y\cdot\xi}\d y, \quad\text{for each $(x,\xi)\in X\times\RR^n$,}
\]
is in $\D(X\times\RR^n)$, where $y\cdot\xi$ is the Euclidean inner product.
\end{defs}

Throughout the rest of the paper, fix a subfield $\KK$ of $\CC$.

\begin{defs}\label{def:intro:OurCategFcts}
For each $X\in\Omega$,
\begin{enumerate}
\item
define $\S_+(X) = \{f\in\S(X) : f > 0\}$,

\item
define $\S^{\KK}(X)$ to be the $\cl_{\CC}(\KK)$-algebra of functions on $X$ generated by
\[
\{f^\gamma : \text{$f\in\S_+(X)$ and $\gamma\in\KK$}\},
\]

\item
define $\C^{\KK}(X)$ to be the ring of functions on $X$ generated by
\[
\S^{\KK}(X)\cup \{\log(f) : f\in\S_+(X)\},
\]

\item
define $\C^{\KK,\i\S}(X)$ to be the ring of functions on $X$ generated by
\[
C^{\KK}(X)\cup \{\e^{\i f} : f\in\S(X)\},
\]

\item
and define
\[
\qquad \C^{\KK,\FTs}(X) = \{\f[g] : \text{$g\in\C^{\KK}(X\times\RR)$ and $g_x\in\L^1(\RR)$ for all $x\in X$}\},
\]
where
\[
\f[g](x) = \int_{\RR} g(x,y)\e^{\i y}\d y \quad\text{for each $x\in X$.}
\]
\end{enumerate}
For each symbol $\D\in\{\S^{\KK},\C^{\KK},\C^{\KK,\i\S},\C^{\KK,\FTs}\}$, write $\D := \{\D(X)\}_{X\in\Omega}$.
\end{defs}

\begin{rems}\label{rem:Categs}
\hfill
\begin{enumerate}
\item\label{rem:Categs:SK}
We have $\S\subseteq\S^{\KK}$, and equality holds if and only if $\KK = \QQ$.

\begin{proof}
For any subanalytic set $X$ and $f\in\S(X)$, the functions
\[
\qquad
g(x) = \begin{cases}
|f(x)|, & \text{if $f(x) \neq 0$,} \\
1,      & \text{if $f(x) = 0$,}
\end{cases}
\quad\text{and}\quad
h(x) = \begin{cases}
3,  & \text{if $f(x) > 0$,} \\
2,  & \text{if $f(x) = 0$,} \\
1,  & \text{if $f(x) < 0$,} \\
\end{cases}
\]
are in $\S_+(X)$, so $f = g(h-2)\in\S^{\QQ}(X)$.  This shows that $\S \subseteq \S^{\QQ}$.  The inclusions $\S^{\QQ}\subseteq\S^{\KK}$ and $\S^{\QQ}\subseteq\S$ are clear, so $\S = \S^{\QQ}\subseteq\S^{\KK}$.  When $\KK\neq\QQ$, we have $\S\neq\S^{\KK}$ because $\QQ$ is the set of definable exponents of $\RR_{\an}$ (for instance, by C. Miller \cite[Corollary 4.7]{ChrisMil:1994} or by Proposition \ref{prop:SubPrep}).
\end{proof}

\item\label{rem:Categs:Hierarchy}
We have strict inclusions
\begin{equation}\label{eq:Background:OurCategs}
\S^{\KK} \subset \C^{\KK} \subset \C^{\KK,\i\S} \subset \C^{\KK,\FTs}.
\end{equation}

\begin{proof}
We have $\S^{\KK}\subseteq\C^{\KK}\subseteq\C^{\KK,\i\S}$ by definition, and both inclusions are strict because the function $t\mapsto\log(t)$ is in $\C^{\KK}(\RR_+)\setminus\S^{\KK}(\RR_+)$ and the function $t\mapsto\e^{\i t}$ is in $\C^{\KK,\i\S}(\RR)\setminus\C^{\KK}(\RR)$.  The remainder of the proof relies on facts that are implicit in the literature but which we shall develop in detail in the remainder of the section and a little bit of the next section.  So the reader may wish to return to this proof after completing Section \ref{s:Background} and the proof of Proposition \ref{prop:CFK:anal}.

Let $X\subseteq\RR^m$ be subanalytic, and let us temporarily write $\C^{\KK,\exp}(X)$ for the $\C^{\KK,\i\S}(X)$-module of functions on $X$ generated by
\begin{equation}\label{eq:gammaFct}
\{\gamma_{h,\ell} : \text{$\ell\in\NN$ and $h\in\S^{\KK}(X\times\RR)$ s.t. $h_x\in \L^1(\RR)$ for all $x\in X$}\},
\end{equation}
where
\[
\gamma_{h,\ell}(x) = \int_{\RR} h(x,y)(\log |y|)^\ell \e^{\i y} \d y
\quad\text{for all $x\in X$.}
\]
We claim that $\C^{\KK,\exp}(X) = \C^{\KK,\FTs}(X)$.  Indeed, the inclusion $\C^{\KK,\exp}(X) \subseteq \C^{\KK,\FTs}(X)$ follows from the proof of \cite[Lemma 2.6]{AiClRaSe}, and the inclusion $\C^{\KK,\FTs}(X) \subseteq \C^{\KK,\exp}(X)$ follows from starting with any $f = \f[g] \in \S^{\KK,\FTs}(X)$ and applying Proposition \ref{prop:CKprep} to prepare $g$ to see that $f\in\C^{\KK,\exp}(X)$ (similar to what is done in the proof of Proposition \ref{prop:CFK:anal}).  The constant function $1:X\to\RR$ is in the set \eqref{eq:gammaFct} by \cite[Remark 2.6]{ClCoRoSe}, so $\C^{\KK,\i\S}\subseteq \C^{\KK,\FTs}$ by the claim, and this inclusion is strict because $t\mapsto \e^{-|t|}$ is in $\C^{\KK,\FTs}(\RR)\setminus\C^{\KK,\i\S}(\RR)$ by Remark \ref{rem:intro:CategFctsProps}\ref{rem:intro:CategFctsProps:CKF:Schwartz&Bump} and Corollary \ref{cor:CKF:AsymExpKernel} (or when $\KK = \QQ$, by \cite[Section 7]{ClCoRoSe}).
\end{proof}
\end{enumerate}
\end{rems}

The hierarchy of categories \eqref{eq:Background:OurCategs} was first developed in the special case of $\KK=\QQ$ in Cluckers, Comte, Miller, Rolin, and Servi \cite{ClCoRoSe}, where \eqref{eq:Background:OurCategs} was instead denoted by $\S \subset \C \subset \C^{\exp}_{\textrm{naive}} \subset \C^{\exp}$ and $\C$-functions were called constructible functions (which were studied previously in Cluckers and Miller \cite{CluMil:2011,CluMil:2012,CluMil:2013,CluMil}).  These ideas were later extended to \eqref{eq:Background:OurCategs} for any complex subfield $\KK$ in Cluckers, Comte, Rolin, and Servi \cite{ClCoRoSe:preprint} (for $\CC^{\KK}$) and Cluckers, Comte, and Servi in \cite{clcose:fourier_mellin_power_constructible} (for $\C^{\KK,\FTs}$).\footnote{To be precise, \cite{clcose:fourier_mellin_power_constructible} only worked with $\C^{\CC,\FTs}$, but all of that paper's results that are relevant to the present work go through verbatim for $\C^{\KK,\FTs}$ for any subfield $\KK$ of $\CC$.}  All of these developments were, in turn, an outgrowth of the earlier work on preparation of subanalytic functions (in Parusi{\'n}ski \cite{Par:1994} and in Lion and Rolin \cite{LR:1997}) and on integration of subanalytic functions (in Lion and Rolin \cite{LionRol} and in Comte, Lion, and Rolin \cite{CoLiRo}).

\begin{rems}\label{rem:intro:CategFctsProps}
We now clarify and extend some remarks made in the Introduction.
\begin{enumerate}
\item\label{rem:intro:CategFctsProps:CK}
The category $\C^{\KK}$ is the smallest category of rings of function on the subanalytic sets containing $\S^{\KK}$ that is stable under parametric integration. (See \cite[Theorem 2.4]{ClCoRoSe:preprint}.)

\item\label{rem:intro:CategFctsProps:CKF}
The category $\C^{\KK,\FTs}$ is the smallest category of rings of function on the subanalytic sets containing $\S^{\KK}\cup\{\e^{\i f} : f\in\S(X)\}_{X\in\Omega}$ that is stable under parametric integration, and $\C^{\KK,\FTs}$ is also the smallest category of rings of function on the subanalytic sets containing $\S^{\KK}$ that is stable under parametric Fourier transforms.  (See \cite[Theorem 2.9]{clcose:fourier_mellin_power_constructible}.)

\item\label{rem:intro:CategFctsProps:CKF:Schwartz&Bump}
The function $f(t) = (1+t^2)^{-1}$ is in $\S(\RR)$, so $\FT[f](\xi) = \pi\e^{-2\pi|\xi|}$ is in $\C^{\QQ,\FTs}(\RR)$ by the previous remark, from which it follows that $\C^{\QQ,\FTs}(\RR)$ contains the nonzero Schwartz function $t\mapsto \e^{-t^2}$ and many nonzero $C^\infty$ functions with compact support, as shown in Aizenbud, Cluckers, Raibaut, and Servi \cite[Proposition 2.9]{AiClRaSe}.  This allows the notions of distributions and tempered distributions to be adapted to the category $\C^{\QQ,\FTs}$, as also done in \cite{AiClRaSe}, but we shall not make use of that here because our focus is on Fourier transforms of integrable functions, not distributions.
\end{enumerate}
\end{rems}

\subsection{Preparation}\label{ss:Background:Prep}

Our next task is to formulate a version of the subanalytic preparation theorem and a related preparation theorem for $\C^{\KK}$.  Using the variables $(x,y) = (x_1,\ldots,x_m,y)$ on $\RR^{m+1}$, the subanalytic preparation theorem from \cite{LR:1997} considers a finite set $\F$ of subanalytic functions on a subanalytic set $D\subseteq\RR^{m+1}$ and constructs a finite subanalytic cell decomposition $\A$ of $D$ over $\RR^m$ such that for each function $f\in\F$ and each cell $A\in\A$ over $\RR^m$ with $1$-dimensional fibers over $\pi(A)$, the function $f$ may be expressed in the form
\begin{equation}\label{eq:fprep}
f(x,y) = c(x)|y-\theta(x)|^r u(x,y)
\end{equation}
on $A$, where $c,\theta\in\S(\pi(A))$, $r\in\QQ$, and $u$ is a subanalytic unit on $A$ of a special form.  Our use here of the phrase ``over $\RR^m$'' requires some explanation: calling $A$ a cell over $\RR^m$ means that $A$ is required to have the structure of a cell in the $y$-variable but not necessarily in the $x$-variables, so that $\pi(A)$ is only required to be a subanalytic set (not necessarily a cell); and accordingly, calling $\A$ a cell decomposition of $D$ over $\RR^m$ means that $\A$ is a partition of $D$ into cells over $\RR^m$ such that $\{\pi(A) : A\in\A\}$ is a partition of $\pi(D)$.

We have no use in this paper for the members of $\A$ with $0$-dimensional fibers over $\pi(A)$, so it is convenient to formulate our notion of a decomposition of $D$ to only include the members of $\A$ with $1$-dimensional fibers over $\pi(A)$.  And, the data used to specify the prepared form of $f$ on $A$ given in \eqref{eq:fprep}, including the special form of the unit $u$, depends only on the choice of $A\in\A$ and not on the choice of $f\in\F$.  So it is also convenient to consider the set $A$ equipped with its associated preparation data to be the basic object in our decompositions, not just the set $A$ alone.  We now codify these ideas in some technical terminology given in Definitions \ref{def:PrepData}-\ref{def:Decomp}.

\begin{defs}\label{def:PrepData}
We use the phrase \DEF{prepared coordinate data over $\RR^m$} to mean a tuple $(A,\theta,d,\psi)$ such that $A$ is a subanalytic subset of $\RR^{m+1}$, the function $\theta\in\S(\pi(A))$ is analytic, $d$ is a positive integer, and $\psi:A\to[0,1]^{N+2}$ is an analytic subanalytic map, for some $N\in\NN$, with the following three properties.
\begin{enumerate}[label = (P\arabic*), ref = (P\arabic*)]
\item\label{def:PrepData:Prop1}
There exist $\sigma\in\{+,-\}$, $\tau\in\{-1,1\}$, and analytic functions $a,b\in\S(\pi(A))$ such that
\[
A = \{(x,y)\in \pi(A)\times\RR : a(x) < y_{\ST} < b(x)\},
\]
where $a(x) < b(x)$ on $\pi(A)$ and
\begin{equation}\label{eq:PrepData:yS}
y_{\ST} = \sigma(y-\theta(x))^{\tau},
\end{equation}
and where either $a(x) = 0$ on $\pi(A)$, or $a(x) > 0$ on $\pi(A)$ and $\tau = 1$.

\item\label{def:PrepData:Prop2}
The range of $\theta$ is contained in $(-\infty,0)$, $\{0\}$, or $(0,+\infty)$; and when $\theta$ is nonzero, there exists a $\lambda > 1$ such that $\lambda^{-1} < \frac{y}{\theta(x)} < \lambda$ on $A$.

\item\label{def:PrepData:Prop3}
The map $\psi$ is of the form
\[
\qquad
\psi(x,y) = \left(\psi_1(x),\ldots,\psi_N(x), \left(\frac{a(x)}{y_{\ST}}\right)^{1/d}, \left(\frac{y_{\ST}}{b(x)}\right)^{1/d} \right)
\]
on $A$ for some $\psi_1,\ldots,\psi_N\in\S(\pi(A))$.
\end{enumerate}
We respectively call $A$, $\theta$, and $\pi(A)$ the \DEF{domain}, \DEF{center}, and \DEF{base set} of $(A,\theta,d,\psi)$; we say that $(A,\theta,d,\psi)$ \DEF{lies above $\pi(A)$}; and we say that $(A,\theta,d,\psi)$ is \DEF{contained} in a set $D$ when $A\subseteq D$.
The conditions on $a$ and $\tau$ divide into the following cases.
\begin{enumerate}[leftmargin = *]
\item[]
Case 1: $a > 0$ and $\tau = 1$.
\begin{quoting}[leftmargin = 2ex]
In this case, for each $x\in\pi(A)$, $\theta(x)$ is not in the closure of the fiber of $A_x$ in $\RR$.
\end{quoting}

\item[]
Case 2.A: $a = 0$ and $\tau = 1$.
\begin{quoting}[leftmargin = 2ex]
In this case, for each $x\in\pi(A)$, $\theta(x)$ is in the closure of the fiber $A_x$ in $\RR$.
\end{quoting}

\item[]
Case 2.B: $a = 0$ and $\tau = -1$.
\begin{quoting}[leftmargin = 2ex]
In this case, $\theta = 0$ by Property \ref{def:PrepData:Prop2}, and for each $x\in\Pi(A)$, $\sigma\infty$ is in the closure of the fiber $A_x$ in the extended real line.
\end{quoting}
\end{enumerate}
We define the \DEF{limiting boundary of $(A,\theta,d,\psi)$} to be $\theta^\sigma$ in Case 2.A and to be $\sigma\infty$ in Case 2.B.  We do not define a limiting boundary in Case 1.  (See Remark \ref{rem:PrepData}\ref{rem:PrepData:LimBnd} below for an explanation of this choice of terminology.)  For $y_{\ST}$ defined in \eqref{eq:PrepData:yS}, we call $(x,y_{\ST})$ the \DEF{standard coordinates for $(A,\theta,d,\psi)$}.   We also write
\begin{equation}\label{eq:PrepData:yA}
y_{\AF} = \sigma(y-\theta(x))
\end{equation}
and call $(x,y_{\AF})$ the \DEF{affine coordinates for $(A,\theta,d,\psi)$}.
\end{defs}

\begin{defs}\label{def:psiFct}
Let $(A,\theta,d,\psi)$ be as in the Definitions \ref{def:PrepData}.  If $F:[0,1]^{N+2}\to\CC$ extends to a holomorphic function $F:[0,1]^N\times D_1(0)^2\to\CC$, then we call $F\circ\psi$ a \DEF{$\psi$-function}.  If $\epsilon \in (0,1)$ and $U:[0,1]^{N+2} \to [1-\epsilon,1+\epsilon]$ extends to a holomorphic function $U:[0,1]^N \times D_1(0)^2 \to D_{\epsilon}(1)$, then we call $U\circ\psi$ a \DEF{$(\psi,\epsilon)$-unit}. A \DEF{$\psi$-unit} is a $(\psi,\epsilon)$-unit for some $\epsilon\in(0,1)$.
\end{defs}

\begin{rems}\label{rem:PrepData}
We now remark on some intended uses and variations of the notions from Definitions \ref{def:PrepData} and \ref{def:psiFct}.  As such, consider again $(A,\theta,d,\psi)$ as in the Definitions \ref{def:PrepData}.
\begin{enumerate}
\item\label{rem:PrepData:LimBnd}
We shall frequently be in the situation where $(A,\theta,d,\psi)$ is either in Case 2.A or 2.B, we have a function $f:A\to\CC$ that is of a certain prepared form determined from the data $(A,\theta,d,\psi)$, and for each $x\in\pi(A)$, we wish to study the asymptotic behavior of $f_x$ as $y\to\theta(x)^{\sigma}$ in Case 2.A or as $y\to\sigma\infty$ in Case 2.B.  Thus, depending on the sign $\sigma$ and on the case 2.A or 2.B, the object $\theta^\sigma$ or $\sigma\infty$ specifies either the upper or lower boundary points of each of the fibers in $\{A_x:x\in\pi(A)\}$ in the extended real line to which we let $y$ approach in our limits, along with the direction of approach.  This is why we call $\theta^\sigma$ or $\sigma\infty$ the ``limiting boundary'' of $(A,\theta,d,\psi)$.  Both Cases 2.A and 2.B are equivalent to saying that we study $f_x$ as $y_{\ST}\to 0^+$, and that is why we call $(x,y_{\ST})$ the standard coordinates for $(A,\theta,d,\psi)$, for these coordinates allow us to ``standardize'' our notation by phrasing all of our asymptotics at $0^+$.  We do not define a limiting boundary in Case 1 because $y_{\ST}$ cannot approach $0^+$ in the fibers $\{A_x : x\in\pi(A)\}$ in that case.

\item\label{rem:PrepData:Case2}
We are interested in working with $\psi$-functions, not really with the map $\psi$ in and of itself.  In both Cases 2.A and 2.B, a $\psi$-function does not depend on the component $(a(x)/y_{\ST})^{1/d}$ of $\psi$ because this component is identically zero.  To emphasize this fact, we shall instead write
\[
\psi(x,y) = \left(\psi_1(x),\ldots,\psi_N(x), \left(\frac{y_{\ST}}{b(x)}\right)^{1/d} \right)
\]
on $A$ in situations when we are only considering the Cases 2.A or 2.B.

\item\label{rem:PrepData:yA}
We have $y_{\AF} = y_{\ST}$ in Cases 1 and 2.A, and we have $\theta = 0$ and $y_{\AF} = y_{\ST}^{-1} = \sigma y$ in Case 2.B.  Thus, the only distinction between $(x,y_{\ST})$ and $(x,y_{\AF})$ occurs in Case 2.B.  The fact that $y_{\ST}$ is not an affine function of $y$ in Case 2.B is inconvenient in certain applications, in which case we will opt to use the affine coordinates $(x,y_{\AF})$ instead.  In situations when we are only considering the Case 2.B and are working with the coordinates $(x,y_{\AF})$, we shall modify our notation for $A$ and $\psi$ accordingly by writing
\[
A = \{(x,y) \in \pi(A)\times\RR : y_{\AF} > b(x)\}
\]
and
\[
\psi(x,y) = \left(\psi_1(x),\ldots,\psi_N(x), \left(\frac{b(x)}{y_{\AF}}\right)^{1/d} \right)
\]
on $A$, where $b\in\S_+(\pi(A))$.

\item\label{rem:PrepData:LimBndCoord}
Consider the Case 2.A or 2.B, and for the moment write $\ell$ for the limiting boundary of $(A,\theta,d,\psi)$.  (Thus, $\ell = \theta^\sigma$ or $\ell = \sigma\infty$.)  The formula \eqref{eq:PrepData:yS} for $y_{\ST}$ is determined by $\ell$ alone, so we may also refer to $(x,y_{\ST})$ as the \DEF{standard coordinates for $\ell$}.  In the same way, we may also refer to $(x,y_{\AF})$ as the \DEF{affine coordinates for $\ell$}.  In the natural way, we may use the coordinates $(x,y_{\ST})$ and $(x,y_{\AF})$ on the set $\{(x,y)\in\pi(A)\times\RR: y_{\ST} > 0\}$, not just on the subset $A$.

\item\label{rem:PrepData:LimBndCoordBIS}
In our forthcoming discussion of preparation theorems and asymptotic expansions, when $\D\in\{\S^{\KK},\C^{\KK},\C^{\KK,\i\S},\C^{\KK,\FTs}\}$ and a $\psi$-function $f = F\circ\psi$ on $A$ appears in either a prepared form or an asymptotic expansion for a function in $\D(A)$, one should assume that $f(A)\subseteq\cl_{\CC}(\KK)$.
\end{enumerate}
\end{rems}

\begin{defs}\label{def:Decomp}
Consider a subanalytic set $D\subseteq\RR^{m+1}$, and consider a finite set $\A$ of prepared coordinate data over $\RR^m$ contained in $D$ with disjoint domains, as defined in the Definitions \ref{def:PrepData}.  We say that $\A$ is \DEF{compatible} with a set $\X$ of subsets of $D$ when for each $X\in\X$ and $(A,\theta,d,\psi)\in\A$, either $A\subseteq X$ or $A\cap X = \emptyset$.  We say that $\A$ \DEF{partitions $D$ over $\pi(D)$} to mean that for each $x\in\pi(D)$, the fiber $(D\setminus\bigcup_{(A,\theta,d,\psi)\in\A} A)_x$ is a finite subset of $\RR$.  If $\A$ partitions $D$ over $\pi(D)$ and is such that its collection $\{\pi(A) : (A,\theta,d,\psi)\in\A\}$ of base sets is a partition of $\pi(D)$, then we call $\A$ a \DEF{decomposition of $D$ over $\pi(D)$}.
\end{defs}

The next proposition is a formulation of the subanalytic preparation theorem from \cite{LR:1997}.

\begin{prop}[$\S$-Preparation]\label{prop:SubPrep}
Let $D\subseteq\RR^{m+1}$ be subanalytic, $\F\subseteq\S(D)$ be finite, and $\epsilon\in(0,1)$.  Then there exists a finite set $\A$ of prepared coordinate data partitioning $D$ over $\pi(D)$ such that for each $(A,\theta,d,\psi)\in\A$ and $f\in\F$,
\begin{equation}\label{eq:SubPrep}
f(x,y) = c(x)y_{\ST}^{r} u(x,y)
\end{equation}
on $A$ for some analytic function $c\in\S(\pi(A))$, exponent $r\in\QQ$ with $dr\in\ZZ$, and $(\psi,\epsilon)$-unit $u$.
\end{prop}

We say that we have \DEF{$\epsilon$-prepared} $f$ when we apply Proposition \ref{prop:SubPrep}, and we call \eqref{eq:SubPrep} the \DEF{$\epsilon$-prepared form} of $f$ on $A$.  Usually the choice of $\epsilon$ is unimportant, in which case we will just say \emph{prepared} rather than \emph{$\epsilon$-prepared}.

\begin{prop}[$\S$-Rectilinearization]\label{prop:SubRect}
Let $D\subseteq\RR^m$ be subanalytic and $\F\subseteq\S(D)$ be finite.  Then there exists a finite set $\Phi$ of subanalytic maps with the following properties.
\begin{enumerate}
\item
The ranges of the maps in $\Phi$ partition $D$.

\item
Each $\varphi\in\Phi$ is a bi-analytic map from $(0,1)^{m_\varphi}$  onto its range for some $m_\varphi\in\{0,\ldots,m\}$.

\item
For each $\varphi\in\Phi$ and $f\in\F$, we have
\[
f\circ\varphi(x) = x^\alpha u(x) \,\,\text{on $(0,1)^{m_\varphi}$}
\]
for some $\alpha\in\ZZ^{m_\varphi}$ and analytic unit $u:D_1(0)^{m_\varphi}\to D_{1/2}(1)$ with $u((0,1)^{m_\varphi}) \subseteq \RR$.
\end{enumerate}
\end{prop}

Proposition \ref{prop:SubRect} can be proven through repeated application of Proposition \ref{prop:SubPrep}.  For a proof, see \cite[Theorem 1.5]{CluMil:2013} for a more detailed parametric version of Proposition \ref{prop:SubRect}, and see \cite[Proposition 6.4]{AiClRaSe} for a proof of a nonparametric version from the parametric version.

The next proposition is a formulation of the $\C^{\KK}$-preparation theorem of \cite[Proposition 4.2]{ClCoRoSe:preprint}.  We opt to include a proof so that we may then adapt it to prove Lemma \ref{lem:CKiSprep} below, which is not explicitly stated in the literature.

\begin{prop}[$\C^{\KK}$-Preparation]\label{prop:CKprep}
Let $D\subseteq\RR^{m+1}$ be subanalytic and $\F\subseteq\C^{\KK}(D)$ be finite.  Then there exists a finite set $\A$ of prepared coordinate data partitioning $D$ over $\pi(D)$ such that for each $(A,\theta,d,\psi)\in\A$, each $f\in\F$ may be expressed as a finite sum
\begin{equation}\label{eq:CKprep}
f(x,y) = \sum_{i\in I} c_i(x) y_{\ST}^{r_i}(\log y_{\ST})^{s_i} f_i(x,y)
\end{equation}
on $A$, where $I$ is a finite index set and for each $i\in I$, the function $c_i\in\C^{\KK}(\pi(A))$ is analytic, $r_i\in\KK$, $s_i\in\NN$, and $f_i$ is a $\psi$-function.  In addition, when $f\in\F\cap\S^{\KK}(D)$, we have $c_i\in\S^{\KK}(\pi(A))$ and $s_i = 0$ for all $i\in I$.
\end{prop}

\begin{proof}
Each $f\in\F$ is a finite sum of terms of the form $w f_{0}^{\gamma} \prod_{i=1}^{k}\log f_i$, where $w\in\cl_{\CC}(\KK)$, $\gamma\in\KK$, and $f_0,\ldots,f_k\in\S_+(X)$.  Let $\Gamma$ be the set of all such $\gamma$ for all terms of members of $\F$, and fix an $\epsilon\in(0,1)$ such that $\{z^\gamma : z\in D_\epsilon(1)\} \subseteq D_{1/2}(1)$ for all $\gamma\in\Gamma$. Let $\A$ be the set of prepared coordinate data over $\RR^m$ given by applying Proposition \ref{prop:SubPrep} to $\epsilon$-prepare the set of all such $f_0,\ldots,f_k$ for all terms of members of $\F$.  Consider any $(A,\theta,d,\psi)\in\A$, and write
\[
f_i(x,y) = c_i(x) y_{\ST}^{r_i}u_i(x,,y), \quad\text{for each $i\in\{0,\ldots,k\}$,}
\]
for the prepared forms for all such functions $f_0,\ldots,f_k$.  Now distribute the powers and expand the logarithms, by which we mean that we write
\[
f_0(x,y)^\gamma = c_0(x)^\gamma y_{\ST}^{r_0\gamma} u_0(x,y)^{\gamma}
\]
and
\[
\log f_i(x,y) = \log c_i(x) + r_i\log y_{\ST} + \log u_i(x,y)
\]
for each $i\in\{1,\ldots,k\}$.  Then $u_{0}^{\gamma}$ and each $\log u_i$ is a $\psi$-function, so distributing multiplication across addition expresses each member of $\F$ in the form \eqref{eq:CKprep}.  The last sentence of the proposition follows from the fact that no logarithms are present when $f\in\F\cap\S^{\KK}(D)$.
\end{proof}

In Propositions \ref{prop:SubPrep} and \ref{prop:CKprep}, if we are additionally given a finite set $\X$ of subanalytic subsets of $D$, then we may construct $\A$ to be a decomposition of $D$ over $\pi(D)$ that is compatible with $\X$.  Indeed, just apply the proposition, either \ref{prop:SubPrep} or \ref{prop:CKprep}, to the restrictions of the functions in $\F$ to each atom of the boolean algebra of subsets of $D$ generated by $\X$ (which guarantees compatibility with $\X$), and then further partition the base sets of the members of $\A$ by the atoms of the boolean algebra of subsets of $\pi(D)$ that they generate (which creates a decomposition, rather than just a partition, of $D$ over $\Pi(D)$).

\subsection{Asymptotic Expansions}\label{ss:Background:AsymExp}

Our next task is to specify the forms of the asymptotic expansions for functions in each of the categories of the hierarchy $\S^{\KK}\subset \C^{\KK} \subset \C^{\KK,\i\S} \subset \C^{\KK,\FTs}$.  We quote \cite{clcose:fourier_mellin_power_constructible} to obtain the expansions for $\C^{\KK,\FTs}$-functions, but we supply proofs for the other categories since those expansions are easily obtained and are only implicit in the work \cite{clcose:fourier_mellin_power_constructible}.  To simplify the exposition, we focus on asymptotic behavior at $0^+$, for this can be easily adapted to give one-sided asymptotic expansions at any other point of the extended real line, as explained in the Remarks \ref{rem:AsymExp}.  The next definition specifies the form of the phase functions that arise in terms of the asymptotic expansions for functions in $\C^{\KK,\i\S}$ and $\C^{\KK,\FTs}$.

\begin{defn}\label{def:PrincLaurentPoly}
Consider $(A,\theta,d,\psi)$ in Case 2.A or 2.B of the Definitions \ref{def:PrepData}, and employ the notation from Remark \ref{rem:PrepData}\ref{rem:PrepData:Case2}.  If the function $c\in\S(B)$ is analytic and $P:[0,1]^N\times(\CC\setminus\{0\})\to\CC$ is of the form
\[
P(w_1,\ldots,w_N,z) = \sum_{j=1}^{\mu}P_j(w_1,\ldots,w_N) z^{-j}
\]
for some $\mu\in\NN$ and analytic functions $P_1,\ldots,P_\mu:[0,1]^N\to\CC$, then we call the function
\[
(x,y)\mapsto c(x) P\circ\psi(x,y) = c(x)\sum_{j=1}^{\mu} P_j(\psi_j(x),\ldots,\psi_N(x))\left(\frac{y_{\ST}}{b(x)}\right)^{-j/d}
\]
on $A$ a \DEF{principal Laurent $\psi$-polynomial}.
\end{defn}

\begin{lem}[$\C^{\KK,\i\S}$-Preparation at $0^+$]\label{lem:CKiSprep}
Consider a finite set $\F\subseteq\C^{\KK,\i\S}(X\times\RR_+)$, where $X\subseteq\RR^m$ is subanalytic.  There exists a finite partition $\B$ of $X$ into subanalytic sets such that for each $B\in\B$, there exist prepared coordinate data $(A,0,d,\psi)$ lying over $B$ with limiting boundary $0^+$ such that each $f\in\F$ may be written in the prepared form
\begin{equation}\label{eq:CKiSprep}
f(x,y) = \sum_{i\in I} c_i(x) y^{r_i}(\log y)^{s_i} f_i(x,y) \e^{\i \phi_i(x,y)}
\end{equation}
on $A$, where $I$ is a finite index set and for each $i\in I$, we have $c_i\in\C^{\KK,\i\S}(B)$, $r_i\in\KK$, $s_i\in\NN$, $f_i$ is a $\psi$-function, and $\phi_i$ is a principal Laurent $\psi$-polynomial.  In addition:
\begin{enumerate}
\item\label{lem:CKiSprep:phase}
For all $i,j\in I$, either $\phi_i=\phi_j$ or else $(\phi_i)_x\neq(\phi_j)_x$ for all $x\in B$.

\item\label{lem:CKiSprep:CK}
For each $\D\in\{\S^{\KK},\C^{\KK}\}$, $f\in\F\cap\D(X\times\RR_+)$, and $i\in I$, we have $c_i\in\D(B)$ and $\phi_i = 0$, and also $s_i = 0$ when $\D=\S^{\KK}$.
\end{enumerate}
\end{lem}

\begin{proof}
Each $f\in\F$ is a finite sum of terms of the form
\begin{equation}\label{eq:CKiSprep:Terms}
w f_{0}^{\gamma}\left(\prod_{i=1}^{k}\log f_i\right) \e^{\i\phi},
\end{equation}
where $w\in\FF$, $\gamma\in\KK$, $f_0,\ldots,f_k\in\S_+(X\times\RR_+)$, and $\phi\in\S(X\times\RR_+)$.  Applying Proposition \ref{prop:SubPrep} to the set of all such phase functions $\phi$ produces a finite set of prepared coordinate data decomposing $X\times\RR_+$ over $X$.  Let $\A'$ be the members of this decomposition with limiting boundary $0^+$.  Consider any $(A,\theta,d,\psi)\in\A'$, and write $A = \{(x,y)\in \pi(A)\times\RR : 0 < y < b(x)\}$.  Using the notation from Remark \ref{rem:PrepData}\ref{rem:PrepData:Case2}, for each prepared phase function $\phi$, we may expand its unit as a power series in $(\frac{y}{b(x)})^{1/d}$ to write
\begin{equation}\label{eq:CKisprep:LaurentForm}
\phi(x,y) = \phi_-(x,y) + \phi_0(x) + \left(\frac{y_{\ST}}{b(x)}\right)^{1/d}\phi_+(x,y)
\end{equation}
on $A$, where $\phi_-$ is a principal Laurent $\psi$-polynomial, $\phi_0\in\S(\pi(A))$, and $\phi_+$ is a $\psi$-function.  We may replace $b$ with a suitable $b'\in\S_+(\pi(A))$ with $b'\leq b$ (thereby replacing $A$ with a smaller cell) to assume that the third term in the right side of \eqref{eq:CKisprep:LaurentForm} is bounded.  So by writing
\[
\e^{\i\phi(x,y)} = \e^{\i\phi_-(x,y)} \e^{\i\phi_0(x)} \e^{\i\left(\frac{y_{\ST}}{b(x)}\right)^{1/d}\phi_+(x,y)}
\]
and using the fact that the real and imaginary parts of the third exponential factor in the right side of this equation are subanalytic, we reduce to the case where each function in $\F$ is a finite sum of terms of the form \eqref{eq:CKiSprep:Terms} on $A$, except now $w$ is a function in $\C^{\KK,\i\S}(\pi(A))$ and the phase function $\phi$ is a principal Laurent $\psi$-polynomial.  To finish, partition the base sets of the members of $\A'$ to achieve the property \ref{lem:CKiSprep:phase}, and then for each $(A,0,d,\psi)\in\A'$, proceed as in the proof of Proposition \ref{prop:CKprep} to prepare all such functions $f_0,\ldots,f_k$ on $A$. Let $\A$ be the members of the resulting decompositions of $A$ over $\pi(A)$, for each $(A,0,d,\psi)\in\A'$, that have limiting boundary $0^+$.  Then $\B := \{\pi(A) : (A,0,d,\psi)\in\A\}$ is the desired partition of $X$.  Indeed, for each $B\in\B$, there is exactly one  $(A,0,d,\psi)\in \A$ lying over $B$, and each $f\in\F$ is of the form \eqref{eq:CKiSprep} on $A$.  Property \ref{lem:CKiSprep:CK} follows from the fact that $\phi = 0$ for each term \eqref{eq:CKiSprep:Terms} of $f\in\F$ when $f\in\C^{\KK}(X\times\RR_+)$ and that also $k=0$ in \eqref{eq:CKiSprep:Terms} when $f\in\S^{\KK}(X\times\RR_+)$.
\end{proof}

\begin{defs}\label{def:AsympExp}
Consider prepared coordinate data $(A,0,d,\psi)$ with limiting boundary $0^+$, write $A = \{(x,y)\in B\times\RR : 0 < y < b(x)\}$, and use the notation from Remark \ref{rem:PrepData}\ref{rem:PrepData:Case2}.  In general, we call a set $\mathbf{AS}$ of complex-valued functions on $A$ an \DEF{asymptotic scale at $0^+$} when $\mathbf{AS}$ is linearly ordered by the relation $\prec$ defined by setting, for each $f,g\in\mathbf{AS}$, $f\prec g$ if and only if $f_x(y) = o(g_x(y))$ as $y\to 0^+$ for all $x\in B$.  If $\CF$ is a set of bounded complex-valued functions on $A$, then an \DEF{asymptotic expansion in $\AS$ with coefficient functions in $\CF$} is a series of the form $\sum_{i\in\N} c_is_i$, where $\N$ is $\NN$ or an initial segment of $\NN$, $\{c_i\}_{i\in \N}\subseteq\CF$, and $\{s_i\}_{i\in \N} \subseteq\AS$ is such that $s_0 \succ s_1 \succ s_2 \succ \cdots $.

To introduce the specific asymptotic expansions we work with in this paper, let $\Gamma$ be a finite subset of $\KK$, let $S$ be a finite initial segment of $\NN$, and let $\Phi$ be a finite set of principal Laurent $\psi$-polynomials satisfying the following distinctness condition.
\begin{equation}\label{eq:AsymExp:LaurPolyCond}
\text{For all $\phi_1,\phi_2\in\Phi$, either $\phi_1=\phi_2$ or $(\phi_1)_x\neq(\phi_2)_x$ for all $x\in B$.}
\end{equation}
The data $(A,0,d,\psi)$, $\Gamma$, and $S$ determine the asymptotic scale
\begin{equation}\label{eq:AsymExp:Scale}
\textstyle \{y^{r}(\log y)^{s} : (r,s)\in (\RE(\Gamma) + \frac{1}{d}\NN)\times S\}
\end{equation}
on $A$ at $0^+$.  Since $\Gamma$ and $S$ are finite, the linear order
\[
\textstyle (\{y^{r}(\log y)^{s} : (r,s)\in (\RE(\Gamma) + \frac{1}{d}\NN)\times S\},\succ)
\]
is isomorphic to $(\NN,<)$.  (Notice the reversal of the directions of $\succ$ and $<$.)

The type of coefficient functions we work with depends on the choice of category $\D\in\{\S^{\KK},\C^{\KK},\C^{\KK,\i\S},\C^{\KK,\FTs}\}$. Define $\E[\D](\Gamma,\Phi)$ to be the set of all functions $E:A\to\CC$ of the form
\begin{equation}\label{eq:AsymExp:EKF}
E(x,y) = \sum_{(\gamma,\phi)\in\IM(\Gamma)\times\Phi} c_{\gamma,\phi}(x) y^{\i\gamma} \e^{\i\phi(x,y)}
\end{equation}
with $c_{\gamma,\phi}\in\D(B)$ for all $(\gamma,\phi)\in\IM(\Gamma)\times\Phi$, with the stipulation that $\Phi = \{0\}$ when $\D\in\{\S^{\KK},\C^{\KK}\}$.
\end{defs}

The next lemma is proven in \cite[Proposition 3.4]{ClCoRoSe:preprint}\footnote{However, Lemma \ref{lem:noncomp} corrects the statement of \cite[Proposition 3.4]{ClCoRoSe:preprint} in three ways.  First, the formulation of Lemma \ref{lem:noncomp}\ref{lem:noncomp:Seq} in \cite[Proposition 3.4]{ClCoRoSe:preprint} unnecessarily assumes that $(\gamma_j,\phi_j)\neq (0,0)$.  Second, the formulation of Lemma \ref{lem:noncomp}\ref{lem:noncomp:TwoSeq} in \cite[Proposition 3.4]{ClCoRoSe:preprint} neglects to mention the required hypothesis that $c_j\neq 0$ and $(\gamma_j,\phi_j)\neq (0,0)$ for some $j\in J$.  And third, the formulation of Lemma \ref{lem:noncomp}\ref{lem:noncomp:L1} in \cite[Proposition 3.4]{ClCoRoSe:preprint} does not include the function $g$, but its proof goes through with $g$ present, and $g$ is needed in most applications of \ref{lem:noncomp:L1}.} using the theory of continuously uniformly distributed functions modulo $1$.  Loosely speaking, the lemma shows that when $E$ is of the form \eqref{eq:AsymExp:EKF},  if $x\in B$ is such that $c_{\gamma,\phi}(x)\neq 0$ for some $(\gamma,\phi)\in\IM(\Gamma)\times\Phi$, then $E_x$ is a nonzero function whose zero set is not ``too large'' near $0^+$.

\begin{lem}[Noncompensation]\label{lem:noncomp}
Define
\[
E(y) = \sum_{j\in J}c_j y^{\i\gamma_j} \e^{\i\phi_j(y^{-1/d})}
\]
on $\RR_+$ for some finite index set $J$, family $\{c_j\}_{j\in J}$ of complex coefficients, family $\{(\gamma_j,\phi_j)\}_{j\in J}$ of distinct members of $\RR\times\{\phi\in\RR[y] : \phi(0) = 0\}$, and positive integer $d$.
\begin{enumerate}
\item\label{lem:noncomp:Seq}
If $c_j\neq 0$ for some $j\in J$, then there exist an $\epsilon > 0$ and a sequence $\{s_i\}_{i\in\NN}$ in $\RR_+$ converging to $0$ such that $|E(s_i)|\geq \epsilon$ for all $i\in\NN$.

\item\label{lem:noncomp:TwoSeq}
If $c_j\neq 0$ and $(\gamma_j,\phi_j)\neq (0,0)$ for some $j\in J$, then there exist an $\epsilon > 0$ and sequences $\{s_i\}_{i\in\NN}$ and $\{t_i\}_{i\in\NN}$ in $\RR_+$ converging to $0$ such that $|E(s_i) - E(t_i)|\geq \epsilon$ for all $i\in\NN$.

\item\label{lem:noncomp:L1}
Let $b > 0$, and define $f(y) = y^r(\log y)^s E(y)+ g(y)$ on $(0,b)$ for some $r\leq -1$, $s\in\NN$, and function $g:(0,b)\to\RR$ with $g(y) = o(y^r(\log y)^s)$ as $t\to 0^+$.  If $f\in \L^1((0,b))$, then $c_j = 0$ for all $j\in J$.
\end{enumerate}
\end{lem}

\begin{prop}[Asymptotic Expansions]\label{prop:AsymExp}
Let $\D\in\{\S^{\KK},\C^{\KK},\C^{\KK,\i\S},\C^{\KK,\FTs}\}$, and consider any finite set $\F\subseteq\D(X\times\RR_+)$, where $X\subseteq\RR^m$ is subanalytic.  There exists a finite partition $\B$ of $X$ into subanalytic sets such that for each $B\in\B$, there exist prepared coordinate data $(A,0,d,\psi)$ lying over $B$ with limiting boundary $0^+$, a finite set $\Gamma\subseteq\KK$, a finite initial segment $S$ of $\NN$, and a finite set $\Phi$ of principal Laurent $\psi$-polynomials satisfying the distinctness condition \eqref{eq:AsymExp:LaurPolyCond} such that for each $f\in\F$, there exists a unique asymptotic expansion
\begin{equation}\label{eq:AsymExp}
\sum_{i\in \N} y^{r_i}(\log y)^{s_i} E_i(x,y)
\end{equation}
in the scale \eqref{eq:AsymExp:Scale} with coefficient functions $\{E_i\}_{i\in\N}$ in $\E[\D](\Gamma,\Phi)$ such that for each $\rho\in\RR$, there exist $n(\rho)\in \N$ and a function $C_\rho:B\to\RR_+$ such that
\[
\left|f(x,y) - \sum_{i=0}^{n(\rho)} y^{r_i}(\log y)^{s_i} E_i(x,y)\right| \leq C_\rho(x) y^{\rho}
\]
on $A$, with the understanding that $\Phi = \{0\}$ when $\D\in\{\S^{\KK},\C^{\KK}\}$ and that $S = \{0\}$ when $\D=\S^{\KK}$.  In addition, the expansion \eqref{eq:AsymExp} converges to $f$ on $A$ when $\D\in\{\S^{\KK},\C^{\KK},\C^{\KK,\i\S}\}$.
\end{prop}

\begin{proof}
The uniqueness assertion follows from Lemma \ref{lem:noncomp}\ref{lem:noncomp:Seq} and the proof of \cite[Lemma 7.2]{ClCoRoSe}, and the existence assertion for $\D = \C^{\KK,\FTs}$ is proven in \cite[Theorem 7.6]{clcose:fourier_mellin_power_constructible}.  So we need only prove the existence assertion when $\D\in\{\S^{\KK},\C^{\KK},\C^{\KK,\i\S}\}$.  Apply Lemma \ref{lem:CKiSprep}, and let $B\in\B$ and $(A,0,d,\psi)$ be as in the statement of the lemma.  Let $f\in\F$, expand each $\psi$-function in the prepared form \eqref{eq:CKiSprep} of $f$ on $A$ as a power series in $(\frac{y_{\ST}}{b(x)})^{1/d}$ to write $f$ as the sum of a convergent series
\begin{equation}\label{eq:CKiS:preAsympExp}
f(x,y) = \sum_{(r,s,\phi)\in R\times S\times\Phi} c_{r,s,\phi}(x) y_{\ST}^{r}(\log y_{\ST})^s \e^{\i\phi(x,y)}
\end{equation}
on $A$, where $R = \Gamma + \frac{1}{d}\NN$ and $\{c_{(r,s,\phi)}\}_{(r,s,\phi)\in R\times S\times\Phi} \subseteq\D(B)$, with $\Phi = \{0\}$ when $\D\in\{\S^{\KK},\C^{\KK}\}$ and with $S = \{0\}$ when $\D = \S^{\KK}$.  Writing $y_{\ST}^{r} = y_{\ST}^{\RE(r)} y_{\ST}^{\i\IM(r)}$ and combining like terms in the power-log monomials $y_{\ST}^{\RE(r)}(\log y_{\ST})^s$ gives \eqref{eq:AsymExp}.
\end{proof}

\begin{rem}\label{rem:CKiS:FiniteAsymExp}
Consider the situation of Proposition \ref{prop:AsymExp} with $\D\in\{\S^{\KK},\C^{\KK},\C^{\KK,\i\S}\}$, and let $f\in\F$ and $\rho\in\RR$.  Instead of fully expanding the $\psi$-functions in \eqref{eq:CKiSprep} to obtain \eqref{eq:CKiS:preAsympExp}, we may only partially expand \eqref{eq:CKiSprep} to write $f$ in the form $f = f_\rho + R_\rho$ on $A$, where $f_\rho$ is a finite sum of terms from the convergent expansion \eqref{eq:AsymExp} and the remainder $R_\rho$ is a finite sum of the form \eqref{eq:CKiSprep} with $\RE(r_i)\geq\rho$ for all $i\in I$.
\end{rem}

\begin{defs}\label{def:AsymExp:CKF}
Let $f\in\C^{\KK,\FTs}(\RR_+)$, and apply Proposition \ref{prop:AsymExp} to $\{f\}$, where $m = 0$ (so $x$ is the empty tuple of variables).  Call the resulting series given in \eqref{eq:AsymExp} the \DEF{asymptotic expansion for $f$ at $0^+$}.  Suppose this expansion is nonzero, and let $i\in\N$ be minimal such that $E_i\neq 0$.  Respectively call
\begin{equation}\label{eq:AsympExp:data}
\left\{
\begin{aligned}
\LE(f) &= (r_i,s_i), \\
\LC(f) &= E_i(y), \\
\LM(f) &= y^{r_i}(\log y)^{s_i}, \\
\LF(f) &= \LC(f)\LM(f)
\end{aligned}\right.
\end{equation}
the \DEF{leading exponents}, \DEF{leading coefficient function}, \DEF{leading monomial}, and \DEF{leading form of $f$ at $0^+$}.  We also write $\LE(f) = (\LE_1(f),\LE_2(f))$ and respectively call
\begin{equation}\label{eq:AsympExp:LE}
\LE_{1}(f) = r_i \quad\text{and}\quad \LE_{2}(f) = s_i
\end{equation}
the \DEF{leading exponent of $y$ for $f$} and the \DEF{leading exponent of $\log y$ for $f$}.  In general, $\LC(f)$ is a finite sum of the form \eqref{eq:AsymExp:EKF}.  If $\LC(f)$ happens to have only one nonzero term, then so does $\LF(f)$, and in this situation we may also write $\LT(f) := \LF(f)$ and call this the \DEF{leading term} of $f$.
\end{defs}

The next corollary follows readily from Lemma \ref{lem:noncomp}\ref{lem:noncomp:Seq} and Proposition \ref{prop:AsymExp}.

\begin{cor}\label{cor:CKF:AsymExpKernel}
The following two statements are equivalent for any $f\in\C^{\KK,\FTs}(\RR_+)$.
\begin{enumerate}
\item\label{cor:CKF:AsymExpKernel:zero}
The asymptotic expansion for $f$ at $0^+$ is zero.

\item\label{cor:CKF:AsymExpKernel:rapidDecay}
For every $\rho > 0$, $f(y) = o(y^\rho)$ as $y\to 0^+$.
\end{enumerate}
If we further assume that $f\in\C^{\KK,\i\S}(\RR_+)$, then statements \ref{cor:CKF:AsymExpKernel:zero} and \ref{cor:CKF:AsymExpKernel:rapidDecay} are also equivalent to the following statement.

\begin{enumerate}{\setcounter{enumi}{2}
\item\label{cor:CKF:AsymExpKernel:f=zero}
The germ of $f$ at $0^+$ is zero.
}\end{enumerate}
\end{cor}

Since, for example, the function $y\mapsto \e^{-1/y}$ is in $\C^{\KK,\FTs}(\RR_+)$, statement \ref{cor:CKF:AsymExpKernel:f=zero} is not equivalent to \ref{cor:CKF:AsymExpKernel:zero} and \ref{cor:CKF:AsymExpKernel:rapidDecay} in the larger category $\C^{\KK,\FTs}$.

\begin{rems}\label{rem:AsymExp}
The next two remarks explain how the Definitions \ref{def:AsymExp:CKF} and Proposition \ref{prop:AsymExp} readily adapt to give one-sided asymptotic expansions at any point in the extended real line and that, for parametric families of functions, these asymptotic expansions can be expressed in a piecewise manner over a finite partition of the parameter space into subanalytic sets.
\begin{enumerate}
\item\label{rem:AsymExp:Nonparam}
Let $f\in\C^{\KK,\FTs}(\RR)$, and suppose that $\ell$ is the object specified in one of the following two cases.
\begin{enumerate}
\item[]
Case A: $\ell = \theta^\sigma$ for some $\theta\in\RR$ and $\sigma\in\{+,-\}$.

\item[]
Case B: $\ell = \sigma\infty$ for some $\sigma\in\{+,-\}$.
\end{enumerate}
Let $y_{\ST}$ and $y_{\AF}$ be the standard and affine coordinates for $\ell$, as defined in Remark \ref{rem:PrepData}\ref{rem:PrepData:LimBndCoord}.  The pullback of $f$ by the map $(x,y_{\ST})\mapsto (x,y)$ has an asymptotic expansion at $0^+$, which we may then pullback by the inverse map $(x,y)\mapsto (x,y_{\ST})$ to obtain an asymptotic expansion
\begin{equation}\label{eq:AsymExp:yS}
\sum_{i\in \N} y_{\ST}^{r_i}(\log y_{\ST})^{s_i} E_i(y_{\ST})
\end{equation}
for $f(y)$ as $y\to\ell$.  We call \eqref{eq:AsymExp:yS} the \DEF{asymptotic expansion for $f$ at $\ell$ in standard coordinates}.  By substituting $y_{\ST} = y_{\AF}$ in  \eqref{eq:AsymExp:yS} in Case A and substituting $y_{\ST} = y_{\AF}^{-1}$ in  \eqref{eq:AsymExp:yS} in Case B, we may also express this series in the form
\begin{equation}\label{eq:AsymExp:yA}
\sum_{i\in \N} y_{\AF}^{\tld{r}_i}(\log y_{\AF})^{s_i} \tld{E}_i(y_{\AF}),
\end{equation}
which we call the \DEF{asymptotic expansion for $f$ at $\ell$ in affine coordinates}.  The series \eqref{eq:AsymExp:yS} and \eqref{eq:AsymExp:yA} determine the data \eqref{eq:AsympExp:data} in the same way as when $\ell = 0^+$ in the Definitions \ref{def:AsymExp:CKF}.  However, this data depends on the choice of coordinates in Case B, for then $y_{\ST} = y_{\AF}^{-1}$.  In particular, note that in Case B the asymptotic scales for \eqref{eq:AsymExp:yS} and \eqref{eq:AsymExp:yA}
are respectively
\[
\textstyle \{y_{\ST}^{r}(\log y_{\ST})^{s} : (r,s)\in (\RE(\Gamma) + \frac{1}{d}\NN)\times S\}
\]
and
\[
\textstyle \{y_{\AF}^{r}(\log y_{\AF})^{s} : (r,s)\in (\RE(\Gamma) - \frac{1}{d}\NN)\times S\},
\]
as denoted in \eqref{eq:AsymExp:Scale}.

\item\label{rem:AsymExp:Param}
Now let us consider the parametric case.  Consider a subanalytic set $X\subseteq\RR^m$ and a finite set $\F\subseteq\C^{\KK,\FTs}(X\times\RR)$, and let $L$ be a finite set of functions $\ell$ on $X$ of the following two possible forms.
\begin{enumerate}
\item[]
Case A: There exist $\theta\in\S(X)$ and $\sigma\in\{+,-\}$ such that $\ell(x) = \theta(x)^\sigma$ for all $x\in X$.

\item[]
Case B: There exists a $\sigma\in\{+,-\}$ such that $\ell(x) = \sigma\infty$ for all $x\in X$.
\end{enumerate}
Then there exists a finite partition $\B$ of $X$ into subanalytic sets such that for each $B\in\B$ and $\ell\in L$, there exist prepared coordinate data $(A_{B}^{\ell},\theta_{B}^{\ell},d_{B}^{\ell},\psi_{B}^{\ell})$ with base $B$ and limiting boundary $\ell\restriction{B}$ such each $f\in\F$ has an asymptotic expansion of the form
\[
\sum_{i\in \N} y_{\ST}^{r_i}(\log y_{\ST})^{s_i} E_i(x,y_{\ST})
\]
on $A_{B}^{\ell}$ as $y\to\ell$. For each $x\in B$, the resulting series in $y$ is the asymptotic expansion for $f_x$ at $\ell(x)$ in standard coordinates, as just defined in Remark \ref{rem:AsymExp}\ref{rem:AsymExp:Nonparam}.  When $\ell(x)\neq\ell'(x)$ for all $x\in X$ and distinct $\ell,\ell'\in L$, we may construct $\{(A_{B}^{\ell},\theta_{B}^{\ell},d_{B}^{\ell},\psi_{B}^{\ell})\}_{(B\times\ell)\in\B\times L}$ to have disjoint domains, and of course, we may also express these expansions in affine coordinates $(x,y_{\AF})$.
\end{enumerate}
\end{rems}

\begin{prop}[Pointwise Limits]\label{prop:Limits}
Let $\D\in\{\S^{\KK},\C^{\KK},\C^{\KK,\i\S},\C^{\KK,\FTs}\}$, and consider a function $f\in\D(X\times\RR_+)$ for some subanalytic set $X\subseteq\RR^m$.  If the limit $\displaystyle g(x) := \lim_{y\to 0^+} f(x,y)$ exists for all $x\in X$, then $g\in\D(X)$.
\end{prop}

\begin{proof}
The case of $\D = \C^{\KK,\FTs}$ is proven in \cite[Theorem 7.11]{clcose:fourier_mellin_power_constructible} (and in the same manner in \cite[Proposition 7.7]{ClCoRoSe} for $\KK=\QQ$) using Proposition \ref{prop:AsymExp} with $\D=\C^{\KK,\FTs}$ and statements \ref{lem:noncomp:Seq} and \ref{lem:noncomp:TwoSeq} of Lemma \ref{lem:noncomp}.  The other cases follow from the same argument.  Indeed, the only difference between between \cite{clcose:fourier_mellin_power_constructible} and the present situation is that we have formulated Proposition \ref{prop:AsymExp} for all four of our categories rather than for $\C^{\KK,\FTs}$ alone, which then allows the present proposition to go through verbatim for all four of these categories.
\end{proof}

\section{Consequences of Preparation and Asymptotic Expansions}\label{s:ConseqPrepAsym}

This section applies the results of Section \ref{s:Background} to prove a number of technical statements that are used throughout the rest of the paper.  We extensively use the leading data notation for asymptotic expansions from the Definitions \ref{def:AsymExp:CKF} (such as $\LE_1(f)$, etc.) in both the original form at $0^+$ and in the generalized forms at $\theta^\sigma$ or $\sigma\infty$ as described in the Remarks \ref{rem:AsymExp}.

\begin{lem}\label{lem:LE1}
Consider $f\in\C^{\KK,\FTs}(X\times\RR_+)$, where $X\subseteq\RR^m$ is subanalytic.
\begin{enumerate}
\item\label{lem:LE1:min}
The number $\overline{r} := \min\{\LE_1(f_x) : x\in X\}$ exists.

\item\label{lem:LE1:BigO}
There exists a function $b\in\S_+(X)$ such that for each $r < \overline{r}$, there exists a function $C_r:X\to(0,+\infty)$ such that
\[
|f(x,y)| \leq C_r(x) y^r
\]
for all $(x,y)\in X\times\RR_+$ with $y < b(x)$.

\item\label{lem:LE1:lim0}
We have $\overline{r} > 0$ if and only if $\lim_{y\to 0^+}f(x,y) = 0$ for all $x\in X$.

\item\label{lem:LE1:L1}
We have $\overline{r} > -1$ if and only if $f_x$ is locally integrable at $0^+$ for all $x\in X$.
\end{enumerate}
\end{lem}

\begin{proof}
Apply Proposition \ref{prop:AsymExp} to $\{f\}$, consider any $B\in\B$ and $(A,0,d,\psi)$ as in the statement of the proposition, and consider the asymptotic expansion for $f$ given in \eqref{eq:AsymExp}.  It suffices to prove the lemma with $B$ in place of $X$.  Statement \ref{lem:LE1:min} follows from the fact that $\{\LE_1(f_x) : x\in B\}$ is contained in the well-ordered set $\RE(\Gamma) + \frac{1}{d}\NN$ for some finite $\Gamma\subseteq\KK$.  Statement \ref{lem:LE1:BigO} follows from the fact that $f_x(y) = O(\LM(f_x)) = o(y^r)$ as $y\to 0^+$ for each $x\in B$ and $r < \overline{r}$, where the function $b\in\S_+(X)$ in \ref{lem:LE1:BigO} is used to define $A = \{(x,y)\in B\times\RR : 0 < y < b(x)\}$.  Statement \ref{lem:LE1:lim0} follows from the fact that for each $x\in X$, it is clear that $\lim_{y\to 0^+}f(x,y) = 0$ if $\LE_1(f_x) > 0$, and the converse follows from Lemma \ref{lem:noncomp}\ref{lem:noncomp:Seq}.  Statement \ref{lem:LE1:L1} follows similarly, except apply Lemma \ref{lem:noncomp}\ref{lem:noncomp:L1}.
\end{proof}

With the obvious adjustments, we may apply Lemma \ref{lem:LE1} more generally with $\theta^\sigma$ or $\sigma\infty$ in place of $0^+$, as in Case A or B of Remark \ref{rem:AsymExp}\ref{rem:AsymExp:Param}, in either standard coordinates $(x,y_{\ST})$ or affine coordinates $(x,y_{\AF})$.

\begin{lem}\label{lem:CK:L1}
Let $f\in\C^{\KK}(\RR)$.
\begin{enumerate}
\item\label{lem:CK:L1:Infty}
Let $\sigma\in\{+,-\}$.  If $f$ is integrable at $\sigma\infty$, then so is $f^{(k)}$ for all $k\in\NN$.

\item\label{lem:CK:L1:Point}
If $f$ is continuous at a point $a\in\RR$, then $f'$ is locally integrable at $a$.

\item\label{lem:CK:L1:Cn}
Let $n\in\NN$.  If $f\in \L^1(\RR)\cap C^n(\RR)$, then $f^{(k)}\in \L^1(\RR)$ for all $k\in\{0,\ldots,n+1\}$.
\end{enumerate}
\end{lem}

\begin{proof}
To prove \ref{lem:CK:L1:Infty}, let $\sigma\in\{+,-\}$, assume that $f\in\C^{\KK}(\RR)$ is integrable at $\sigma\infty$, and let
\begin{equation}\label{eq:CK:L1}
f(y) = \sum_{i\in \N} y_{\AF}^{r_i}(\log y_{\AF})^{s_i} E_i(y)
\end{equation}
be the convergent asymptotic expansion for $f$ at $\sigma\infty$ using the affine coordinates $y_{\AF} = \sigma y$.  Then
$\LE_1(f) < -1$ by Lemma \ref{lem:LE1}\ref{lem:LE1:L1}, where $\LE_1(f)$ is the leading exponent of $y_{\AF}$ for $f$ at $\sigma\infty$.  Differentiating \eqref{eq:CK:L1} preserves its general form but decreases the leading exponent of $y_{\AF}$, which preserves integrability at $\sigma\infty$.  This prove \ref{lem:CK:L1:Infty}.

To prove \ref{lem:CK:L1:Point}, assume that $f$ is continuous at $a$, and let $f_{a^\sigma}(y) = f(a+\sigma y) - f(a)$ for each $\sigma\in\{+,-\}$.  Then for each $\sigma\in\{+,-\}$, we have $\lim_{y\to 0^+} f_{a^{\sigma}}(y) = 0$, so $\LE_1(f_{a^\sigma}) > 0$ by Lemma \ref{lem:LE1}\ref{lem:LE1:lim0}, where this leading exponent is of $y$ at $0^+$.  Hence, $\LE_1(f'_{a^{\sigma}}) > -1$, so $f'$ is locally integrable at $a$ by Lemma \ref{lem:LE1}\ref{lem:LE1:L1}, which proves \ref{lem:CK:L1:Point}.

Finally, to prove \ref{lem:CK:L1:Cn}, note that if $n\in\NN$, $f\in C^n(\RR)$, and $k\in\{0,\ldots,n+1\}$, then $f^{(k)}\in \L^1(\RR)$ because $f^{(k)}$ is integrable on any compact interval by the continuity of $f^{(k)}$ when $k\leq n$ and by \ref{lem:CK:L1:Point} when $k = n+1$, and because $f^{(k)}$ is integrable at $\pm\infty$ by \ref{lem:CK:L1:Infty}.
\end{proof}

\begin{exam}\label{exam:CKiS:notL1}
Lemma \ref{lem:CK:L1} fails in the category $\C^{\KK,\i\S}$ because integrability at $\pm\infty$ does not, in general, impose a restriction on the behavior of the phase functions in a formula defining a function $f\in\C^{\KK,\i\S}(\RR)$, and the derivatives of these phase functions appear outside of the phase in the formula for $f'$ because of the chain rule.  For example, for any given $n\in\NN$, the function $f\in\C^{\KK,\i\S}(\RR)$ defined by
\[
f(t) = \frac{\e^{\i t^n}}{t^2+1}, \quad\text{for all $t\in\RR$,}
\]
is in $\L^1(\RR)$, but $f'\not\in\L^1(\RR)$ when $n\geq 2$.
\end{exam}

We use the next lemma in Section \ref{s:HolomMeromExt} to show that smooth $\C^{\KK,\i\S}$-functions are analytic.

\begin{lem}\label{lem:CKis:smoothExt}
Let $b > 0$, and consider a function $f\in\C^{\KK,\i\S}((0,b))$ whose asymptotic expansion at $0^+$ converges to $f$ on $(0,b)$.  If $f$ extends to a $C^\infty$ function on $[0,b)$, then the asymptotic expansion of $f$ at $0^+$ is a power series that converges on the disk $\{z\in\CC : |z| < b\}$.
\end{lem}

\begin{proof}
Assume that $f$ extends to a $C^\infty$ function on $[0,b)$.  By \eqref{eq:CKiS:preAsympExp} we may write the asymptotic expansion for $f$ at $0^+$ as a convergent series
\begin{equation}\label{eq:CKis:smoothExt:AsymExp}
f(y) = \sum_{(r,s,\phi)\in R\times S\times\Phi} f_{(r,s,\phi)}(y)
\end{equation}
on $(0,b)$ with
\[
f_{(r,s,\phi)}(y) = c_{(r,s,\phi)}y^r(\log y)^s\e^{\i\phi(y^{-1/d})}
\]
for each $(r,s,\phi)\in R\times S\times \Phi$, where $d$ is a positive integer, $R = \Gamma + \frac{1}{d}\NN$ for some finite $\Gamma\subseteq\KK$, $S$ is a finite initial segment of $\NN$, $\Phi$ is a finite subset of $\{\phi\in\RR[y] : \phi(0) = 0\}$ with $0\in\Phi$, and $c_{r,s,\phi}\in\CC$ for all $(r,s,\phi)\in R\times S\times\Phi$.  Our goal is to show that the support
\[
\Supp(f) = \{(r,s,\phi)\in R\times S\times\Phi : c_{(r,s,\phi)} \neq 0\}
\]
of the series \eqref{eq:CKis:smoothExt:AsymExp} is contained in $\NN\times\{0\}\times\{0\}$, for that means that \eqref{eq:CKis:smoothExt:AsymExp} is a power series, which then necessarily converges on the disk $\{z\in\CC : |z| < b\}$ since it converges on the interval $(0,b)$.  We divide the proof into two steps, first showing that $\Supp(f) \subseteq R\times S\times\{0\}$ and then showing that $\Supp(f)\subseteq \NN\times\{0\}\times\{0\}$, and we use notation from the Definitions \ref{def:AsymExp:CKF}.

Suppose for a contradiction that $\Supp(f) \not\subseteq R\times S\times\{0\}$. Consider any $(r,s,\phi)\in \Supp(f)$, and write $\phi(y) = \sum_{j\in J} \phi_j y^j$ for some finite $J\subseteq\ZZ_+$ and coefficients $\{\phi_j\}_{j\in J}\subseteq\RR$.  Then
\[
f_{(r,s,\phi)}(y) = c_{(r,s,\phi)}y^r(\log y)^s\left(\prod_{j\in J} \e^{\i\phi_j y^{-j/d}}\right),
\]
so the product rules gives
\begin{align*}
f'_{(r,s,\phi)}(y)
    &= c_{(r,s,\phi)}r y^{r-1}(\log y)^s\left(\prod_{j\in J} \e^{\i\phi_j y^{-j/d}},\right) \\
    &\phantom{=}+ c_{(r,s,\phi)}s y^{r-1}(\log y)^{s-1}\left(\prod_{j\in J} \e^{\i\phi_j y^{-j/d}},\right) \\
    &\phantom{=}+ \sum_{j\in J} c_{(r,s,\phi)}\left(-\frac{\i j\phi_j}{d}\right)y^{r-(1+j/d)}(\log y)^s\left(\prod_{k\in J} \e^{\i\phi_k y^{-k/d}}\right).
\end{align*}
We see that
\[
\LE_1(f'_{(r,s,\phi)}) =
\begin{cases}
\LE_1(f_{(r,s,\phi)}) - \left(1 + \frac{\deg(\phi)}{d}\right),
    & \text{if $\phi\neq 0$,} \\
\LE_1(f_{(r,s,\phi)}) - 1,
    &\text{if $\phi = 0$.} \\
\end{cases}
\]
Therefore, by letting $\mu = \max\{\deg(\phi) : (r,s,\phi)\in\Supp(f)\}$ (where $\mu\geq 1$ because $\Supp(f) \not\subseteq R\times S\times\{0\}$) and by defining
\[
\Supp_{\mu}(f) = \{(r,s,\phi)\in\Supp(f) : \deg(\phi) = \mu\}
\]
and
\[
f_{\mu}(y) = \sum_{(r,s,\phi)\in\Supp_{\mu}(f)} f_{(r,s,\phi}(y),
\]
we see that for some sufficiently large $k\in\NN$,
\[
\LE_1(f^{(k)}) = \LE_1(f^{(k)}_{\mu}) < 0.
\]
So by Lemma \ref{lem:noncomp}\ref{lem:noncomp:Seq}, there exists a sequence $\{t_i\}_{i\in\NN}$ in $(0,b)$ converging to $0$ such that $\lim_{i\to\infty}|f^{(k)}(t_i)| = +\infty$, which is contrary to our assumption that $f$ extends to a $C^\infty$ function on $[0,b)$.  So, in fact, $\Supp(f) \subseteq R\times S\times\{0\}$.

We may now write $f$ in the simpler form
\[
f(y) = \sum_{(r,s)\in R\times S} f_{(r,s)}(y)
\]
with $f_{(r,s)}(y) = c_{(r,s)} y^r(\log y)^s$ for each $(r,s)\in R\times S$, and accordingly, we now define $\Supp(f) = \{(r,s)\in R\times S : c_{(r,s)}\neq 0\}$.  If $\LE_1(f) \not\in\NN$, then for some sufficiently large $k\in\NN$, we have $\LE_1(f^{(k)}) < 0$, and we again apply Lemma \ref{lem:noncomp}\ref{lem:noncomp:Seq} to arrive at the contradictory conclusion that $\lim_{i\to\infty} |f^{(k)}(t_i)| = +\infty$ for some sequence $\{t_i\}_{i\in\NN}$ in $(0,b)$ converging to $0$.  So, in fact, $\LE_1(f)\in\NN$.  If $\LE_2(f) > 0$, then $\LE(f^{(\LE_1(f))}) = (0,\LE_2(f))$, and we yet again arrive at the contradictory conclusion that $\lim_{i\to\infty} |f^{(\LE_1(f))}(t_i)| = +\infty$ for some sequence $\{t_i\}_{i\in\NN}$ in $(0,b)$ converging to $0$.  So, in fact, $\LE_2(f) = 0$.

To summarize, we have shown that $f$ has a leading term of the form $\LT(f) = ct^r$ for some $c\neq 0$ and $r\in\NN$.  Letting $f_0 := f$ and recursively applying the same fact to $f_{i+1} := f_{i} - \LT(f_i)$, for each $i\in\NN$, shows that all nonzero terms of $f$ are of this form and thereby establishes that $\Supp(f)\subseteq\NN\times\{0\}$.
\end{proof}

Versions of the next proposition are found in nonparametric form in \cite[Remark 5.8]{ClCoRoSe} and \cite[Proposition 2.8]{AiClRaSe} for $\C^{\QQ,\FTs}$, albeit with less explicit proofs.  The proposition and its corollary use the Notation \ref{notn:Background:CoordProjFiber}.

\begin{prop}\label{prop:CFK:anal}
Let $\D\in\{\S^{\KK},\C^{\KK},\C^{\KK,\i\S},\C^{\KK,\FTs}\}$ and $f\in\D(D)$ for a subanalytic set $D\subseteq\RR^{m+n}$ such that $D_x$ has nonempty interior in $\RR^n$ for each $x\in\pi(D)$.
\begin{enumerate}
\item\label{prop:CFK:anal:Deriv}
If $E\subseteq D$ is a subanalytic set such that for each $x\in\pi(E)$, $E_x$ is open in $\RR^n$ and $f_x$ is $C^\infty$ on $E_x$, then $\frac{\partial^{|\alpha|} f\restriction{E}}{\partial y^{\alpha}} \in \D(E)$ for all $\alpha\in\NN^n$.

\item\label{prop:CFK:anal:set}
There exists a subanalytic set $E\subseteq D$ such that $\pi(E) = \pi(D)$ and such that for each $x\in \pi(E)$, the fiber $E_x$ is open $\RR^n$, $\dim(D_x\setminus E_x) < n$, and $f_x$ is analytic on $E_x$.
\end{enumerate}
\end{prop}

\begin{proof}
Since partial derivatives are defined as limits of difference quotients, statement \ref{prop:CFK:anal:Deriv} follows immediately from repeated application of Proposition \ref{prop:Limits}.  Statement \ref{prop:CFK:anal:set} for $\S$ (namely, $\S^{\KK}$ with $\KK= \QQ$) follows from the fact that $\RR_{\an}$ has analytic cell decomposition (which, for instance, follows from Proposition \ref{prop:SubPrep}), and statement \ref{prop:CFK:anal:set} then follows readily for $\S^{\KK}$, $\C^{\KK}$, and $\C^{\KK,\i\S}$ from the way these categories are defined from $\S$.  It remains to prove statement \ref{prop:CFK:anal:set} for $\C^{\KK,\FTs}$.

So assume that $\D = \C^{\KK,\FTs}$, and write
\[
f(x,y) = \f[g](x,y) = \int_{\RR} g(x,y,t)\e^{\i t}\d t
\]
on $D$ for some $g\in\C^{\KK}(D\times\RR)$ with $g_{(x,y)}\in \L^1(\RR)$ for all $(x,y)\in D$.  Apply Proposition \ref{prop:CKprep} to prepare $g$ with respect to the variable $t$, let $\A$ be the resulting set of prepared coordinate data decomposing $D\times\RR$ over $D$, and let $\B$ be the collection of base sets of $\A$.  By further partitioning these base sets, we may assume that each member of $\B$ is a cell over $\RR^m$.  Let $E$ be the union of the cells in $\B$ whose fibers over $\RR^m$ are open in $\RR^n$.  Then $\pi(E) = \pi(D)$ and for each $x\in \pi(E)$, $E_x$ is open $\RR^n$ and $\dim(D_x\setminus E_x) < n$.  We claim that $f_x$ is also analytic on $E_x$ for each $x\in\pi(E)$.

Indeed, consider any $(A,\theta,d,\psi)\in\A$ whose base set $B$ is contained in $E$.  It suffices to prove the claim with $g$ replaced by $g\cdot\chi_A$ and $D$ replaced by $B$ in the defining equation $f=\f[g]$ for $f$, where $\chi_A$ is the characteristic function of $A$.  Write the $\C^{\KK}$-prepared form of $g$ on $A$ using the affine coordinates $(x,y,t_{\AF})$ for $(A,\theta,d,\psi)$, where $t_{\AF} = \sigma(t-\theta(x,y))$, and observe from Remark \ref{rem:CKiS:FiniteAsymExp} and Lemma \ref{lem:LE1}\ref{lem:LE1:L1} that if $(A,\theta,d,\psi)$ has a limiting boundary (as in Case 2.A or 2.B of the Definitions \ref{def:PrepData}), then we may write this prepared form so that each of the powers of $t_{\AF}$ are in $(-1,+\infty)$ or $(-\infty,-1)$ according to whether the limiting boundary of $(A,\theta,d,\psi)$ is finite or infinite (i.e., $\theta^\sigma$ or $\sigma\infty$), respectively.  Therefore, each term in the prepared form of $g$ is integrable, so it suffices to prove the claim with $g$ replaced by any one of these terms in the defining equation $f = \f[g]$ for $f$. Change coordinates in the integral by pulling back by $(x,y,t)\mapsto (x,y,\sigma t + \theta(x,y))$, which is the inverse of the map $(x,y,t)\mapsto(x,y,t_{\AF})$.  We may thereby assume that $(A,0,d,\psi)$ has center $0$, and this changes $\e^{\i t}$ to $\e^{\i(\sigma t + \theta(x,y))}$ in the integrand.  Factor $\e^{\i\theta(x,y)}$ out of the integral.

Thus, we have reduced to the case that
\[
f(x,y) = c(x,y)\e^{\i\theta(x,y)} \int_{a(x,y)}^{b(x,y)} t^r(\log t)^s h(x,y,t) \e^{\i\sigma t} \d t
\]
on $B$ for some analytic functions $c\in\C^{\KK}(B)$ and $\theta\in\S(B)$, exponents $r\in\frac{1}{d}\ZZ$ and $s\in\NN$, $\psi$-function $h$, and sign $\sigma\in\{+,-\}$, where $a\in\S(B)$ and $b\in\S(B)\cup\{+\infty\}$ are analytic functions with $A = \{(x,y,t)\in B\times\RR: a(x,y) < t < b(x,y)\}$, and exactly one of the following cases holds.

Case 1: $0 < a(x,y) < b(x,y) < +\infty$ on $B$.

Case 2.A: $0 = a(x,y) < b(x,y) < +\infty$ on $B$ and $r > -1$.

Case 2.B: $0 < a(x,y) < b(x,y) = +\infty$ on $B$ and $r < -1$.

\noindent
Write $\overline{A}$ for the closure of $A$ in $B\times\RR_+$.  In Cases 1 and 2.A, we may define a holomorphic function
\[
F(x,y,z) = c(x,y) e^{\i\theta(x,y)} \int_{a(x,y)}^{z} w^r(\log w)^s h(x,y,w) e^{\i\sigma w} \d w
\]
on some neighborhood $U$ of $\overline{A}$ in $\{(x,y,z)\in B\times\CC : \RE(z) > 0\}$, with the understanding that this integral is defined via a path-independent choice of contour from $a(x,y)$ to $z$ in $U_{(x,y)}$, for each $(x,y)\in B$, and that in Case 2.A the integral is improper at $0$ and convergent since $r > -1$.  In Case 2.B, instead define the holomorphic function
\[
F(x,y,z) = -c(x,y) e^{\i\theta(x,y)} \int_{z}^{+\infty} w^r(\log w)^s h(x,y,w) e^{\i\sigma w} \d w
\]
on $U$, which converges since $r<-1$.  Then in Case 1 we have $f(x,y) = F(x,y,b(x,y)) - F(x,y,a(x,y))$ on $B$, which is analytic, and likewise in Case 2.A with $F(x,y,0) := 0$ and in Case 2.B with $F(x,y,+\infty) := 0$.  It follows that $f$ is analytic on $B$, which proves the claim and so also the proposition.
\end{proof}

\begin{defn}\label{def:a.e.Property}
Consider an open set $U\subseteq\RR^m$ and a function $f:U\to\CC$.  We say that the a.e.-equivalence class of $f$,
\[
[f] = \{g:U\to\CC \mid \text{$g(x)=f(x)$ for almost all $x\in U$}\},
\]
\DEF{has a property $P$} when there exists a function in $[f]$ with the property $P$.
\end{defn}

Observe that when $[f]$ is continuous, the equivalence class $[f]$ contains exactly one continuous function.

\begin{cor}\label{cor:aeClass}
Let $\D\in\{\S^{\KK},\C^{\KK},\C^{\KK,\i\S},\C^{\KK,\FTs}\}$, and consider a subanalytic set $D\subseteq\RR^{m+n}$ whose fibers over $\RR^m$ are open in $\RR^n$ and a function $f\in\D(D)$ such that $[f_x]$ is continuous for each $x\in\pi(D)$.  Then there exists a $g\in\D(D)$ such that for each $x\in\pi(D)$, the function $g_x$ is continuous and contained in $[f_x]$.
\end{cor}

A simple but noteworthy consequence of this corollary is that when $f\in\D(\RR)$ and $P$ is a property that implies continuity (e.g., continuous, $C^k$, smooth, analytic), then $[f]$ has property $P$ if and only if there exists a $g\in\D(\RR)$ with property $P$ such that $f(y) = g(y)$ for all but finitely many $y\in\RR$.

\begin{proof}
Let $g:D\to\CC$ be the function such that for each $x\in\pi(D)$, $g_x$ is continuous and contained in $[f_x]$, and let $E$ be the subset of $D$ given by applying Proposition \ref{prop:CFK:anal}\ref{prop:CFK:anal:set} to $f$.  Then $g\restriction{E} = f\restriction{E} \in \D(E)$, so it suffices to show that $g\restriction{(D\setminus E)}\in\D(D\setminus E)$.  For each $x\in\pi(D)$, we have that $E_x\subseteq D_x$, $\dim(D_x\setminus E_x) < n$, and $D_x$ is open in $\RR^n$, so $E_x$ is dense in $D_x$.  Therefore for each $(x,y)\in D\setminus E$ and $\epsilon > 0$, there exists a $z\in E_x$ such that $|y-z| < \epsilon$.  So by definable choice for the o-minimal structure $\RR_{\an}$ (see Van den Dries \cite[Sec. 6.1, Prop 1.2]{vdD:Book:1998}), there exists a subanalytic map $h:(D\setminus E)\times\RR_+ \to E$ such that for all $(x,y,\epsilon)\in (D\setminus E)\times\RR_+$, we have $h(x,y,\epsilon) \in E_x$ and $|y-h(x,y,\epsilon)| < \epsilon$.  Therefore for all $(x,y)\in D\setminus E$,
\[
g(x,y) = \lim_{\epsilon\to 0^+} g(x,h(x,y,\epsilon)) = \lim_{\epsilon\to 0^+} f(x,h(x,y,\epsilon)),
\]
so $g\restriction{(D\setminus E)} \in \D(D\setminus E)$ by Proposition \ref{prop:Limits}.
\end{proof}

\section{Holomorphic and Meromorphic Extensions} \label{s:HolomMeromExt}

This section studies an assortment of interrelated theorems about the existence and properties of holomorphic and meromorphic extensions of parametric families of functions in the three categories $\S^{\KK}$, $\C^{\KK}$, and $\C^{\KK,\i\S}$.  The theorems are divided into three subsections discussing (1) holomorphic extensions at infinity, (2) holomorphic extensions to horizontal strips about the real axis, and (3) meromorphic extensions to the complex plane.  For the remainder of the section, we fix a subanalytic set $X\subseteq\RR^m$, and we extensively use the notation $\RR_\sigma$ and $\CC_\sigma$ from Notation \ref{notn:complex}, for each $\sigma\in\{+,-\}$, along with the leading data notation for asymptotic expansion (such as $\LE_1(f)$, etc.) defined in the Definitions \ref{def:AsymExp:CKF} and Remarks \ref{rem:AsymExp}.  The proofs in this section also make use of the following elementary remarks.

\begin{rem}\label{rem:ComplexCat}
\hfill
\begin{enumerate}
\item\label{rem:ComplexCat:power}
For any $\gamma\in\KK$, the principal power function $z\mapsto z^{\gamma}$ is an $\S^{\KK}$-function on $\CC_+$, being the product of the $\S^{\KK}$-function $z\mapsto |z|^{\gamma}$ and the subanalytic function $z\mapsto \e^{\i\Arg(z)}$.  If $\gamma\in\QQ$, then $z\mapsto|z|^\gamma$ is also subanalytic, so $z\mapsto z^\gamma$ is as well.

\item\label{rem:ComplexCat:Log}
The principal logarithm $z\mapsto\Log(z) = \log|z| + \i\Arg(z)$ is a $\C^{\QQ}$-function on $\CC_+$, being the sum of the $\C^{\QQ}$-function $z\mapsto \log|z|$ and the subanalytic function $z\mapsto \i\Arg(z)$.

\item\label{rem:ComplexCat:HolomCompactDef}
If $f$ is a holomorphic function on an open subset $U$ of $\CC^m$, and if $K\subseteq U$ is compact and subanalytic, then $f\restriction{K}$ is subanalytic.
\end{enumerate}
\end{rem}

\subsection{Holomorphic Extensions at Infinity} \label{s:HolomMeromExt:Infty}

The next theorem shows that for each $\D\in\{\S^{\KK},\C^{\KK},\C^{\KK,\i\S}\}$, any $\D$-function $f$ may be extended piecewise over a parameter space to analytic $\D$-functions $f^+$ and $f^-$ at $+\infty$ and $-\infty$, by which we mean on neighborhoods of $\infty$ in the Riemann sphere minus appropriate branch cuts.  The first additional statement in the theorem shows that the asymptotic bounds for $\C^{\KK,\FTs}$-functions given in Lemma \ref{lem:LE1} also hold for  the extensions $f^+$ and $f^-$ when $\D\in\{\S^{\KK},\C^{\KK}\}$.

\begin{thm}\label{thm:HolomExtInfty}
Let $f\in\D(X\times\RR)$, where $\D\in\{\S^{\KK},\C^{\KK},\C^{\KK,\i\S}\}$.  There exist a function $b\in\S_+(X)$ and a finite partition $\B$ of $X$ into subanalytic sets such that for each $B\in\B$ and $\sigma\in\{+,-\}$, the restriction of $f$ to $\{(x,y)\in B\times\RR : \sigma y > b(x)\}$ extends to an analytic $\D$-function
\[
f_{B}^{\sigma}:\{(x,z)\in B\times\CC_\sigma : |z| > b(x)\}\to\CC.
\]
For each $\sigma\in\{+,-\}$, define the $\D$-function
\[
f^\sigma :\{(x,z)\in X\times\CC_\sigma : |z| > b(x)\}\to\CC
\]
to be the union of the functions $\{f_{B}^{\sigma}\}_{B\in\B}$.  The following additionally hold when $\D\in\{\S^{\KK},\C^{\KK}\}$.
\begin{enumerate}
\item\label{thm:HolomExtInfty:CK}
If $\sigma\in\{+,-\}$ and we define
\[
\overline{r} = \max\{\LE_1(f_x) : x\in X\},
\]
where $\LE_1(f_x)$ is the leading exponent of $y_{\AF} = \sigma y$ for $f_x$ at $\sigma\infty$, then for each $r > \overline{r}$ there exists a function $C:X\to\RR_+$ such that
\begin{equation}\label{eq:HolomExtInfty:CK}
|f^{\sigma}(x,z)| \leq C(x)|z|^r
\end{equation}
on $\{(x,z)\in X\times\CC_\sigma : |z| > b(x)\}$.

\item\label{thm:HolomExtInfty:Prep}
Let $B\in\B$.  If the functions $f_{B}^{+}$ and $f_{B}^{-}$ glue together to form a single analytic function $f_B$ on $\{(x,z)\in B\times\CC : |z|>b(x)\}$, then we may write $f_B$ in the prepared form
\begin{equation}\label{eq:S:HolomExtInfty:MeromPrep}
f_B(x,z) = \sum_{i\in I} c_i(x) z^{r_i} f_i(x,z)
\end{equation}
on $\{(x,z)\in B\times\CC : |z|>b(x)\}$, where $I$ is a finite index set and for each $i\in I$, the function $c_i\in\D(B)$ is analytic, $r_i\in\ZZ$, and $f_i$ is a $\psi$-function with
\[
\psi(x,z) = \left(\psi_1(x),\ldots,\psi_N(x),\frac{b(x)}{z}\right).
\]
\end{enumerate}
\end{thm}

\begin{proof}
By applying Lemma \ref{lem:CKiSprep} at $\pm\infty$ rather than at $0^+$, as explained in Remark \ref{rem:AsymExp}\ref{rem:AsymExp:Param}, we obtain a function $b\in\S_+(X)$ and a finite partition $\B$ of $X$ into subanalytic sets such that for each $B\in\B$ and $\sigma\in\{+,-\}$, there exists prepared coordinate data $(A^\sigma,0,d,\psi^\sigma)$ lying over $B$ for which we may write $f$ in the prepared form
\begin{equation}\label{eq:CKiS:HolomExtInf:Prep}
f(x,y) = \sum_{i\in I^\sigma} c_i(x) y_{\AF}^{r_i}(\log y_{\AF})^{s_i} f_i(x,y)\e^{\i \phi_i(x,y)}
\end{equation}
on $A^\sigma := \{(x,y)\in B\times\RR : y_{\AF} > b(x)\}$, where \eqref{eq:CKiS:HolomExtInf:Prep} is as in equation \eqref{eq:CKiSprep}, except we are using the affine coordinates $(x,y_{\AF}) = (x,\sigma y)$ for $(A^\sigma,0,d,\psi^\sigma)$ and
\[
\psi^\sigma(x,y) = \left(\psi_1(x),\ldots,\psi_N(x), \left(\frac{b(x)}{y_{\AF}}\right)^{1/d}\right).
\]
Recall from Lemma \ref{lem:CKiSprep} that this means that $I^\sigma$ is a finite index set and that for each $i\in I^\sigma$, we have $c_i\in\D(B)$, $r_i\in\KK$, $s_i\in\NN$ (with $s_i=0$ when $\D=\S^{\KK}$), $f_i$ is a $\psi^\sigma$-function, and $\phi_i$ is a principal Laurent $\psi^\sigma$-polynomial (with $\phi_i=0$ when $\D\in\{\S^{\KK},\C^{\KK}\}$).  Notice that we are constructing the preparation so that the positive integer $d$ and the functions $b$ and $\psi_1,\ldots,\psi_N$ are the same for both $\sigma = +$ and $\sigma = -$, which is easily achievable.

For each $B\in\B$ and $\sigma\in\{+,-\}$, by using principal powers and the principal logarithm, the formula on the right side of \eqref{eq:CKiS:HolomExtInf:Prep} extends to define an analytic function $f^{\sigma}_{B}$ on $\{(x,z)\in B\times\CC_{\sigma} : |z| > b(x)\}$ that is easily seen to be a $\D$-function by the Remarks \ref{rem:ComplexCat}.  This completes the construction of $f^\sigma := \bigcup_{B\in\B}f_{B}^{\sigma}$ for each $\sigma\in\{+,-\}$ and thereby proves the main assertion of the theorem.  Additionally observe from this construction that the following holds:
\begin{equation}\label{eq:HolomExtInfty:Bnd}
\text{\parbox{4in}{For each $x\in X$ and $R > 0$, the restriction of the function $f^{\sigma}_{x}$ to the set $\{z\in\CC_\sigma : b(x) < |z| < b(x)+R\}$ is bounded.}}
\end{equation}

Now assume that $\D\in\{\S^{\KK}, \C^{\KK}\}$, define $\overline{r}$ as in the statement of the theorem, and let $B\in\B$.  Remark \ref{rem:CKiS:FiniteAsymExp} shows that for each $\sigma\in\{+,-\}$, the prepared form \eqref{eq:CKiS:HolomExtInf:Prep} may be constructed so that $\RE(r_i)\leq\overline{r}$ for all $i\in I^\sigma$.  It follows that for any $r > \overline{r}$ and $i\in I^\sigma$,
\[
(\sigma z)^{r_i}(\log \sigma z)^{s_i} = O(|z|^{\RE(r_i)}(\log|z|)^{s_i}) = o(|z|^r)
\]
as $z\to\infty$ in the set $\CC_\sigma$, which when combined with \eqref{eq:HolomExtInfty:Bnd} gives \eqref{eq:HolomExtInfty:CK}.  This proves Statement \ref{thm:HolomExtInfty:CK}.

Now additionally assume that $f_{B}^{+}$ and $f_{B}^{-}$ glue together to form a single analytic function $f_B$ on $\{(x,z)\in B\times\CC : |z|>b(x)\}$.  For each $x\in B$, we have two convergent asymptotic expansions for $(f_B)_x$ on $(b(x),+\infty)$ at $+\infty$, the first of which is constructed by expanding the $\psi$-functions in \eqref{eq:CKiS:HolomExtInf:Prep} and the second of which is the Laurent series representation of $(f_B)_x$ on $\{z\in\CC : |z| > b(x)\}$.  The form of the second expansion on $(b(x),+\infty)$ is a special case of the form of the first, which we know to be uniquely determined by $(f_B)_x$, so these two expansions must be identical.  Since this holds for all $x\in B$, equation \eqref{eq:CKiS:HolomExtInf:Prep} can be expressed in the simpler form \eqref{eq:S:HolomExtInfty:MeromPrep}, which proves Statement \ref{thm:HolomExtInfty:Prep}.
\end{proof}

\begin{cor}\label{cor:CK:EntirePoly}
A function $f\in\C^{\KK}(\RR)$ extends to an entire function if and only if $f$ is a polynomial with coefficients in $\cl_{\CC}(\KK)$. Therefore, the only $f\in\C^{\KK}(\RR)\cap\L^1(\RR)$ that extends to an entire function is the function $f=0$.
\end{cor}

\begin{proof}
Every polynomial with coefficients in $\cl_{\CC}(\KK)$ is a function in $\C^{\KK}(\RR)$ that extends to an entire function.  Conversely, Theorem \ref{thm:HolomExtInfty}\ref{thm:HolomExtInfty:CK} shows that every function in $\C^{\KK}(\RR)$ that extends to an entire function is polynomially bounded and is therefore a polynomial by Liouville's extended theorem, whose coefficients are necessarily in $\cl_{\CC}(\KK)$ by the definition of $\C^{\KK}(\RR)$.  This proves the first sentence of the corollary, and the second sentence follows immediately from the first since the only polynomial that is integrable on $\RR$ is the zero polynomial.
\end{proof}

\begin{exam}\label{exam:vertical sides in CKiS}
Theorem \ref{thm:HolomExtInfty}\ref{thm:HolomExtInfty:CK} fails in the larger category $\C^{\KK,\i\S}$.  For example, the function $f:\RR\to\CC$ defined by
\[
f(y) = \frac{\e^{\i y}}{1+y^2}, \quad\text{for all $y\in\RR$,}
\]
is in $\C^{\KK,\i\S}(\RR)$, and $\LE_1(f) = -2$ with respect to the affine coordinates $y_{\AF} = \pm y$ at both $+\infty$ and $-\infty$.  However, $f$ extends to a meromorphic function on $\CC$ with poles at $\pm\i$, and
\[
|f(-\i y)| = \frac{\e^y}{|1-y^2|}
\]
grows exponentially as $y\to+\infty$, so the bound given in Theorem \ref{thm:HolomExtInfty}\ref{thm:HolomExtInfty:CK} fails for $f$ for all $r\in\RR$.
\end{exam}

\begin{exam}\label{exam:CKis:notEntirePoly}
Corollary \ref{cor:CK:EntirePoly} also fails in the larger category $\C^{\KK,\i\S}$ because  $t\mapsto \e^{\i t}$ is a nonpolynomial function in $\C^{\KK,\i\S}(\RR)$ that extends to an entire function on $\CC$.
\end{exam}

Consider again a function $f\in\C^{\KK,\i\S}(X\times\RR)$ and a sign $\sigma\in\{+,-\}$.  The function $f^{\sigma} : \{(x,z)\in X\times\CC_\sigma : |z| > b(x)\}\to\CC$ constructed in Theorem \ref{thm:HolomExtInfty} is defined using the ray $\{y\in\RR : \sigma y\leq 0\}$ as a branch cut in the last coordinate.  The need for a branch cut can be eliminated by instead extending $f$ in a similar manner to a function defined on $\{(x,z)\in X\times\LL : |z| > b(x)\}$, where $\LL$ is the Riemann surface of the logarithm, but we have chosen to formulate Theorem \ref{thm:HolomExtInfty} in the complex plane with branch cuts to keep its statement more elementary.  We instead give the next proposition as a replacement for a version of Theorem \ref{thm:HolomExtInfty} involving $\LL$.

\begin{prop}\label{thm:HolomExtInfty:exp}
Let $f\in\C^{\KK,\i\S}(X\times\RR)$, $b\in\S_+(X)$, and $\B$ be as in Theorem \ref{thm:HolomExtInfty}, and define $E[f]:X\times\RR\to\CC$ by
\[
E[f](x,y) = f(x,\e^y) \quad\text{for all $(x,y)\in X\times\RR$.}
\]
\begin{enumerate}
\item\label{thm:HolomExtInfty:exp:CKiS}
For each $B\in\B$, the restriction of $E[f]$ to $\{(x,y)\in B\times\RR : y > \log b(x)\}$ extends to an analytic function
\[
E[f]_B:\{(x,z)\in B\times\CC : \RE(z) > \log b(x)\}\to\CC.
\]
Define $E[f]:\{(x,z)\in X\times\CC : \RE(z) > \log b(x)\}\to\CC$ to be the union of the functions in $\{E[f]_{B}\}_{B\in\B}$.

\item\label{thm:HolomExtInfty:exp:CK}
If $f\in\C^{\KK}(X\times\RR)$ and we define $\overline{r} = \max\{\LE_1(f_x) : x\in X\}$, where $\LE_1(f_x)$ is the leading exponent of $y_{\AF} = y$ for $f_x$ at $+\infty$, then for each $r > \overline{r}$ and $a > 0$, there exists a function $C:X\to\RR_+$ such that
\[
|E[f](x,z)| \leq C(x)\e^{r\RE(z)}
\]
on $\{(x,z)\in X\times\CC : \text{$|\IM(z)| < a$ and $\RE(z) > \log b(x)$}\}$.
\end{enumerate}
\end{prop}

\begin{proof}
Let $B\in\B$.  Equation \eqref{eq:CKiS:HolomExtInf:Prep} with $\sigma = +$ gives
\[
E[f](x,y) =  \sum_{i\in I^+} c_i(x)\e^{r_i y} y^{s_i} f_i(x,\e^y) \e^{\i \phi_i(x,\e^y)}
\]
on $\{(x,y)\in B\times\RR : y > \log b(x)\}$.  This function extends analytically to the set $\{(x,z)\in B\times\CC : \RE(z) > \log b(x)\}$, with the constraint on $\RE(z)$ arising from the fact that we need $\e^{\RE(z)} = |\e^z| > b(x)$ to hold for each function $(x,z)\mapsto f_i(x,\e^z)$ to be defined and analytic.  This proves statement \ref{thm:HolomExtInfty:exp:CKiS}.

Now assume that $f\in\C^{\KK}(X\times\RR)$, and let $a > 0$.  Then $E[f]$ is of the simpler form
\[
E[f](x,y) =  \sum_{i\in I^+} c_i(x)\e^{r_i y} y^{s_i} f_i(x,\e^y)
\]
on $\{(x,y)\in B\times\RR : y > \log b(x)\}$.  For any $R > 0$ and $x\in B$, the function $E[f]_x$ is bounded on
\[
\{z\in\CC : \text{$|\IM(z)| < a$ and $\log b(x) < \RE(z) < \log b(x) + R$}\}.
\]
And for each $r > \overline{r}$, $i\in I^+$, and $x\in B$,
\[
\e^{r_i z} z^{s_i} = O(\e^{r_i\RE(z)}|z|^{s_i}) = o(\e^{r\RE(z)})
\]
as $\RE(z)\to+\infty$ in $\{z\in\CC : |\IM(z)| < a\}$.  Statement \ref{thm:HolomExtInfty:exp:CK} follows.
\end{proof}

\subsection{Holomorphic Extensions to Horizontal Strips}\label{subsection:holomorphic extension to horizontal strips}

For any $\D\in\{\S^{\KK},\C^{\KK},\C^{\KK,\i\S}\}$, if $f\in\D(\RR)$ is smooth, then $f$ is in fact analytic and extends analytically to a $\D$-function on an open horizontal strip about the real axis in the complex plane.  This assertion can be easily proven as follows.
\vspace*{1ex}
\begin{addmargin}[3ex]{0ex}
\begin{proof}
Let $f\in\D(\RR)$ be smooth.  Applying Lemma \ref{lem:CKis:smoothExt} to $t\mapsto f(y+\sigma t)$, for each $y\in\RR$ and $\sigma\in\{+,-\}$, shows that $f$ is analytic.  Theorem \ref{thm:HolomExtInfty} shows that there is some $b>0$ such that the restriction of $f$ to $\{t\in\RR : |t| \geq b\}$ extends analytically to a $\D$-function on $\{z\in\CC : |\RE(z)| \geq b\}$.  Applying Theorem \ref{thm:HolomExtInfty} to $t\mapsto f(y+t^{-1})$ for each $y\in [-b,b]$, along with the compactness of $[-b,b]$, shows that the restriction of $f$ to $[-b,b]$ extends analytically to a $\D$-function on $\{z\in\CC : |\RE(z)| \leq b, |\IM(z)| < a\}$ for some $a>0$, from which it follows that $f$ extends analytically to a $\D$-function on the horizontal strip $\{z\in\CC : |\IM(z)| < a\}$.
\end{proof}
\end{addmargin}
\vspace*{1ex}
The next theorem formulates this assertion more generally for parametric families of functions.  Its proof requires careful use of cell decomposition arguments in place of the compactness of $[-b,b]$ used in the nonparametric case.

\begin{thm}\label{thm:CKiS:HolomExtStrip}
Let $\D\in\{\S^{\KK},\C^{\KK},\C^{\KK,\i\S}\}$, and assume that $f\in\D(X\times\RR)$ is such that $f_x\in C^\infty(\RR)$ for all $x\in X$.   There exist a function $a\in\S_+(X)$ and a finite partition $\B$ of $X$ into subanalytic sets such that for each $B\in\B$, the restriction of $f$ to $B\times\RR$ is analytic and extends to an analytic $\D$-function on $\{(x,z)\in B\times\CC : |\IM(z)|<a(x)\}$.
\end{thm}

Notice from Remark \ref{rem:intro:CategFctsProps}\ref{rem:intro:CategFctsProps:CKF:Schwartz&Bump} that $\C^{\KK,\FTs}$ contains smooth functions that are not analytic, so Theorem \ref{thm:CKiS:HolomExtStrip} does not hold in the larger class $\C^{\KK,\FTs}$.  The proof of Theorem \ref{thm:CKiS:HolomExtStrip} is based on the following lemma.

\begin{lem}\label{lem:CKiS:HolomExtStrip}
Let $\D\in\{\S^{\KK},\C^{\KK},\C^{\KK,\i\S}\}$, assume that $f\in\D(X\times\RR)$ is such that $f_x\in C^\infty(\RR)$ for all $x\in X$, and define $g:X\times\RR\times\RR\to\CC$ by
\[
g(x,y,t) = f(x,y+t)
\]
for all $(x,y,t)\in X\times\RR\times\RR$.  There exist a function $\rho\in\S_+(X\times\RR)$ and a finite partition $\A$ of $X\times\RR$ into subanalytic sets such that for each $A\in\A$, the restriction of $g$ to $\{(x,y,t)\in A\times\RR : |t| < \rho(x,y)\}$ is analytic and extends to an analytic $\D$-function on $\{(x,y,z)\in A\times\CC : |z| < \rho(x,y)\}$.
\end{lem}

\begin{proof}
Let $f$ be as hypothesized.  Applying Theorem \ref{thm:HolomExtInfty} to the function $(x,y,t)\mapsto g(x,y,t^{-1})$ and then pulling back by the transformation $(x,y,t)\mapsto(x,y,t^{-1})$ constructs a function $\rho\in\S_+(X\times\RR_+)$ and a finite partition $\A$ of $X\times\RR$ into subanalytic sets such that for each $A\in\A$ and $\sigma\in\{+,-\}$, the restriction of $g$ to $\{(x,y,t)\in A\times\RR_{\sigma} : |t| < \rho(x,y)\}$ extends to an analytic $\D$-function on $\{(x,y,z)\in A\times\CC_{\sigma} : |z| < \rho(x,y)\}$.  Consider any $A\in\A$, $(x,y)\in A$, and $\sigma\in\{+,-\}$.  By construction, the asymptotic expansion for $g_{(x,y)}$ at $0^{\sigma}$ converges to $g_{(x,y)}$ on $\{t\in\RR_{\sigma} : |t| < \rho(x,y)\}$, and since $g_{(x,y)}\in C^\infty(\RR)$, Lemma \ref{lem:CKis:smoothExt} shows that this expansion is a power series.  These power series expansions for $g_{(x,y)}$ at $0^+$ and $0^-$ must both be the Taylor series for $g_{(x,y)}$ at $0$ and are therefore identical.  It follows that the restriction of $g$ to $\{(x,y,t)\in A\times\RR : |t| < \rho(x,y)\}$ extends to an analytic $\D$-function on $\{(x,y,z)\in A\times \CC : |z| < \rho(x,y)\}$.
\end{proof}

\begin{proof}[Proof of Theorem \ref{thm:CKiS:HolomExtStrip}]
Let $f$ be as hypothesized in Theorem \ref{thm:CKiS:HolomExtStrip}.   By Theorem \ref{thm:HolomExtInfty} there exist a function $b\in\S_+(X)$ and a finite partition $\B$ of $X$ into subanalytic sets such that for each $B\in\B$, the restriction of $f$ to $\{(x,t)\in B\times\RR : |t| > b(x)\}$ extends to an analytic $\D$-function on $\{(x,z)\in B\times\CC : |\RE(z)| > b(x)\}$.  Define the interval $I(x) = [-b(x),b(x)]$ for each $x\in X$.  To finish, it suffices to show that by further partitioning the sets in $\B$, we may construct a function $a\in\S_+(X)$ such that for each $B\in\B$, the restriction of $f$ to $\{(x,t)\in B\times\RR: \dist(t,I(x)) < a(x)\}$ extends to an analytic $\D$-function on $\{(x,z)\in B\times\CC: \dist(z,I(x)) < a(x)\}$, where $\dist(z,I(x))$ is the distance from $z$ to the set $I(x)$.

Apply Lemma \ref{lem:CKiS:HolomExtStrip} to $f$.  Let $g$, $\rho$, and $\A$ be as in the statement of the lemma, and write
\[
G:\{(x,y,z)\in X\times\RR\times\CC : |z|<\rho(x,y)\} \to \CC
\]
for the extension of $g:X\times\RR\times\RR\to\CC$ that is constructed in the lemma on $\{(x,y,z)\in A\times\CC : |z| < \rho(x,y)\}$ for each $A\in\A$.  By further partitioning the sets in $\A$ and $\B$, we may assume that $\A$ is a subanalytic cell decomposition of $X\times\RR$ over $X$ with $\B = \{\pi_m(A) : A\in\A\}$, that $\A$ is compatible with the graphs of the functions $\pm b$ on $X$, and that $\rho$ restricts to a continuous function on each member of $\A$.

Consider any $x\in X$.  For each $y\in\RR$, the restriction of $f_x(t) = g(x,y,t-y)$ to the set $\{t\in\RR : |t-y| < \rho(x,y)\}$ extends to an analytic function $F_{x}^{y}(z) := G(x,y,z-y)$ on $\{z\in\CC : |z-y| < \rho(x,y)\}$.  For any $y_1,y_2\in\RR$, the functions $F_{x}^{y_1}$ and $F_{x}^{y_2}$ both agree with $f_x$ on the intersection of their domains in $\RR$ and therefore agree on the intersection of their domains in $\CC$.  It follows that $f_x$ extends to an analytic function on $\{z\in\CC : \text{$|z-y| < \rho(x,y)$ for some $y\in\RR$}\}$.  By performing this construction for each $x\in X$, we obtain a function $F$ on the set
\[
D := \{(x,z)\in X\times\CC : \text{$|z-y| < \rho(x,y)$ for some $y\in\RR$}\}
\]
that extends $f$.

Now let $B\in\B$.  To finish, it suffices to define a function $a\in\S_+(B)$ such that
\begin{equation}\label{eq:F:RestrSet}
\{(x,z)\in B\times\CC : \dist(z,I(x)) < a(x)\}
\end{equation}
is a subset of $D$ and the restriction of $F$ to the set \eqref{eq:F:RestrSet} is an analytic $\D$-function.  We shall accomplish this by using the cell decomposition $\A$ to define a finite open cover of $\{(x,y)\in B\times\RR : |y|\leq b(x)\}$ in $D\cap(B\times\CC)$ by subanalytic sets such that $F$ restricts to an analytic $\D$-function on each subanalytic set from the cover and therefore restricts to an analytic $\D$-function on the union of the cover.  The specific way we shall define the cover will enable us to then easily define the desired function $a\in\S_+(B)$ so that the set \eqref{eq:F:RestrSet} is contained in the union of the cover, which will then complete the proof.

To that end, let
\[
-b = h_0 < h_1 < \cdots < h_{k-1} < h_k = b
\]
be the subanalytic functions on $B$ for which the sets
\[
\{(x,y)\in B\times\RR : y=h_i(x)\}, \,\,\text{for each $i\in\{0,\ldots,k\}$,}
\]
and
\[
\{(x,y)\in B\times\RR : h_{j-1}(x) < y< h_j(x)\}, \,\,\text{for each $j\in\{1,\ldots,k\}$,}
\]
are the cells in $\A$ contained in $\{(x,y)\in B\times\RR : |y| \leq b(x)\}$.  Consider any $i\in\{0,\ldots,k\}$.  The restriction of $F$ to the set
\[
D_i := \{(x,z)\in B\times\CC : |z-h_i(x)| < \rho(x,h_i(x))\} \subseteq D
\]
is an analytic $\D$-function given by the formula
\[
F(x,z) = G(x,h_i(x),z-h_i(x))
\]
on $D_i$.  Define $\delta_i(x) = 2^{-1/2}\rho(x,h_i(x))$ on $B$, and observe from elementary geometry that
\[
\{(x,z)\in B\times\CC : \text{$|\RE(z) - h_i(x)| < \delta_i(x)$ and $|\IM(z)| < \delta_i(x)$}\} \subseteq D_i.
\]
Now consider any $j\in\{1,\ldots,k\}$.  Fix functions $\alpha_j,\beta_j\in\S(B)$ such that for all $x\in B$,
\[
h_{j-1}(x) < \alpha_j(x) < \beta_j(x) < h_j(x)
\]
and also
\[
\alpha_j(x) <  h_{j-1}(x) + \delta_{j-1}(x) \,\,\text{and}\,\, h_j(x) - \delta_j(x) < \beta_j(x).
\]
For each $x\in B$, the function $\rho_x$ is continuous on the interval $I_j(x) := [\alpha_j(x),\beta_j(x)]$, so we may define
\[
\epsilon_j(x) = \min\{\rho(x,y) : y\in I_j(x)\}.
\]
The restriction of $F$ to the set
\[
E_j := \{(x,z)\in B\times\CC : \text{$\RE(z)\in \textrm{int}(I_j(x))$ and $|\IM(z)| < \epsilon_j(x)$}\} \subseteq D
\]
is given by the formula
\begin{equation}\label{eq:CKiS:HolomExtStrip:E_j}
F(x,z) = G(x,\RE(z),\i\IM(z)).
\end{equation}
The set $E_j$ is contained in $\{(x,y+z)\in B\times\CC : \text{$(x,y)\in A$ and $|z| < \rho(x,y)$}\}$ for some $A\in\A$, so equation \eqref{eq:CKiS:HolomExtStrip:E_j} shows that $F\restriction{E_j}$ is a continuous $\D$-function that is analytic in $x$ for each choice of $z$.  By construction, we also have that $F$ is analytic in $z$ for each choice of $x$, so $F\restriction{E_j}$ is analytic by Osgood's lemma.

To summarize, we have shown that the restriction of $F$ to each of the individual subanalytic sets $D_i$ and $E_j$ is an analytic $\D$-function.  The family $\{D_i\}_{i=0}^{k} \cup \{E_j\}_{j=1}^{k}$ is an open cover of $\{(x,y)\in B\times\RR : |y|\leq b(x)\}$ in $D\cap(B\times\CC)$, so the restriction of $F$ to $(\bigcup_{i=0}^{k} D_i) \cup(\bigcup_{j=1}^{k}E_j)$ is an analytic $\D$-function.  Define $a:B\to\RR_+$ by setting
\[
a(x) = \min\left(\{\delta_i(x)\}_{i=0}^{k} \cup \{\epsilon_j(x)\}_{j=1}^{k}\right) \,\,\text{for all $x\in B$,}
\]
and observe that $a$ is subanalytic because each function in $\{\delta_i\}_{i=0}^{k} \cup \{\epsilon_j\}_{j=1}^{k}$ is subanalytic. Also observe that the set \eqref{eq:F:RestrSet} defined from this function $a$ is contained in the union $(\bigcup_{i=0}^{k} D_i) \cup(\bigcup_{j=1}^{k}E_j)$, so $F$ further restricts to an analytic $\D$-function on the set \eqref{eq:F:RestrSet}.
\end{proof}

It follows from the proof of Theorem \ref{thm:CKiS:HolomExtStrip} that if let $\D\in\{\S^{\KK},\C^{\KK},\C^{\KK,\i\S}\}$ and consider a function $f\in\D(X\times\RR_{\geq 0})$ with $f_x\in C^\infty(\RR_{\geq 0})$ for all $x\in X$, then there exists a function $a\in\S_+(X)$ and a finite partition $\B$ of $X$ into subanalytic sets such that for each $B\in\B$, the restriction of $f$ to $B\times\RR_{\geq 0}$ extends to an analytic $\D$-function on the set $\{(x,z)\in B\times \CC : \dist(z,\RR_{\geq 0}) < a(x)\}$.  We now combine this observation with Theorem \ref{thm:HolomExtInfty} to extract a simple corollary about extending $f$ analytically to sectors.

\begin{cor}\label{cor:SectorExp}
Let $\D\in\{\S^{\KK},\C^{\KK},\C^{\KK,\i\S}\}$, and assume that $f\in\D(X\times\RR_{\geq 0})$ is such that $f_x\in C^{\infty}(\RR_{\geq 0})$ for all $x\in X$.  Then there exists a subanalytic function $\Theta:X\to(0,\frac{\pi}{2}]$ and a finite partition $\B$ of $X$ into subanalytic sets such that for each $B\in\B$, the restriction of $f$ to $B\times\RR_{\geq 0}$ extends to an analytic $\D$-function on the set
\[
S(\Theta) := \{(x,r\e^{\i\theta}) :  \text{$x\in B$, $r\geq 0$, and $|\theta| < \Theta(x)$}\}.
\]
\end{cor}

\begin{proof}
By Theorem \ref{thm:HolomExtInfty} and the adaptation of Theorem \ref{thm:CKiS:HolomExtStrip} to $X\times\RR_{\geq 0}$ stated directly above, there exist functions $a,b\in\S_+(X)$ and a finite partition $\B$ of $X$ into subanalytic sets such that for each $B\in\B$, the restriction of $f$ to $B\times\RR_{\geq 0}$ extends to an analytic $\D$-function on
\begin{equation}\label{eq:SectorHolomParam}
\textstyle \{(x,z)\in B\times \CC : \text{$\dist(z,\RR_{\geq 0}) < a(x)$, or $|z| > b(x)$ and $|\Arg(z)| < \pi$}\},
\end{equation}
with the understanding that by replacing $b$ with the function $x\mapsto \max\{a(x),b(x)\}$, we may assume that $a(x) \leq b(x)$ for all $x\in X$. 

\begin{figure}[htbp]
 \label{fig:Sector}
 \vskip0.0cm
    \centering
  \hskip0cm
  \includegraphics[scale=0.48]{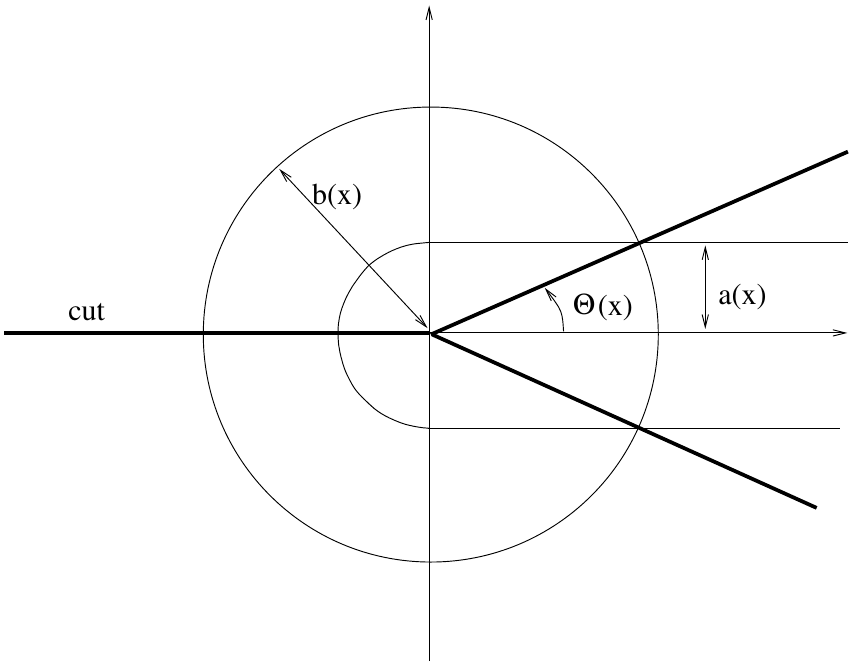}
\caption{}
\end{figure}

Then the function $\Theta:X\to(0,\frac{\pi}{2}]$ defined by
\[
\Theta(x) = \arcsin\left(\frac{a(x)}{b(x)}\right), \quad\text{for all $x\in X$,}
\]
is subanalytic, and $S(\Theta)$ is contained in the set \eqref{eq:SectorHolomParam}.  (See Figure \ref{fig:Sector}.)
\end{proof}

\subsection{Meromorphic Extensions to the Complex Plane}\label{subsection:Meromorphic extension to the complex plane}

Because of the great flexibility afforded by the Mittag-Leffler theorem, if all one knows about a function $f:\RR\to\RR$ is that it is analytic and extends to a meromorphic function on $\CC$, then little can be concluded about the form of its meromorphic extension.  The current subsection is devoted to proving the following theorem, which shows that the situation simplifies radically if we additionally know that $f$ is in $\C^\KK$.

\begin{thm}\label{thm:CK:MeromExt}
Let $\FF\in\{\RR,\CC\}$, $\D\in\{\S^{\KK}_{\FF},\C^{\KK}_{\FF}\}$, and $f\in\D(X\times\RR)$ for some subanalytic set $X\subseteq\RR^m$, and assume that for each $x\in X$, the function $f_x$ is analytic on $\RR$ and extends to a meromorphic function on $\CC$.  Then there exist polynomials $p\in\D(X)[y]$ and $q\in\S(X)[y]$ such that
\begin{equation}\label{eq:CK:MeromExt}
f(x,y) = \frac{p(x,y)}{q(x,y)} \quad\text{for all $(x,y)\in X\times\RR$.}
\end{equation}
\end{thm}

Observe from Example \ref{exam:CKis:notEntirePoly} that Theorem \ref{thm:CK:MeromExt} fails for $\D = \C^{\KK,\i\S}$.  Before proving Theorem \ref{thm:CK:MeromExt}, we remark on its consequences about the poles of $f$, on how it compares to some other theorems in the literature, and why for each $\D\in\{\S^{\KK},\C^{\KK}\}$ we formulate the theorem for both $\D_{\RR}$ and $\D_{\CC}$.

\begin{rems}\label{rem:CK:MeromExt}
Let $\D$, $f$, $p$, and $q$ be as in Theorem \ref{thm:CK:MeromExt}.
\begin{enumerate}
\item\label{rem:CK:MeromExt:SCdenom}
When $\D=\S_{\CC}$, requiring in \eqref{eq:CK:MeromExt} that $p\in\S_{\CC}(X)[y]$ and $q\in\S(X)[y]$ is equivalent to requiring that $p,q\in\S_{\CC}[y]$.

\begin{proof}
If $f= \frac{p}{q}$ for some $p,q\in\S_{\CC}(X)[y]$, then multiplying $\frac{p}{q}$ by $\frac{\,\,\overline{q}\,\,}{\overline{q}}$ expresses $f$ as a quotient with numerator in $\S_{\CC}(X)[y]$ and denominator in $\S(X)[y]$.  The converse is obvious.
\end{proof}

\item\label{rem:CK:MeromExt:partition}
Extend $f$ to $\{(x,z)\in X\times\CC : q(x,z)\neq 0\}$ via the formula \eqref{eq:CK:MeromExt}.  There exists a finite partition $\B$ of $X$ into subanalytic sets such that for each $B\in\B$,
\begin{enumerate}
\item
the polynomials $p$ and $q$ define analytic functions in $\D(B\times\CC)$ and $\S(B\times\CC)$, respectively;

\item
the number of complex zeros of $q$ is constant on $B$ in the following sense: there exist functions $\zeta^{B}_{1},\ldots,\zeta^{B}_{d(B)}\in\S_{\CC}(B)$ and positive integers $\mu^{B}_{1},\ldots,\mu^{B}_{d(B)}$ such that for each $x\in B$, the numbers $\zeta^{B}_{1}(x),\ldots,\zeta^{B}_{d(B)}(x)$ are the distinct complex zeros of $q_x$ and these zeros respectively have multiplicities $\mu^{B}_{1},\ldots,\mu^{B}_{d(B)}$.
\end{enumerate}
The significance of this is that for each $B\in\B$ and $x\in B$, we have
\begin{equation}\label{eq:CK:MeromExt:poleSet}
\{z\in\CC : \text{$z$ is a pole of $f_x$}\} \subseteq \{\zeta^{B}_{1}(x),\ldots,\zeta^{B}_{d(B)}(x)\},
\end{equation}
and for each pole $z$ of $f_x$, say with $z = \zeta_{i}^{B}(x)$, we have
\begin{equation}\label{eq:CK:MeromExt:poleOrder}
\mu(z) \leq \mu_{B}^{i}(x),
\end{equation}
where $\mu(z)$ is the order of the pole $z$ of $f_x$.  When $\D=\S_{\FF}$ for some $\FF\in\{\RR,\CC\}$, we may additionally construct $\B$ so that the number of complex zeros of $p$ is constant on each set in $\B$ in the same sense, and by canceling common linear factors of the numerator and denominator of $\frac{p}{q}$, we may express the quotient $\frac{p}{q}$ in its reduced form, meaning that $p,q\in\S_{\FF}(X)[y]$ and that $\gcd(p_x,q_x) = 1$ for all $x\in X$. This reduction forces equality to hold in \eqref{eq:CK:MeromExt:poleSet} and \eqref{eq:CK:MeromExt:poleOrder}.
\end{enumerate}
\end{rems}

\begin{rems}
Theorem \ref{thm:CK:MeromExt} is somewhat comparable to a number of theorems in the literature.  Perhaps the most classical comparative theorem is due to Borel \cite{Bor}, which states that if a single-variable analytic function $f$ at $0$ has a Maclaurin series with integer coefficients and extends to a meromorphic function on $D_r(0)$ for some $r > 1$, then $f$ is rational.  However, Borel's theorem has nothing to do with tame geometry.

Theorem \ref{thm:CK:MeromExt} more closely resembles Peterzil and Starchenko \cite[Theorem 1.2]{PetStar:2001}, which states that if $\R$ is an o-minimal expansion of a real closed field $R$, $C_0$ is a finite subset of the algebraically closed field $C := R + \i R$, and $f:C\setminus C_0 \to C$ is a $C$-analytic function definable in $\R$, then $f$ is a rational function over $C$.  When $R=\RR_{\an}$ and $\D=\S_{\FF}$ for some $\FF\in\{\RR,\CC\}$, this is similar to Theorem \ref{thm:CK:MeromExt}, except that Theorem \ref{thm:CK:MeromExt} only assumes the definability of $f:X\times\RR\to\FF$ in $\RR_{\an}$, while the definability in $\RR_{\an}$ of its piecewise meromorphic extension is actually a consequence of Theorem \ref{thm:CK:MeromExt}.  Another point of contrast between Theorem \ref{thm:CK:MeromExt} and \cite[Theorem 1.2]{PetStar:2001} is that when $\KK\not\subseteq\RR$, the categories $\S^{\KK}$ and $\C^{\KK}$ contain oscillatory functions that are not definable in any o-minimal structure, as remarked in the Introduction.

Theorem \ref{thm:CK:MeromExt} and \cite[Theorem 1.2]{PetStar:2001} are two theorems among a handful of algebraization results arising out of the study of tame geometry.  For instance, algebraization results similar to Chow's theorem hold, more generally, for analytic sets of $\CC^n$ with monomial growth volume (see Bishop \cite{Bis}, Chirka \cite{Chirka}, Rubel \cite{Ru}, and Stoll \cite{Sto}), or even directly based on o-minimality, a property preventing wild phenomena such as described by the Casorati-Weierstra\ss{}
theorem (see Fortuna and Lojasiewicz \cite{ForLoj} and \cite{PetStar:2001}).
\end{rems}

\begin{rems}
For each $\D\in\{\S^{\KK},\C^{\KK},\C^{\KK,\i\S},\C^{\FTs,\i\S}\}$, every theorem we have discussed so far about preparation, asymptotic expansions, and analytic continuation of $\D$-functions goes through verbatim for $\D_{\CC}$-functions as well, and we shall freely use this fact.  We only limited ourselves to $\D$-functions for simplicity of formulation.  However, $\D_{\CC}$ plays a more significant role in the proof of Theorem \ref{thm:CK:MeromExt}.  Indeed, even if we only cared about establishing Theorem \ref{thm:CK:MeromExt} for $\D$ for each $\D\in\{\S^{\KK},\C^{\KK}\}$, the theorem's proof requires us to first prove the theorem for $\S_{\cl_{\CC}(\KK)}$, which equals $\S^{\QQ}_{\FF}$ for some $\FF\in\{\RR,\CC\}$, so it makes for a more uniform statement of Theorem \ref{thm:CK:MeromExt} to explicitly formulate it for both $\D_{\RR}$ and $\D_{\CC}$.
\end{rems}

We now turn to the proof of Theorem \ref{thm:CK:MeromExt}, which we divide into a number of lemmas.  The first lemma is a version of the theorem without parameters.

\begin{lem}\label{lem:CK:MermExt:NoParam}
If $f\in\C^{\CC}(\RR)$ is analytic and extends to a meromorphic function on $\CC$, then $f$ is a rational function.
\end{lem}

\begin{proof}
Let $f$ be as hypothesized, and also write $f$ for its meromorphic extension to $\CC$.  Theorem \ref{thm:HolomExtInfty} shows that there exists an $r>0$ such that $f\restriction{\CC\setminus D_r(0)}$ is holomorphic.  The function $f$ has finitely many poles $z_1,\ldots,z_k$ in the compact set $D_r(0)$, say with orders $m_1,\ldots,m_k$.  Letting $q(z) = (z-z_1)^{m_1}\cdots(z-z_k)^{m_k}$, the function $qf$ extends to an entire function $p\in\C^{\KK}(\CC)$, which is necessarily a polynomial by Corollary \ref{cor:CK:EntirePoly}. So $f = \frac{p}{q}$, which is rational.
\end{proof}

Let $f$ be as hypothesized in Theorem \ref{thm:CK:MeromExt}.  For each $x\in X$, applying Lemma \ref{lem:CK:MermExt:NoParam} to $f_x$ shows that $f_x = \frac{p_x}{q_x}$ for some polynomials $p_x,q_x\in\CC[y]$ with $\gcd(p_x,q_x) = 1$.  A key step towards proving Theorem \ref{thm:CK:MeromExt} is to give a uniform bound on the degree
\[
\deg(f_x) := \max\{\deg(p_x),\deg(q_x)\}
\]
of the rational function $f_x$ for all $x\in X$.  We do this by using the next lemma, which is a formulation of the well-known fact that a power series in a single variable is a Taylor series
of a rational function if and only if its coefficients satisfy a linear recurrence relation.

\begin{lem}\label{lem:RationalChar}
Let $\FF$ be any field, consider a formal power series $f(y) = \sum_{j=0}^{\infty} f_jy^j$ in $\FF[[y]]$, and let $K\in\NN$.  Then the following three statements are equivalent.
\begin{enumerate}
\item\label{lem:RationalChar:Rational}
The series $f$ is a rational function with $\deg(f)\leq K$.

\item\label{lem:RationalChar:Span}
The set $\{(f_{j},\ldots,f_{j+K}) : j\in\ZZ_+\}$ does not span the vector space $\FF^{K+1}$.

\item\label{lem:RationalChar:det}
For all $J\subseteq\ZZ_+$ with $|J| = K+1$, we have $\det(f_{j+k})_{(j,k)\in J\times\{0,\ldots,K\}} = 0$.
\end{enumerate}
\end{lem}

\begin{proof}
For any polynomials $p(y) = \sum_{j=0}^{K}p_{j}y^{j}$ and $q(y) = \sum_{j=0}^{K} q_{j}y^{j}$ in $\FF[y]$, we have
\[
q(y)f(y)
= \left(\sum_{j=0}^{K}q_{j}y^{j}\right)\left(\sum_{j=0}^{\infty} f_jy^{j}\right)
= \sum_{k=0}^{\infty}\left(\sum_{j=0}^{\min\{k,K\}}q_{j}f_{k-j}\right)y^{k},
\]
so $p = qf$ if and only if
\begin{equation}\label{eq:RationalChar:RecRel}
\text{$\sum_{j=0}^{k}q_{j}f_{k-j} = p_k$ for all $k\leq K$ and  $\sum_{j=0}^{K}q_{j}f_{k-j} = 0$ for all $k > K$.}
\end{equation}
Therefore if \ref{lem:RationalChar:Rational} holds, say with $f=\frac{p}{q}$ for some polynomials $p,q\in\FF[y]$ with $\deg(p),\deg(q)\leq K$ as written above, then \eqref{eq:RationalChar:RecRel} shows that the span of $\{(f_{j},\ldots,f_{j+K}) : j\in\ZZ_+\}$ is contained in the null space of the nonzero linear functional $(x_0,\ldots,x_K) \mapsto q_K x_0 + \cdots + q_0x_K$ on $\FF^{K+1}$, which gives \ref{lem:RationalChar:Span}.  Conversely, if \ref{lem:RationalChar:Span} holds, then there exists a nonzero linear functional $(x_0,\ldots,x_K) \mapsto q_0x_K + \cdots + q_Kx_0$ on $\FF^{K+1}$ whose null space contains the span of $\{(f_{j},\ldots,f_{j+K}) : j\in\ZZ_+\}$, and defining $p_0,\ldots,p_K$ from $q_0,\ldots,q_K$ as in \eqref{eq:RationalChar:RecRel} gives polynomials $p$ and $q$ such that $\deg(p),\deg(q)\leq K$ and $f=\frac{p}{q}$, which proves \ref{lem:RationalChar:Rational}.  To finish, observe that \ref{lem:RationalChar:Span} holds if and only if for every $J\subseteq\ZZ_+$ with $|J| = K+1$, the family of vectors $\{(f_{j},\ldots,f_{j+K})\}_{j\in J}$ is linearly dependent, which is equivalent to \ref{lem:RationalChar:det}.
\end{proof}

\begin{lem}\label{lem:CK:MeromExt:DegBnd}
Let $f\in\C^{\CC}(X\times\RR)$ for a subanalytic set $X\subseteq\RR^m$, and assume that for each $x\in X$, the function $f_x$ is analytic on $\RR$ and extends to a meromorphic function on $\CC$.   Then there exists a $K\in\NN$ such that $\deg(f_x)\leq K$ for all $x\in X$.
\end{lem}

\begin{proof}
For each $x\in X$, the function $f_x:\RR\to\RR$ is analytic, so the equation $f_x(y) = \frac{p_x(y)}{q_x(y)}$ holds on $\RR$ for some polynomials $p_x,q_x\in\CC[y]$ if and only if the same equation holds as germs at $0^+$, so it suffices to study representatives of the germs of the functions $f_x$ at $0^+$.  Each of these germs extend analytically to a neighborhood of $0$, so we may apply Theorem \ref{thm:HolomExtInfty}\ref{thm:HolomExtInfty:Prep} to $f(x,1/y)$ and then pullback by $(x,y) \mapsto (x,1/y)$ to prepare $f$ in $y$ at $0^+$.  Let $B\in\B$, where $\B$ is the partition of $X$ given by Theorem \ref{thm:HolomExtInfty}\ref{thm:HolomExtInfty:Prep}.  Then
\begin{equation}\label{eq:CK:MeromExt:DegBnd:prep}
f(x,y) = \sum_{i\in I} c_i(x) y^{r_i} \tld{f}_i(x,y)
\end{equation}
on $\{(x,y)\in B\times\RR : 0 < y < b(x)\}$, where $b\in\S_+(B)$, $I$ is a finite index, and for each $i\in I$, the function $c_i\in\C^{\CC}(B)$ is analytic, $r_i\in\ZZ$, and $\tld{f}_i = F_i\circ\psi$ is a $\psi$-function with $\psi(x,y) = (\overline{\psi}(x), \frac{y}{b(x)})$, where $\overline{\psi}(x) = (\psi_1(x),\ldots,\psi_N(x))$ and where
\[
F_i(w,z) = \sum_{j=0}^{\infty} F_{i,j}(w)z^j
\]
is an analytic function of $(w,z) = (w_1,\ldots,w_N,z)$ on $[0,1]^N \times D_1(0)$.  It suffices to prove the conclusion of the lemma with $B$ in place of $X$.

By further partitioning $B$, we may assume that $B$ is a cell.  By pulling back by the transformation $(x,y)\mapsto (x,b(x)y)$, we may assume that $b = 1$.  And by multiplying both sides of \eqref{eq:CK:MeromExt:DegBnd:prep} by a suitable integer power of $y$, we may assume that $r_i \geq 0$ for all $i\in I$.  We may therefore absorb each power function $y^{r_i}$ into the $\psi$-function $\tld{f}_i$ to express \eqref{eq:CK:MeromExt:DegBnd:prep} more simply as
\begin{equation}\label{eq:CK:MeromExt:DegBnd:prepSimple}
f(x,y) = \sum_{i\in I} c_i(x) \tld{f}_i(x,y)
\end{equation}
on $B\times(0,1)$.  Using the variables $W = (W_i)_{i\in I}$ on $\CC^{|I|}$, we may write
\[
f(x,y) = G((c_i(x))_{i\in I}, \overline{\psi}(x), y)
\]
on $B\times(0,1)$ for the analytic function $G$ defined by
\[
G(W,w,z) = \sum_{i\in I} W_i F_i(w,z)
\]
on $\CC^{|I|}\times [0,1]^N \times D_1(0)$.  For each $j\in J$, define
\[
G_j(W,w) = \sum_{i\in I}W_iF_{i,j}(w)
\quad\text{and}\quad
f_j(x) = G_j((c_i(x))_{i\in I}, \overline{\psi}(x))
\]
on the sets $\CC^{|I|}\times [0,1]^N$ and $B$, respectively.  Then
\[
G(W,w,z) = \sum_{j=0}^{\infty}G_j(W,w)z^j
\quad\text{and}\quad
f(x,y) = \sum_{j=0}^{\infty}f_j(x)y^{j}
\]
on the sets $\CC^{|I|}\times [0,1]^N\times D_1(0)$ and $B\times(0,1)$, respectively.

For each $K\in\NN$, define
\[
V_K = \{x\in B : \deg(f_x) \leq K\},
\]
and observe from Lemma \ref{lem:RationalChar} that
\[
V_K = \VV(\det(f_{j+k})_{(j,k)\in J\times\{0,\ldots,K\}} : \text{$J\subseteq\ZZ_+$ with $|J| = K+1$})
\]
where we are using the notation $\VV$ to denote the common zero set of a family of functions.  Since each $G_j$ is a polynomial in $W$ with analytic coefficient functions on the compact set $D_1(0)^N$, topological Noetherianity shows that
for each $K\in\NN$,
\begin{align*}
&\VV(\det(G_{j+k})_{(j,k)\in J\times\{0,\ldots,K\}} : \text{$J\subseteq\ZZ_+$ with $|J| = K+1$}) \\
&\quad = \VV(\det(G_{j+k})_{(j,k)\in J\times\{0,\ldots,K\}} : J\in\J_K)
\end{align*}
for some finite set $\J_K \subseteq \{J\subseteq\ZZ_+ : |J| = K+1\}$, so
\[
V_K = \VV(\det(f_{j+k})_{(j,k)\in J\times\{0,\ldots,K\}} : J\in\J_K),
\]
which shows that $V_K$ is an analytic subset of the cell $B$.  We now finish the proof by showing that $V_K = B$ for some $K\in\NN$.

Suppose for a contradiction that $V_K\neq B$ for each $K\in\NN$.  Then for each $K\in\NN$, $V_K$ is a proper analytic subset of the cell $B$ and is therefore closed and nowhere dense in $B$.  So the Baire category theorem supplies an $x\in B\setminus\bigcup_{K\in\NN}V_K$.  But this is impossible because $f_x$ is a rational function by Lemma \ref{lem:CK:MermExt:NoParam}, so $x\in V_K$ for $K = \deg(f_x)$.  This contradiction completes the proof.
\end{proof}

We now divide the proof of Theorem \ref{thm:CK:MeromExt} into three cases according to whether $\D=\S_{\FF}$ (i.e., $\D=\S^{\KK}_{\FF}$ with $\KK=\QQ$), $\D = \S^{\KK}_{\FF}$, or $\D = \C^{\KK}_{\FF}$, where $\FF\in\{\RR,\CC\}$.

\begin{proof}[Proof of Theorem \ref{thm:CK:MeromExt} when $\D = \S_{\FF}$ for some $\FF\in\{\RR,\CC\}$.]
\hfill

Let $f$ be as hypothesized in the theorem.  By Lemma \ref{lem:CK:MeromExt:DegBnd}, we may fix a $K\in\NN$ such that $\deg(f_x)\leq K$ for each $x\in X$.  Therefore, the first-order sentence
\begin{equation}\label{eq:S:MeromExt:Sentence}
\forall x\in X \exists p\exists q\forall y\left((q_Ky^K + \cdots + q_0)f(x,y) = p_Ky^K + \cdots + p_0\right)
\end{equation}
holds in $\RR_{\an}$, where $p=(p_0,\ldots,p_K)$ and $q=(q_0,\ldots,q_K)$ are variables on $\FF^{K+1}$ and $y$ is a variable on $\RR$.  By definable choice for the o-minimal structure $\RR_{\an}$ (see \cite[Sec. 6.1, Prop 1.2]{vdD:Book:1998}), there exist functions $p_0,\ldots,p_K$ and $q_0,\ldots,q_K$ in $\S_{\FF}(X)$ for which
\begin{equation}\label{eq:MeromExt:defChoice}
(q_K(x)y^K+\cdots+q_0(x))f(x,y) = p_K(x)y^K + \cdots + p_0(x)
\end{equation}
holds for all $(x,y)\in X\times\RR$.  This gives us polynomials $p(x,y) = \sum_{j=0}^{K}p_j(x)y^{j}$ and $q(x,y) = \sum_{j=0}^{K}q_j(x)y^{j}$ in $\S_{\FF}(X)[y]$ such that $f = \frac{p}{q}$ on $X\times\RR$.  Keeping in mind Remark \ref{rem:CK:MeromExt}\ref{rem:CK:MeromExt:SCdenom}, this proves the theorem when $\D=\S_{\FF}$.
\end{proof}

The next remarks review some basic concepts about convex polytopes that we use to prove Theorem \ref{thm:CK:MeromExt} when $\D=\S^{\KK}_{\FF}$.

\begin{rems}\label{rem:polytopes}
Consider a convex \DEF{polytope} $P$ contained in a vector space $V$ over an ordered field $\FF$.  By definition, this means that $P = \CONV(G)$ for some finite set $G\subseteq P$, where
\[
\CONV(G) = \left\{\sum_{x\in G} c_x x : \text{$(c_x)_{x\in G} \in\FF_{\geq 0}^{G}$ and $\sum_{x\in G} c_x = 1$}\right\}
\]
is the convex hull of $G$.  We call $G$ a generating set for $P$.  We also call
\[
\VERT(P) = \{x\in P : \text{$P\setminus\{x\}$ is convex}\}
\]
the set of \DEF{vertices} of $P$.
\begin{enumerate}
\item\label{rem:polytopes:VertChar}
There is clearly a minimal subset $M$ of $G$ satisfying $\CONV(M) = P$.  It turns out that $M = \VERT(P)$, as proven below.  It follows that $\VERT(P)$ is necessarily nonempty when $P$ is nonempty.

\begin{proof}
The fact that $M = \VERT(P)$ holds for the ordered field $\RR$ is shown in many sources, such as Br{\o}ndsted \cite[Theorem 7.2]{Bron}.  It follows that $M = \VERT(P)$ must hold for any ordered field $\FF$ because of the quantifier elimination theorem for semilinear sets over ordered fields given in Van den Dries \cite[Chapter 1, Section 7]{vdD}.  (In particular, see \cite[Chapter 1, Remark 7.9]{vdD}.)
\end{proof}

\item\label{rem:polytopes:RemoveVert}
For each $v\in\VERT(P)$, we have $v\not\in\CONV(G\setminus\{v\})$.

\begin{proof}
For each $v\in\VERT(P)$,  the set $G\setminus\{v\}$ is contained in the convex set $P\setminus\{v\}$, so $\CONV(G\setminus\{v\}) \subseteq P\setminus\{v\}$.
\end{proof}

\item\label{rem:polytopes:ConvChar}
The following two statements are equivalent for any finite set $S\subseteq V$ and $v\in V$.
\begin{enumerate}
\item
$v\in\CONV(S)$.

\item
$\sum_{x\in S} c_x(x-v) = 0$ for some nonzero tuple $(c_x)_{x\in S}\in\FF^{S}_{\geq 0}$.
\end{enumerate}

\begin{proof}
When (b) holds, we may divide the tuple $(c_x)_{x\in S}$ by the sum $\sum_{x\in S}c_x$ to additionally assume that $\sum_{x\in S}c_x = 1$.  Thus, the equivalence of (a) and (b) follows from the equivalence of the equations $v = \sum_{x\in S} c_x x$ and $\sum_{x\in S} c_x(x-v) = 0$ when $(c_x)_{x\in S}\in\FF_{\geq 0}^{S}$ satisfies $\sum_{x\in S}c_x = 1$.
\end{proof}
\end{enumerate}
\end{rems}

\begin{lem}\label{lem:ExpansionZero}
Let $f:(0,1)^m\to\CC$ be given as the sum of a convergent series
\begin{equation}\label{eq:ExpansionZero}
f(x) = \sum_{(r,s)\in R\times S} c_{r,s} x^r(\log x)^s
\end{equation}
on $(0,1)^m$, where $R = \Gamma+\NN^m := \{\gamma+\alpha : (\gamma,\alpha)\in\Gamma\times\NN^m\}$ for some finite set $\Gamma\subseteq\KK^m$, the set $S\subseteq\NN^m$ is finite, and $\{c_{r,s}\}_{(r,s)\in R\times S} \subseteq\CC$, and where we are writing $x=(x_1,\ldots,x_m)$ and $\log x = (\log x_1,\ldots,\log x_m)$.  Then $f = 0$ if and only if $c_{r,s} = 0$ for all $(r,s)\in R\times S$.
\end{lem}

\begin{proof}
Clearly $f=0$ if the series is zero, since the series converges to $f$.

Conversely, assume the series is nonzero.  We may assume that $\Gamma = \prod_{i=1}^{m}\Gamma_i$ and $S = \prod_{i=1}^{m}S_i$ for some finite $\Gamma_1,\ldots,\Gamma_m\subseteq\KK$ and $S_1,\ldots,S_m\subseteq\NN^m.$  If $m=1$, then Lemma \ref{lem:noncomp}\ref{lem:noncomp:Seq} shows that $f$ is asymptotic with its leading form along some sequence in $(0,1)$ converging to $0$, so $f\neq 0$.  So suppose that $m>1$ and that the lemma holds with $m-1$ in place of $m$.  Write $f$ as a series in $x_m$ with coefficient functions in $(x_1,\ldots,x_{m-1})$.  The induction hypothesis supplies an $a\in(0,1)^{m-1}$ such that some coefficient function of $f$ is nonzero at $a$, so $f(a,x_m)$ is the sum of a nonzero series in $x_m$ and is therefore a nonzero function by the base case of the induction.
\end{proof}

Call a series of the form \eqref{eq:ExpansionZero} a \DEF{$\KK$-power-log series}, and call $\Supp(f) = \{(r,s)\in R\times S: c_{r,s}\neq 0\}$ the \DEF{support} of the series.  If $S = \{0\}$, then we may more simply write $f(x) = \sum_{r\in R}c_r x^r$ and $\Supp(f) = \{r\in R : c_R\neq 0\}$, and we call this a \DEF{$\KK$-power series}.  When $\Gamma\subseteq\QQ^m$, we instead say $\QQ$-power-log series or $\QQ$-power series.

In the remaining two cases of the proof of Theorem \ref{thm:CK:MeromExt}, let $\FF$ denote the topological closure of $\KK$ in $\CC$, and let $f$ be as hypothesized in the theorem.

\begin{proof}[Proof of Theorem \ref{thm:CK:MeromExt} when $\D = \S^{\KK}_{\FF}$ for some $\FF\in\{\RR,\CC\}$]
By proceeding as in the proof of Lemma \ref{lem:CK:MeromExt:DegBnd}, we may reduce to studying
\[
f(x,y) = \sum_{i\in I} c_i(x) \tld{f}_i(x,y)
\]
on $B\times(0,1)$, as in equation \eqref{eq:CK:MeromExt:DegBnd:prepSimple} except with $\{c_i\}_{i\in I} \subseteq\S^{\KK}_{\FF}(B)$.  Apply Proposition \ref{prop:SubRect} to rectilinearize the subanalytic functions $\psi_1,\ldots,\psi_N$ from the proof of Lemma \ref{lem:CK:MeromExt:DegBnd} and all of the subanalytic functions used in the formulas defining the functions in $\{c_i\}_{i\in I}$.  Thus, up to further partitioning $B$ and pulling back by a bi-analytic subanalytic map, and potentially reducing the value of $m$ accordingly, we reduce to the situation where $B = (0,1)^m$ and
\begin{equation}\label{eq:SK:MeromExt:rect}
f(x,y) = \sum_{\alpha\in A} x^\alpha f_\alpha\left(x,y\right)
\end{equation}
on $(0,1)^{m+1}$, where $A$ is a finite subset of $\KK^m$ and each $f_\alpha$ is an analytic function on $[0,1]^m\times D_1(0)$ with $f_\alpha([0,1]^{m+1}) \subseteq\FF$.  If $\alpha,\beta\in A$ are such that $\beta - \alpha = \delta$ for some nonzero $\delta\in\QQ^m$, then
\[
x^\alpha f_\alpha + x^\beta f_\beta = x^\alpha(f_\alpha + x^\delta f_\beta)
\]
and $f_\alpha + x^\delta f_\beta$ is an $\S_{\FF}$-function represented as the sum of a convergent power series in $y$ with $\QQ$-power series coefficient functions of $x$ on $(0,1)^m$.  By repeating such computations and replacing $A$ with a subset of itself, we may assume that \eqref{eq:SK:MeromExt:rect} holds, except where each $f_\alpha$ is now an analytic function in $\S_{\FF}((0,1)^m\times D_1(0))$ that is the sum of a power series
\[
f_\alpha(x,y) = \sum_{j=0}^{\infty} f_{\alpha,j}(x)y^j
\]
on $(0,1)^{m+1}$ with coefficient functions $\{f_{\alpha,j}\}_{j\in\NN}$ in $\S_{\FF}((0,1)^m)$ that are sums of convergent $\QQ$-power series, and where for all $\alpha,\beta\in A$ with $\alpha\neq\beta$, we have $\alpha + \QQ^m \neq \beta + \QQ^m$ in the $\QQ$-vector space $\KK^m/\QQ^m$.  Write
\[
f(x,y) = \sum_{j=0}^{\infty}f_j(x)y^j
\]
on $(0,1)^{m+1}$, where for each $j\in\NN$,
\[
f_j(x) = \sum_{\alpha\in A} x^\alpha f_{\alpha,j}(x)
\]
on $(0,1)^m$.
\begin{quote}
{\scshape Claim.}  The function $(f_{\alpha})_{x} : y\mapsto f_\alpha(x,y)$ is rational for each $\alpha\in A$ and $x\in(0,1)^m$.
\end{quote}
The claim and Theorem \ref{thm:CK:MeromExt} for $\S_{\FF}$ together show that each function $f_\alpha$ may be expressed as a quotient of polynomials with numerator in $\S_{\FF}((0,1)^m)[y]$ and denominator in $\S((0,1)^m)[y]$.  The theorem for $\D=\S^{\KK}_{\FF}$ follows from this and equation \eqref{eq:SK:MeromExt:rect}, so it suffices to prove the claim.

Suppose for a contradiction that the claim is false.  For each $\alpha\in A$ for which $(f_{\alpha})_x$ is rational for all $x\in(0,1)^m$, subtract the term $f_\alpha$ from both sides of the equation \eqref{eq:SK:MeromExt:rect}.  By doing this, we may assume in \eqref{eq:SK:MeromExt:rect} that for each $\alpha\in A$, there exists an $x\in(0,1)^m$ such that $(f_{\alpha})_{x}$ is not rational.  Even after this reduction, the set $A$ is necessarily nonempty because of the assumed falsity of the claim.  By Lemma \ref{lem:CK:MeromExt:DegBnd}, we may now fix a $K\in\NN$ such that $\deg(f_x)\leq K$ for each $x\in X$.  By Lemma \ref{lem:RationalChar}, the fact that $f_x$ is rational of degree at most $K$ for all $x\in(0,1)^m$ shows that
\begin{equation}\label{eq:SK:f:notSpan}
\SPAN\{(f_j(x),\ldots,f_{j+K}(x))\}_{j\in\ZZ_+} \neq \FF^{K+1} \,\,\text{for all $x\in(0,1)^m$,}
\end{equation}
and for each $\alpha\in A$, the fact that $(f_\alpha)_x$ is not rational for some $x\in(0,1)^m$ shows that
\begin{equation}\label{eq:SK:falpha:Span}
\SPAN\{(f_{\alpha,j}(x),\ldots,f_{\alpha,j+K}(x))\}_{j\in\ZZ_+} = \FF^{K+1} \,\,\text{for some $x\in(0,1)^m$.}
\end{equation}
We will deduce a contradiction from \eqref{eq:SK:f:notSpan} and \eqref{eq:SK:falpha:Span} via a determinant computation.

By Remark \ref{rem:polytopes}\ref{rem:polytopes:VertChar}, we may choose a $\mu\in A$ such that $\mu+\QQ^m$ is a vertex of the convex hull of $\{\alpha+\QQ^m : \alpha\in A\}$ in the $\QQ$-vector space $\KK^m/\QQ^m$.  By dividing both sides of \eqref{eq:SK:MeromExt:rect} by $x^\mu$ and adjusting our notation accordingly, we may assume in \eqref{eq:SK:MeromExt:rect} that $\mu=0$.  Therefore by Remarks \ref{rem:polytopes}\ref{rem:polytopes:RemoveVert} and \ref{rem:polytopes}\ref{rem:polytopes:ConvChar},
\begin{equation}\label{eq:A-0}
\sum_{\alpha\in A\setminus\{0\}} c_\alpha \alpha \not\in \QQ^m \,\,\text{for every nonzero}\,\, (c_\alpha)_{\alpha\in A\setminus\{0\}} \in \QQ_{\geq 0}^{|A\setminus\{0\}|}.
\end{equation}
Applying \eqref{eq:SK:falpha:Span} to $\alpha = 0$ shows that we may fix some $J\subseteq\ZZ_+$ such that $|J| = K+1$ and
\begin{equation}\label{eq:SK:fmu:detNotZero}
\det(f_{0,j+k}(x))_{(j,k)\in J\times\{0,\ldots,K\}} \neq 0 \quad\text{for some $x\in (0,1)^m$.}
\end{equation}
However, for the function $g:(0,1)^m\to\CC$ defined by
\[
g(x) = \det(f_{j+k}(x))_{(j,k)\in J\times\{0,\ldots,K\}}
\]
on $(0,1)^m$, we know from \eqref{eq:SK:f:notSpan} that
\begin{equation}\label{eq:SK:g:detZero}
g(x) = 0 \quad\text{for all $x\in (0,1)^m$.}
\end{equation}

On the other hand, for each $x\in(0,1)^m$,
\begin{align*}
g(x)
    &= \det\left(\sum_{\alpha\in A} x^\alpha f_{\alpha,j+k}(x)\right)_{(j,k)\in J\times\{0,\ldots,K\}} \\
    &= \sum_{\gamma\in A^{K+1}} x^{|\gamma|} \det(f_{\gamma_k,j+k}(x))_{(j,k)\in J\times\{0,\ldots,K\}},
\end{align*}
where for each $\gamma\in A^{K+1}$, we are writing $\gamma = (\gamma_0,\ldots,\gamma_K)$ and $|\gamma| = \sum_{k=0}^{K}\gamma_k \in \KK^m$ with $\gamma_0,\ldots,\gamma_K\in A$.  Combine like terms in the monomials $x^{|\gamma|}$ by writing
\begin{equation}\label{eq:g:det}
g(x) = \sum_{\delta\in\Delta} x^\delta g_\delta(x),
\end{equation}
where $\Delta = \{|\gamma| : \gamma\in A^{K+1}\}$ and for each $\delta\in\Delta$,
\[
g_\delta(x) = \sum_{\scriptstyle \substack{\gamma\in A^{K+1} \\ \scriptstyle\text{s.t.}\,\,|\gamma|=\delta}} \det(f_{\gamma_k,j+k}(x))_{(j,k)\in J\times\{0,\ldots,K\}}.
\]
Each function $x^\delta g_\delta$ is the sum of a $\KK$-power series.  From \eqref{eq:A-0} we see that $\Supp(x^\delta g_\delta) \subseteq \KK^m\setminus\QQ^m$ for all $\delta\in\Delta$ with $\delta\neq 0$ and that
\[
x^0 g_0(x) = \det(f_{0,j+k}(x))_{(j,k)\in J\times\{0,\ldots,K\}},
\]
which by \eqref{eq:SK:fmu:detNotZero} is the sum of a nonzero series with $\Supp(x^0 g_0) \subseteq \QQ^m$.  So $g$ is the sum of a nonzero $\KK$-power series.  This implies that $g\neq 0$ by Lemma \ref{lem:ExpansionZero}, but that is impossible by \eqref{eq:SK:g:detZero}.  This contradiction proves the claim and thereby establishes the theorem when $\D = \S^{\KK}_{\FF}$.
\end{proof}

\begin{proof}[Proof of Theorem \ref{thm:CK:MeromExt} when $\D = \C^{\KK}_{\FF}$]
The proof for $\C^{\KK}_{\FF}$ is similar in many ways to the proof for $\S^{\KK}_{\FF}$.  Indeed, we may reduce to studying
\[
f(x,y) = \sum_{i\in I} c_i(x) \tld{f}_i(x,y)
\]
on $B\times(0,1)$, except now $\{c_i\}_{i\in I}\subseteq\C^{\KK}_{\FF}(B)$.  Apply Proposition \ref{prop:SubRect} to rectilinearize the subanalytic functions $\psi_1,\ldots,\psi_N$ and all of the subanalytic functions used in the formulas defining the functions in $\{c_i\}_{i\in I}$.  Thus, up to further partitioning $B$ and pulling back by a bi-analytic subanalytic map, and potentially reducing the value of $m$ accordingly, we reduce to the situation where $B = (0,1)^m$.   Moreover, by expanding the logarithms and distributing the powers (as done in the proof of Proposition \ref{prop:CKprep}), distributing multiplication across addition, and combining like terms in the powers of the logarithms, we may write
\begin{equation}\label{eq:CK:MeromExt:rect}
f(x,y) = \sum_{\alpha\in A} (\log x)^\alpha f_\alpha(x,y)
\end{equation}
on $(0,1)^{m+1}$, where $A$ is a finite subset of $\NN^m$ and each $f_\alpha$ is an analytic function in $\S^{\KK}_{\FF}((0,1)^m\times D_1(0))$ given as the sum of a convergent power series
\[
f_\alpha(x,y) = \sum_{j=0}^{\infty} f_{\alpha,j}(x)y^j
\]
on $(0,1)^{m+1}$ with coefficient functions $\{f_{\alpha,j}\}_{j\in\NN}$ in $\S_{\FF}((0,1)^m)$ that are sums of convergent $\KK$-power series.  Write
\[
f(x,y) = \sum_{j=0}^{\infty}f_j(x)y^j
\]
on $(0,1)^{m+1}$, where for each $j\in\NN$,
\[
f_j(x) = \sum_{\alpha\in A} (\log x)^\alpha f_{\alpha,j}(x)
\]
on $(0,1)^m$.
\begin{quote}
{\scshape Claim.}  The function $(f_{\alpha})_{x}:y\mapsto f_\alpha(x,y)$ is rational for each $\alpha\in A$ and $x\in(0,1)^m$.
\end{quote}
The claim and Theorem \ref{thm:CK:MeromExt} for $\S^{\KK}_{\FF}$ together show that each function $f_\alpha$ may be expressed as a quotient of polynomials with numerator in $\S^{\KK}_{\FF}((0,1)^m)[y]$ and denominator in $\S((0,1)^m)[y]$.  The theorem for $\D=\C^{\KK}_{\FF}$ follows from this and equation \eqref{eq:CK:MeromExt:rect}, so it suffices the prove the claim.

Suppose for a contradiction that the claim is false.  For each $\alpha\in A$ for which $(f_\alpha)_x$ is rational for all $x\in(0,1)^m$, subtract the term $(\log x)^\alpha f_\alpha(x,y)$ from both sides of the equation \eqref{eq:CK:MeromExt:rect}.  By doing this, we may assume in \eqref{eq:CK:MeromExt:rect} that for each $\alpha\in A$, there exists an $x\in(0,1)^m$ such that $(f_\alpha)_x$ is not rational.  Even after this reduction, the set $A$ is necessarily nonempty because of the assumed falsity of the claim.  By Lemma \ref{lem:CK:MeromExt:DegBnd}, we may now fix a $K\in\NN$ such that $\deg(f_x)\leq K$ for each $x\in X$.  Thus,
\begin{equation}\label{eq:CK:f:notSpan}
\SPAN\{(f_j(x),\ldots,f_{j+K}(x))\}_{j\in\ZZ_+} \neq \FF^{K+1} \,\,\text{for all $x\in(0,1)^m$,}
\end{equation}
and for each $\alpha\in A$,
\begin{equation}\label{eq:CK:falpha:Span}
\SPAN\{(f_{\alpha,j}(x),\ldots,f_{\alpha,j+K}(x))\}_{j\in\ZZ_+} = \FF^{K+1} \,\,\text{for some $x\in(0,1)^m$,}
\end{equation}
from which we will deduce a contradiction via a determinant computation.

To that end, let $\mu$ be the lexicographically maximum element of $A$.  By \eqref{eq:CK:falpha:Span} we may fix some $J\subseteq\ZZ_+$ such that $|J| = K+1$ and
\begin{equation}\label{eq:Fmu:detNotZero}
\det(f_{\mu,j+k}(x))_{(j,k)\in J\times\{0,\ldots,K\}} \neq 0 \quad\text{for some $x\in (0,1)^m$.}
\end{equation}
However, for the function $g:(0,1)^m\to\CC$ defined by
\[
g(x) = \det(f_{j+k}(x))_{(j,k)\in J\times\{0,\ldots,K\}}
\]
on $(0,1)^m$, we know from \eqref{eq:CK:f:notSpan} that
\begin{equation}\label{eq:CK:g:detZero}
g(x) = 0 \quad\text{for all $x\in (0,1)^m$.}
\end{equation}
On the other hand, on $(0,1)^m$ we have
\begin{align*}
g(x)
    &= \det\left(\sum_{\alpha\in A} (\log x)^\alpha f_{\alpha,j+k}(x)\right)_{(j,k)\in J\times\{0,\ldots,K\}} \\
    &= \sum_{\gamma\in A^{K+1}} (\log x)^{|\gamma|} \det(f_{\gamma_k,j+k}(x))_{(j,k)\in J\times\{0,\ldots,K\}}.
\end{align*}
Combine like terms in the powers of the logarithms by writing
\begin{equation}\label{eq:g:DetLog}
g(x) = \sum_{\delta\in\Delta} (\log x)^\delta G_\delta(x),
\end{equation}
where $\Delta = \{|\gamma| : \gamma\in A^{K+1}\}$ and for each $\delta\in\Delta$,
\[
G_\delta(x) = \sum_{\scriptstyle \substack{\gamma\in A^{K+1} \\ \scriptstyle\text{s.t.}\,\,|\gamma|=\delta}} \det(f_{\gamma_k,j+k}(x))_{(j,k)\in J\times\{0,\ldots,K\}}.
\]
It follows from the choice of $\mu$ that $\gamma\in A^{K+1}$ satisfies $|\gamma| = (K+1)\mu$ if and only if $\gamma_0=\cdots=\gamma_K = \mu$, so $G_{(K+1)\mu}\neq 0$ by \eqref{eq:Fmu:detNotZero}.  Each $G_\delta$ is the sum of a $\KK$-power series, so equation \eqref{eq:g:DetLog} expresses $g$ as the sum of a nonzero $\KK$-power-log series.  Therefore $g\neq 0$ by Lemma \ref{lem:ExpansionZero}, which is impossible by \eqref{eq:CK:g:detZero}.  This contradiction proves the claim and so also the theorem when $\D=\C^{\KK}_{\FF}$.
\end{proof}

This completes the proof of Theorem \ref{thm:CK:MeromExt}.

\begin{exam}\label{exam:a non subanalytic series}
Consider the function $f:\RR\to\RR$ defined by
\[
f(y) = \sum_{n=1}^{+\infty} \frac{1}{(y-n)^2+1/n^2},\,\,\text{for each $y\in\RR$,}
\]
and let $P = \{n\pm\i/n : n\in\ZZ_+\}$.  This series converges locally uniformly on $\CC\setminus P$, so $f$ is an analytic function on $\RR$ that extends meromorphically to $\CC$ with set of poles $P$.  We have that $f\not\in\S(\RR)$, which can be checked by a direct computation showing that the graph of $f$ encounters a horizontal line of $\RR^2$ infinitely many times; indeed, $f(n)\geq n^2$ for all $n\in\ZZ_+$, whereas $f(n+1/2)$ is uniformly bounded from above (for instance, by $8+2\zeta(2)$).  However, we can also use our theorems to conclude even more about $f$ without any computation.  Indeed, since $f$ is meromorphic on $\CC$ but is not rational, we must have that $f\not\in\C^{\CC}(\RR)$ by Theorem \ref{thm:CK:MeromExt} (or more simply, by Lemma \ref{lem:CK:MermExt:NoParam}).  And even more, since $f$ is analytic on $\RR$ but any open horizontal strip about $\RR$ in $\CC$ intersects $P$, we in fact have that $f\not\in\C^{\CC,\i\S}(\RR)$ by Theorem \ref{thm:CKiS:HolomExtStrip}.
\end{exam}

\section{Decay of Fourier Transforms of \texorpdfstring{$\C^{\KK}$} --Functions}\label{s:DecayFT}

Given an integrable function $f\in\C^{\KK}(\RR)$, this section discusses two main theorems that relate the size of a suitably shaped domain in the complex plane to which $f$ may be extended analytically (when possible) with the exponential rate of decay of a corresponding integral transform of $f$.  The first theorem, Theorem \ref{thm:CK:FTexpDecay}, discusses analytic extensions to horizontal strips about the real axis and the Fourier transform, while the second theorem, Theorem \ref{thm:CK:MellinexpDecay}, discusses analytic extensions to sectors about the positive real axis and the Mellin transform.  By pulling back $f$ by the exponential function, the second theorem is easily obtained from two key lemmas about Fourier transforms that are used to prove the first theorem, so the section largely focuses on Fourier transforms.  The section is divided into three subsections discussing (1) the first main theorem, (2) the second main theorem, and (3) examples and remarks elaborating upon various issues related to the two main theorems.

\subsection{Analytic Extensions on Strips and the Fourier Transform}\label{ss:DecayFT:Strips}

To motivate the first main theorem of the section, we begin with the following simple proposition whose content is classical.

\begin{prop}\label{prop:FTmotivate}
Let $f\in\L^1(\RR)$.
\begin{enumerate}
\item\label{prop:FTmotivate:FTdiff}
Let $\sigma\in\{-1,1\}$, define $f_k(t) = (-\sigma\i 2\pi t)^kf(t)$ for all $t\in\RR$ and $k\in\NN$, and let $n\in\NN$ with $f_n\in \L^1(\RR)$.  Then the function
\[
\FT^{\sigma}[f](\xi) = \int_{\RR}f(t)\e^{-\sigma\i 2\pi\xi t}\d t
\]
is in $C^n(\RR)$, and $(\FT^{\sigma}[f])^{(k)} = \FT^{\sigma}[f_k]$ for all $k\in\{0,\ldots,n\}$.

\item\label{prop:FTmotivate:FTfastDecay->Smooth}
If $\FT[f]$ has rapid decay at $\pm\infty$, then $[f]\in C^{\infty}(\RR)$.

\item\label{prop:FTmotivate:Smooth->FTfastDecay}
If $f\in C^{\infty}(\RR)$ and $f^{(n)}\in\L^1(\RR)$ for all $n\in\NN$, then $\FT[f]$ has rapid decay at $\pm\infty$.
\end{enumerate}
\end{prop}

\begin{proof}
For statement \ref{prop:FTmotivate:FTdiff}, see (for instance) Gasquet and Witomski \cite[Proposition 17.2.1]{GasWit}.  For statement \ref{prop:FTmotivate:FTfastDecay->Smooth}, note that if $\FT[f]$ has rapid decay at $\pm\infty$, then for each $n\in\NN$, the functions $\xi\mapsto \FT[f](\xi)$ and $\xi\mapsto (\i2\pi \xi)^n \FT[f](\xi)$ are in $\L^1(\RR)$, so $[f] = \FT^{-1}\circ \FT[f]\in C^n(\RR)$ by the Fourier inversion theorem and statement \ref{prop:FTmotivate:FTdiff}.  Statement \ref{prop:FTmotivate:Smooth->FTfastDecay} follows from repeated integration by parts.
\end{proof}

The next corollary is an immediate consequence of Lemma \ref{lem:CK:L1}\ref{lem:CK:L1:Cn}, Corollary \ref{cor:aeClass}, Theorem \ref{thm:CKiS:HolomExtStrip}, and Proposition \ref{prop:FTmotivate}.

\begin{cor}\label{cor:FTmotivate}
Let $f\in\C^{\KK}(\RR)\cap\L^1(\RR)$.  Up to modifying $f$ at finitely many points, the following three statements are equivalent.
\begin{enumerate}
\item
$\FT[f]$ decays rapidly at $\pm\infty$.

\item
$f\in C^{\infty}(\RR)$.

\item
$f$ is analytic.
\end{enumerate}
\end{cor}

Our first main theorem improves upon Corollary \ref{cor:FTmotivate} by showing that the three conditions of the corollary are also equivalent to $\FT[f]$ decaying exponentially at $\pm\infty$ and by showing that the rate of exponential decay of $\FT[f]$ completely characterizes the maximal horizontal strip about the real axis to which the function $f$ extends analytically.

\begin{defs}\label{def:deltas}
Consider any function $f:\RR\to\CC$ and sign $\sigma\in\{+,-\}$.  When $[f]$ is analytic, define
\[
\delta^{H}_{\sigma}(f) = \sup\{\alpha\in\RR_{\geq 0} : \text{$[f]$ extends analytically to $\RR+\sigma\i[0,\alpha]$}\},
\]
and define $\delta^{H}_{\sigma}(f) = 0$ otherwise.  When $f\in \L^1(\RR)$,  define
\[
\delta^{\FTs}_{\sigma}(f) = \sup\{\alpha\in\RR_{\geq 0}: \text{$\FT[f](\xi) = O(\e^{\sigma 2\pi\alpha\xi})$ as $\xi\to -\sigma\infty$}\},
\]
These suprema are in the sense of the extended real line, so $\delta^{H}_{\sigma}(f)$  and $\delta^{\FTs}_{\sigma}(f)$ take values in $[0,+\infty]$.  (The superscripts $H$ and $\FT$ in this notation stand for \emph{holomorphic} and \emph{Fourier transform}.)
\end{defs}

If $f$ is real-valued, then
\begin{equation}\label{eq:deltaReal}
\text{$\delta_{+}^{H}(f) = \delta_{-}^{H}(f)$ and $\delta_{+}^{\FTs}(f) = \delta_{-}^{\FTs}(f)$}
\end{equation}
because of the Schwartz reflection principle and because $\overline{\FT[f](\xi)} = \FT[f](-\xi)$ for all $\xi\in\RR$.  However, \eqref{eq:deltaReal} does not hold in general when $f$ is complex-valued.

\begin{thm}\label{thm:CK:FTexpDecay}
Let $f\in\C^{\KK}(\RR)\cap \L^1(\RR)$.  Up to modifying $f$ at finitely many points, the following four statements are equivalent.
\begin{enumerate}
\item \label{thm:CK:FTdecay:1}
$\FT[f]$ decays exponentially at $\pm\infty$.

\item \label{thm:CK:FTdecay:2}
$\FT[f]$ decays rapidly at $\pm\infty$.

\item \label{thm:CK:FTdecay:3}
$f\in C^{\infty}(\RR)$.

\item \label{thm:CK:FTdecay:4}
$f$ is analytic.
\end{enumerate}
In addition, $\delta^{\FTs}_{\sigma}(f) = \delta^{H}_{\sigma}(f)$ for each $\sigma\in\{+,-\}$.
\end{thm}

Before proving the theorem, we state a simple corollary for parametric families of $\C^{\KK}$-functions.

\begin{cor}\label{cor:CKFTdecay}
Let $f\in\C^{\KK}(X\times\RR)$ for some subanalytic set $X\subseteq\RR^m$, and assume that $f_x\in \L^1(\RR)\cap C^{\infty}(\RR)$ for all $x\in X$.  Then there exists a function $b\in\S_+(X)$ such that $\FT[f_x](\xi) = o(\e^{-2\pi b(x)|\xi|})$ as $|\xi|\to+\infty$.
\end{cor}

\begin{proof}
Theorems \ref{thm:CKiS:HolomExtStrip} and \ref{thm:CK:FTexpDecay} show that there exists a function $b\in\S_+(X)$ such that $\delta_{\sigma}^{\FTs}(f_x) = \delta_{\sigma}^{H}(f_x) > b(x)$ for all $x\in X$ and $\sigma\in\{+,-\}$.
\end{proof}

In this section we work with single variable $\C^{\KK}$-functions rather than with parametric families of such functions because we do not presently know (nor do we have any compelling reason to expect) that $x\mapsto \delta_{\sigma}^{H}(f_x)$ is a subanalytic function on $X$ when $f\in\C^{\KK}(X\times\RR)$.  We now turn to the proof of Theorem \ref{thm:CK:FTexpDecay}, which is based on the next two lemmas.

\begin{lem}\label{lem:FTDecayImplyHolomExt}
Let $\alpha_+,\alpha_-\in\RR$ with $\alpha_- < 0 < \alpha_+$, and consider a function $f\in\L^1(\RR)$ such that $\FT[f](\xi) = O(\e^{2\pi \alpha_\sigma\xi})$ as $\xi\to-\sigma\infty$ for each $\sigma\in\{+,-\}$.  Then $[f]$ extends to a holomorphic function on the open horizontal strip $\RR+\i(\alpha_-,\alpha_+)$.
\end{lem}

\begin{proof}
Let $f$ be as hypothesized.  Then $\FT[f]\in \L^1(\RR)$, so up to replacing $f$ with the continuous function in $[f]$, we have
\[
f(t) = \int_{\RR}\FT[f](\xi)\e^{\i 2\pi \xi t}\d\xi, \quad\text{for all $t\in\RR$,}
\]
by the Fourier inversion theorem.  There exist constants $C,R>0$ such that for each $z = x+\i y\in\CC$, we have
\[
\left|\FT[f](\xi)\e^{\i 2\pi \xi z}\right| = \left|\FT[f](\xi)\right|\e^{-2\pi\xi y} \leq
\begin{cases}
C\e^{2\pi\xi(\alpha_- -y)}, & \quad\text{for all $\xi > R$,} \\
C\e^{2\pi\xi(\alpha_+-y)},  & \quad\text{for all $\xi < - R$.}
\end{cases}
\]
When $\alpha_- < y < \alpha_+$, these two bounding functions are integrable at $+\infty$ and $-\infty$, respectively, so by holomorphy under the integral sign,
\[
z\mapsto \int_{\RR}\FT[f](\xi)\e^{\i 2\pi \xi z}\d\xi
\]
is a holomorphic extension of $f$ to $\RR+\i(\alpha_-,\alpha_+)$.
\end{proof}

\begin{lem}\label{lem:HolExtImplyFTDecay}
Consider a constant $\alpha>0$, a sign $\sigma\in\{+,-\}$, a finite subset $\Omega$ of the open horizontal strip $\RR+\sigma\i(0,\alpha)$, and a holomorphic function $f:(\RR+\sigma\i[0,\alpha])\setminus\Omega \to\CC$ such that the functions $t\mapsto f(t)$ and $t\mapsto f(t+\sigma\i\alpha)$ are in $\L^1(\RR)$ and
\begin{equation}\label{eq:FTexpDecay:supLimit}
\lim_{t\to\pm\infty}\sup_{\sigma y\in[0,\alpha]}|f(t+iy)| = 0.
\end{equation}
Then there exists a function $\R[f]:\RR\to\CC$ such that for each $\xi\in\RR$,
\begin{equation}\label{eq:FTexpDecay:FTbnd}
\FT[f](\xi) = \sigma\i 2\pi\left(\sum_{\omega\in\Omega} \Res{z=\omega} f(z)\e^{-\i 2\pi\xi z}\right) + \R[f](\xi)
\end{equation}
and
\begin{equation}\label{eq:FTexpDecay:Rbnd}
|\R[f](\xi)| \leq \e^{\sigma 2\pi\alpha\xi}\int_{\RR}|f(t+\i\sigma \alpha)|\d t.
\end{equation}
\end{lem}

\begin{proof}
Let $f$ be as hypothesized, and consider any $\xi\in\RR$ and any $R > 0$ that is sufficiently large so that $\Omega \subseteq (-R,R)+\sigma\i(0,\alpha)$.  By the residue theorem, the integral of the function $z\mapsto f(z)\e^{-\i2\pi\xi z}$ along the contour $[-R,R]$ (as in Notation \ref{notn:complex}) is equal to its integral along the contour
\[
[-R,-R+\sigma\i\alpha]+[-R+\sigma\i\alpha,R+\sigma\i\alpha]+[R+\sigma\i\alpha,R]
\]
plus $\sigma\i 2\pi$ times the sum of its residues over all points of $\Omega$.  Letting $R\to+\infty$, the integrals over $[-R,-R+\sigma\i\alpha]$ and $[R+\sigma\i\alpha,R]$ tend to zero by \eqref{eq:FTexpDecay:supLimit}, so \eqref{eq:FTexpDecay:FTbnd} holds for the function $\R[f]:\RR\to\CC$ given by
\[
\R[f](\xi)
=
\int_{\RR} f(t+\sigma\i\alpha)\e^{-\i2\pi\xi(t+\sigma\i\alpha)} \d t
= \e^{\sigma 2\pi\alpha\xi} \int_{\RR} f(t+\sigma \i\alpha)\e^{-\i 2\pi\xi t} \d t
\]
for each $\xi\in\RR$.  Since $t\mapsto f(t+\sigma\i\alpha)$ is in $\L^1(\RR)$, this gives \eqref{eq:FTexpDecay:Rbnd}.
\end{proof}

The proof of Theorem \ref{thm:CK:FTexpDecay} only uses Lemma \ref{lem:HolExtImplyFTDecay} when the set $\Omega$ is empty, in which case $\FT[f] = \R[f]$.  Allowing a potentially nonempty set $\Omega$ of singular points in the formulation of Lemma \ref{lem:HolExtImplyFTDecay} is useful for obtaining a slightly more explicit bound on the decay of $\FT[f]$ when $f$ extends meromorphically to a horizontal strip, as we discuss later in the Remarks \ref{rems:DecayFT:Sing}.

\begin{proof}[Proof of Theorem \ref{thm:CK:FTexpDecay}]
Consider the following two statements:
\begin{enumerate}[label=($\arabic*'$)]\setcounter{enumi}{2}
\item \label{thm:CK:FTdecay:3'}
$[f]\in\C^{\infty}(\RR)$.

\item \label{thm:CK:FTdecay:4'}
$[f]$ is analytic.
\end{enumerate}
By Corollary \ref{cor:aeClass}, our job is to show that the statements \ref{thm:CK:FTdecay:1}, \ref{thm:CK:FTdecay:2}, \ref{thm:CK:FTdecay:3'}, and \ref{thm:CK:FTdecay:4'}  are equivalent.  The implication \ref{thm:CK:FTdecay:1}$\Rightarrow$\ref{thm:CK:FTdecay:2} is trivial; the statements \ref{thm:CK:FTdecay:2}, \ref{thm:CK:FTdecay:3'}, and \ref{thm:CK:FTdecay:4'} are equivalent by Corollary \ref{cor:FTmotivate}; Lemma \ref{lem:FTDecayImplyHolomExt} shows that $\delta_{\sigma}^{H}(f) \geq \delta_{\sigma}^{\FTs}(f)$ for each $\sigma\in\{+,-\}$; and when $[f]$ is not analytic, we have $\delta_{\sigma}^{H}(f) = 0 \leq \delta_{\sigma}^{\FTs}(f)$ for each $\sigma\in\{+,-\}$.  So it remains to prove the implication \ref{thm:CK:FTdecay:4'}$\Rightarrow$\ref{thm:CK:FTdecay:1} and that $\delta_{\sigma}^{\FTs}(f) \geq \delta_{\sigma}^{H}(f)$ for each $\sigma\in\{+,-\}$ when $[f]$ is analytic.

Assume \ref{thm:CK:FTdecay:4'}.  By Corollary \ref{cor:aeClass} we may replace $f$ with its analytic representative in $[f]$ and thereby assume that $f$ is analytic.  Theorems \ref{thm:HolomExtInfty} and \ref{thm:CKiS:HolomExtStrip} supply constants $a,b > 0$ such that $f$ extends to a holomorphic function
\[
f:\{z\in\CC : \text{$|\IM(z)| < a$ or $|\RE(z)| > b$}\} \to \CC.
\]
Lemma \ref{lem:LE1}\ref{lem:LE1:L1} and Theorem \ref{thm:HolomExtInfty}\ref{thm:HolomExtInfty:CK} supply constants $r > 1$ and $C>0$ such that $|f(z)| \leq C|z|^{-r}$ for all $z\in\CC$ with $|\RE(z)| > b$.  Thus, for all $t,y\in\RR$ with $|t| > b$,
\[
|f(t+iy)| \leq C|t+\i y|^{-r} \leq C|t|^{-r}.
\]
Therefore for each $\alpha\in(0,\delta_{\sigma}^{H}(f))$, the functions $t\mapsto f(t)$ and $t\mapsto f(\alpha+it)$ are in $\L^1(\RR)$ and \eqref{eq:FTexpDecay:supLimit} holds, so $\FT[f](\xi)\in O(\e^{\sigma 2\pi\alpha\xi})$ as $\xi\to-\sigma\infty$ by \eqref{eq:FTexpDecay:Rbnd} in Lemma \ref{lem:HolExtImplyFTDecay} (with $\Omega = \emptyset$).  This establishes the inequality $\delta_{\sigma}^{\FTs}(f) \geq \delta_{\sigma}^{H}(f) \geq a > 0$, for each $\sigma\in\{+,-\}$, which thereby also establishes statement \ref{thm:CK:FTdecay:1}.
\end{proof}

For each $f\in\C^{\KK}(\RR)\cap\L^1(\RR)$, observe from Theorem \ref{thm:CKiS:HolomExtStrip} that we have $\delta_{\sigma}^{H}(f) > 0$ for each $\sigma\in\{+,-\}$ when $[f]\in C^{\infty}(\RR)$, and we have $\delta_{\sigma}^{H}(f) = 0$ for each $\sigma\in\{+,-\}$ when $[f]\not\in C^{\infty}(\RR)$.  The next proposition therefore elaborates upon the behavior of $[f]$ and its Fourier transform in the extreme situations when the equivalent quantities $\delta_{\sigma}^{H}(f)$ and $\delta_{\sigma}^{\FTs}(f)$ are either $0$ or $+\infty$.

\begin{prop}\label{prop:CK:deltaExtremeCases}
Let $f\in\C^{\KK}(\RR)\cap\L^1(\RR)$.
\begin{enumerate}
\item\label{prop:CK:deltaExtremeCases:0}
The equivalence class $[f]$ is not smooth if and only if for each $\sigma\in\{+,-\}$,
\begin{equation}\label{eq:CK:deltaExtremeCases:0}
\FT[f](\xi) = |\xi|^r(\log|\xi|)^s E(\xi) + o(|\xi|^r(\log|\xi|)^s )
\end{equation}
as $\xi\to\sigma\infty$ for some $r\in\KK_-$, $s\in\NN$, and function $E:\RR\to\CC$ of the form
\[
E(t) = \sum_{\alpha\in A}c_\alpha t^{\i\alpha}
\]
for some nonempty finite $A\subseteq\RR$ and coefficients $\{c_\alpha\}_{\alpha\in A} \subseteq \CC\setminus\{0\}$.

\item\label{prop:CK:deltaExtremeCases:infty}
For each $\sigma\in\{+,-\}$, the equivalence class $[f]$ extends analytically to $H_\sigma$ if and only if  $\FT[f](\xi) = 0$ on $\{\xi\in\RR : \sigma\xi\leq 0\}$.
\end{enumerate}
\end{prop}

\begin{proof}
To prove statement \ref{prop:CK:deltaExtremeCases:0}, note that $[f]\not\in C^\infty(\RR)$ if and only if $\FT[f]$ is slowly decaying at $\pm\infty$ by Theorem \ref{thm:CK:FTexpDecay}, and since $\FT[f]\in\C^{\KK,\FTs}(\RR)$ by Remark \ref{rem:intro:CategFctsProps}\ref{rem:intro:CategFctsProps:CKF}, Proposition \ref{prop:AsymExp} shows that this is in turn equivalent to $\FT[f]$ having the asymptotic behavior at $\pm\infty$ specified in equation \eqref{eq:CK:deltaExtremeCases:0}.

To prove statement \ref{prop:CK:deltaExtremeCases:infty}, let $\sigma\in\{+,-\}$, and observe that if $\FT[f](\xi) = 0$ on $\{\xi\in\RR : \sigma\xi\leq 0\}$, then $\delta_{\sigma}^{H}(f) = \delta_{\sigma}^{\FTs}(f) = +\infty$, so $[f]$ extends analytically to $H_\sigma$.  Conversely, assume that $[f]$ extends analytically to $H_\sigma$.  We may replace $f$ with the analytic representative in $[f]$ to assume that $f$ itself extends analytically to $H_\sigma$.  By Lemma \ref{lem:LE1}\ref{lem:LE1:L1} and Theorem \ref{thm:HolomExtInfty}, we may fix constants $b,C>0$ and $r>1$ such that $|f(z)| \leq C|z|^{-r}$ for all $z\in H_\sigma$ with $|z| \geq b$.  Consider any $\xi\in\RR_{-\sigma}$ and any $R\geq b$.  Then for all $t\in\RR$,
\[
|f(t+\i\sigma R)| \leq C|t+\i\sigma R|^{-r} \leq C|t+\i\sigma b|^{-r}.
\]
By arguing as in the proof of Lemma \ref{lem:HolExtImplyFTDecay},
\[
\FT[f](\xi)
= \int_{\RR}f(t+\i\sigma R)\e^{-\i2\pi\xi(t+\i\sigma R)}\d t
= \e^{2\pi\sigma R \xi}\int_{\RR} f(t+\i\sigma R)\e^{-\i2\pi\xi t} \d t,
\]
so
\[
|\FT[f](\xi)|
\leq \e^{2\pi\sigma R \xi} \int_{\RR} |f(t+\i R)|\d t
\leq C\e^{2\pi\sigma R \xi} \int_{\RR} |t+\i\sigma b|^{-r} \d t.
\]
Since this bound holds for all $R\geq b$, we may let $R\to+\infty$ to see that $\FT[f](\xi) = 0$.  This shows that $\FT[f]$ vanishes on $\RR_{-\sigma}$ and therefore vanishes on $\{\xi\in\RR : \sigma\xi \leq 0\}$ by the continuity of $\FT[f]$.  This proves statement \ref{prop:CK:deltaExtremeCases:infty}.
\end{proof}

\begin{rem}\label{rem:CK:EntireL1=0}
Let $f\in\C^{\KK}(\RR)\cap\L^1(\RR)$.  Recall from Corollary \ref{cor:CK:EntirePoly} that if $f$ extends to an entire function, then $f=0$.  The previous proposition gives an alternate proof of this fact, since $f$ extends to an entire function if and only if $f$ is continuous and $\FT[f] = 0$ on $\RR$ by Proposition \ref{prop:CK:deltaExtremeCases}\ref{prop:CK:deltaExtremeCases:infty}, if and only if $f=0$.
\end{rem}

\subsection{Analytic Extensions on Sectors and the Mellin Transform}\label{ss:DecayFT:Sectors}

We now turn to the task of relating the exponential rate of decay of the Mellin transform of an integrable function $f\in\C^{\KK}(\RR_{\geq 0})$ with the sectors to which $f$ may be extended analytically.

\begin{defs}\label{def:DecayFT:Sector}
Consider a function $f:\RR_{\geq 0}\to\CC$, sign $\sigma\in\{+,-\}$, and constant $c\in\RR$.  Define $E^c[f]:\RR\to\CC$ by
\[
E^c[f](t) = \e^{ct}f(\e^t) \quad\text{for all $t\in\RR$,}
\]
define the sector
\[
S_{\sigma}(\theta) = \{r\e^{\i\alpha}: \text{$r\geq 0$ and $0 \leq \sigma\alpha < \theta$}\}
\]
for each $\theta\geq 0$, and define
\[
\theta_{\sigma}^{H}(f) = \sup\{\theta\in\RR_{\geq 0} : \text{$[f]$ extends analytically to $S_{\sigma}(\theta)$}\}.
\]
\end{defs}

The condition defining $\theta_{\sigma}^{H}(f)$ requires some clarification.  Because $S_\sigma(0) = \emptyset$, the statement ``$[f]$ extends analytically to $S_{\sigma}(0)$'' is a slight misnomer that we take to mean the true statement ``$[f\restriction\emptyset]$ extends analytically to $\emptyset$''.  When $\theta > 0$, we indeed have that $\RR_{\geq 0} \subseteq S_\sigma(\theta)$, and the statement ``$[f]$ extends analytically to $S_{\sigma}(\theta)$'' means that $[f]$ extends analytically to an open subset of $\CC$ containing $S_{\sigma}(\theta)$.  Thus, $[f]$ extends analytically to $S_{\sigma}(2\pi)$ if and only if $[f]$ extends to an entire function.  It follows that
\begin{equation}\label{eq:DecayFT:thetaValue}
\theta_{\sigma}^{H}(f) \in [0,2\pi)\cup\{+\infty\}.
\end{equation}
Since the exponential function $\exp:\CC\to\CC$ is locally biholomorphic, and since $0\in\RR_{\geq 0}\subseteq S_\sigma(\theta)$ for all $\theta > 0$, a straightforward application of the monodromy theorem shows that
\begin{equation}\label{eq:DecayFT:ThetaDeltaExp}
\theta_{\sigma}^{H}(f) = \delta_{\sigma}^{H}(E^c[f]) \quad\text{for any $c\in\RR$.}
\end{equation}

\begin{defn}\label{def:DecayFT:Mellin}
Suppose that $f\in\L^1(\RR_{\leq 0})$ and that $\alpha,\beta\in\RR$ are such that $\alpha < 1 < \beta$, $f(t) = O(t^{-r})$ as $t\to 0^+$ for all $r > \alpha$, and $f(t) = O(t^{-r})$ as $t\to+\infty$ for all $r < \beta$.  Denote the vertical strip $\RR+\i(\alpha,\beta)$ by $S_{\alpha,\beta}$.  The \DEF{Mellin transform} of $f$ is the function $\M[f]:S_{\alpha,\beta}\to\CC$ defined by
\[
\M[f](s) = \int_{0}^{+\infty} t^{s-1} f(t) \d t \quad\text{for all $s\in S_{\alpha,\beta}$.}
\]
The assumptions on $f$, $\alpha$, and $\beta$ ensure that $S_{\alpha,\beta}$ is nonempty and that $\M[f]$ is holomorphic on $S_{\alpha, \beta}$.  Choose any $c\in(\alpha,\beta)$, and define
\[
\qquad \delta_{\sigma}^{\M}(f) = \sup\{a\in\RR_{\geq 0}: \text{$\M(c+\i2\pi\xi) = O(\e^{-\sigma 2\pi a\xi})$ as $\xi\to\sigma\infty$}\}
\]
for each $\sigma\in\{+,-\}$.  We suppress the dependence of this quantity on the choice of $c$ in our notation because, as we shall see shortly, $\delta_{\sigma}^{\M}(f)$ does not actually depend on the choice of $c$ when $f\in\C^{\KK}(\RR)$.
\end{defn}

The next lemma, which is classical, relates the Mellin transform of $f$ and the Fourier transform of $E^c[f]$.

\begin{lem}\label{lem:MellinFT}
Let $f$, $\alpha$, and $\beta$ be as hypothesized in Definition \ref{def:DecayFT:Mellin}.  For all $c\in(\alpha,\beta)$ and $\xi\in\RR$,
\[
\M[f](c+\i2\pi\xi) = \FT[E^c[f]](-\xi).
\]
\end{lem}

\begin{proof}
For each $c\in(\alpha,\beta)$ and $\xi\in\RR$, the change of variables $t = \e^u$ gives
\begin{align*}
\M[f](c+\i2\pi\xi)
    &= \int_{0}^{+\infty} t^{c+\i2\pi\xi-1}f(t) \d t \\
    &= \int_{\RR} e^{u(c+\i2\pi\xi)}f(\e^u) \d u \\
    &= \FT[E^c[f]](-\xi).
\end{align*}
\end{proof}

This leads us to our second main theorem of the section.

\begin{thm}\label{thm:CK:MellinexpDecay}
Consider a smooth $\C^{\KK}$-function $f\in\L^1(\RR_{\geq 0})$.  Then $\theta_{\sigma}^{H}(f) = \delta_{\sigma}^{\M}(f)$ for each $\sigma\in\{+,-\}$.
\end{thm}

To be clear, respectively define $\alpha$ and $\beta$ to be the values of $-\LE_1(f)$ at $0^+$ and $+\infty$ in the affine coordinates $y_{\AF} = y$.  Then $\alpha,\beta\in\RE(\KK)$ and $\alpha < 1 < \beta$ by Lemma \ref{lem:LE1}\ref{lem:LE1:L1}, the function $\M[f]$ is defined on the vertical strip $S_{\alpha,\beta}$, and $\delta_{\sigma}^{\M}(f)$ is defined for each $\sigma\in\{+,-\}$ from a choice of $c\in(\alpha,\beta)$.  However, the quantity $\theta_{\sigma}^{H}(f)$ does not depend on the choice of $c$, so neither does $\delta_{\sigma}^{\M}(f)$ by Theorem \ref{thm:CK:MellinexpDecay}.

\begin{proof}
Let $\alpha$, $\beta$, and $c\in(\alpha,\beta)$ be as specified directly above, and let $\sigma\in\{+,-\}$.  Note that $\theta_{\sigma}^{H}(f) = \delta_{\sigma}^{H}(E^c[f])$ by \eqref{eq:DecayFT:ThetaDeltaExp}, that $\delta_{\sigma}^{\M}(f) = \delta_{\sigma}^{\FTs}(E^c[f])$ by Lemma \ref{thm:CK:MellinexpDecay}, and that $\delta_{\sigma}^{H}(E^c[f]) \geq \delta_{\sigma}^{\FTs}(E^c[f])$ by Lemma \ref{lem:FTDecayImplyHolomExt}.  So it suffices to show that $\delta_{\sigma}^{H}(E^c[f]) \leq \delta_{\sigma}^{\FTs}(E^c[f])$.  To accomplish this, consider any $a\in\RR$ with $0 < a < \delta_{\sigma}^{H}(E^c[f])$.  It suffices to prove that $\delta_{\sigma}^{\FTs}(E^c[f]) \geq a$.

The function $E^c[f]$ extends analytically to the strip $\{z\in\CC : 0\leq \sigma\IM(z) \leq a\}$.  Let $r\in(c,\beta)$.  Proposition \ref{thm:HolomExtInfty:exp} supplies constants $C,R > 0$ such that
\begin{equation}\label{eq:Ec[f]:infty}
|E^c[f](z)| \leq C\e^{(c-r)\RE(z)}
\end{equation}
on $\{z\in \CC : \text{$0\leq \sigma\IM(z) \leq a$ and $\RE(z) > R$}\}$.  The function $f$ is analytic by the version of Theorem \ref{thm:CKiS:HolomExtStrip} described just prior to Corollary \ref{cor:SectorExp}, so $-\alpha\in\NN$ and, by making $C$ larger, there exists an $\epsilon > 0$ such that $|f(z)| \leq C|z|^{-\alpha}$ on $\{z\in\CC : |z| < \epsilon\}$.  By making $R$ larger, we may assume that $R > -\log\epsilon$, so
\begin{equation}\label{eq:Ec[f]:0}
|E^c[f](z)| \leq C\e^{(c-\alpha)\RE(z)}
\end{equation}
on $\{z\in \CC : \text{$0\leq \sigma\IM(z) \leq a$ and $\RE(z) < -R$}\}$. Since $c-\alpha > 0$ and $c-r < 0$, the inequalities \eqref{eq:Ec[f]:infty} and \eqref{eq:Ec[f]:0} show that the functions $t\mapsto E^c[f](t)$ and $t\mapsto E^c[f](t+\i a)$ are in $\L^1(\RR)$ and that
\[
\lim_{t\to\pm\infty}\sup_{\sigma y\in[0,a]}|E^c[f](t+iy)| = 0.
\]
Therefore by Lemma \ref{lem:HolExtImplyFTDecay}, $\FT[E^c[f]](\xi) = O(\e^{\sigma 2\pi a \xi})$ as $\xi\to-\sigma\infty$, which gives $\delta_{\sigma}^{\FTs}(E^c[f]) \geq a$.
\end{proof}

\subsection{Examples and Remarks}\label{ss:DecayFT:ExamRmk}

We now expand upon the content of Subsections \ref{ss:DecayFT:Strips} and \ref{ss:DecayFT:Sectors} in a series of examples and remarks.  The first example merely serves to concretely demonstrate behavior described in Theorems \ref{thm:CK:FTexpDecay} and \ref{thm:CK:MellinexpDecay}.

\begin{exam}\label{exam:DecayFT:MainThmsDemo}
Consider the analytic function $f\in\S(\RR)$ defined by
\[
f(t) = \frac{1}{1+t^2} \quad\text{for all $t\in\RR$.}
\]
The function $f$ extends meromorphically to $\CC$ with simple poles at $\pm\i$, so $\delta^{H}_{\sigma}(f)=1$ and $\theta_{\sigma}^{H}(f) = \frac{\pi}{2}$ for each $\sigma\in\{+,-\}$.  Also, $\FT[f](\xi)=\pi \e^{-2\pi \vert \xi \vert}$, so we can directly observe that $\delta^{\FTs}_{\sigma}(f) = 1$ for each $\sigma\in\{+,-\}$, in agreement with Theorem \ref{thm:CK:FTexpDecay}.  We have $\LE_1(f) = 0$ at $0^+$ and $\LE_1(f) = -2$ at $+\infty$ in the affine coordinates $t_{\AF} = t$, and a computation shows that $\M[f]:S_{0,2}\to\CC$ is given by
\[
\M[f](s) = \frac{\pi}{2\sin(\pi s/2)} \quad\text{for all $s\in S_{0,2}$.}
\]
By choosing $c\in(0,2)$ to be $c=1$, we have
\[
\mathcal{M}[f](1+2\i\pi \xi)=
\frac{\pi}{\e^{\pi^2 \xi}+\e^{-\pi^2\xi}} \sim \pi \e^{-\pi^2|\xi|}
\,\,\text{as $\xi\to\pm\infty$,}
\]
so we can directly observe that $\delta_{\sigma}^{\M}(f) = \frac{\pi}{2}$ for each $\sigma\in\{+,-\}$, in agreement with Theorem \ref{thm:CK:MellinexpDecay}.
\end{exam}

We know from Theorem \ref{thm:CK:FTexpDecay} that the Fourier transform of a smooth and integrable function $f\in\C^{\KK}(\RR)$ decays exponentially at $\pm\infty$.  The next remark explains how the same can be said about derivatives of $\FT[f]$ when $f$ has a sufficiently fast decay at both $+\infty$ and $-\infty$.

\begin{rem}\label{rem:DecayFT:FTderiv}
Consider a nonzero, smooth, $\C^{\KK}$-function $f\in \L^1(\RR)$.  Let $\overline{r}$ be the minimum value of $-\LE_1(f)$ at both $+\infty$ and $-\infty$ in the affine coordinates $y_{\AF} = \pm y$, and note that $\overline{r} > 1$ by Lemma \ref{lem:LE1}\ref{lem:LE1:L1}.  Then for all $k\in\NN$ with $k < \overline{r}-1$, the function $\FT[f]^{(k)}$ exists, is continuous, and decays exponentially at $\pm\infty$.
\begin{proof}
Let $k\in\NN$ with $k < \overline{r} - 1$.  Then $k-\overline{r} < -1$, so the smooth $\C^{\KK}$-function $f_k(t) = (-i2\pi t)^kf(t)$ is in $\L^1(\RR)$.  Therefore, $\FT[f_k]$ is continuous and decays exponentially at $\pm\infty$ by Theorem \ref{thm:CK:FTexpDecay}, and $\FT[f]^{(k)} = \FT[f_k]$ by Proposition \ref{prop:FTmotivate}\ref{prop:FTmotivate:FTdiff}.
\end{proof}
\noindent
For the function $f$ of Example \ref{exam:DecayFT:MainThmsDemo}, we have  $\overline{r} = 2$ and $\FT[f] \in C^0(\RR) \setminus C^1(\RR)$.  This shows that the bound $k < \overline{r}-1$ given in Remark \ref{rem:DecayFT:FTderiv} is the best possible, in general.
\end{rem}

The next remarks explain how when an integrable analytic function $f:\RR\to\CC$ extends meromorphically to a horizontal strip about the real axis,  we can slightly improve upon the description of the asymptotic behavior of $\FT[f]$ given in Theorem \ref{thm:CK:FTexpDecay} for $\C^{\KK}$-functions.

\begin{rems}\label{rems:DecayFT:Sing}
Consider again the situation hypothesized in Lemma \ref{lem:HolExtImplyFTDecay}.  Namely, we consider a constant $\alpha>0$, a sign $\sigma\in\{+,-\}$, a finite subset $\Omega$ of the open horizontal strip $\RR+\sigma\i(0,\alpha)$, and a holomorphic function $f:(\RR+\sigma\i[0,\alpha])\setminus\Omega \to\CC$ such that the functions $t\mapsto f(t)$ and $t\mapsto f(t+\sigma\i\alpha)$ are in $\L^1(\RR)$ and
\[
\lim_{t\to\pm\infty}\sup_{\sigma y\in[0,\alpha]}|f(t+iy)| = 0.
\]
Additionally assume that $\Omega$ is nonempty and that for each $\omega\in\Omega$, the singularity of $f$ at $\omega$ is not removable.  Then $\delta_{\sigma}^{H}(f) = \min\{|\IM(\omega)| : \omega\in\Omega\}$.
\begin{enumerate}
\item\label{rems:DecayFT:Sing:Residues}
For each $\omega\in\Omega$,
\begin{equation}\label{eq:FTDecayImplyHolomExt:Residues}
\Res{z=\omega} f(z)\e^{-\i 2\pi\xi z} = \e^{-\i2\pi\omega\xi} \sum_{j=1}^{+\infty} P_{\omega}^{j}[f]\frac{(-\i 2\pi\xi)^{j-1}}{(j-1)!},
\end{equation}
where
\[
P_{\omega}[f](z) = \sum_{j=1}^{+\infty}\frac{P_{\omega}^{j}[f]}{(z-\omega)^j}
\]
is the principal part of $f$ at $\omega$.

\begin{proof}
We have
\begin{align*}
\qquad
P_\omega[f](z) \e^{-\i2\pi\xi z}
    &= \e^{-\i2\pi\omega\xi}\left(\sum_{j=1}^{+\infty}\frac{P_{\omega}^{j}[f]}{(z-\omega)^j}\right)\left(\sum_{k=0}^{+\infty}\frac{(-\i2\pi\xi(z-\omega))^k}{k!}\right) \\
    &= \e^{-\i2\pi\omega\xi} \sum_{(j,k)\in\ZZ_+\times\NN} P_{\omega}^{j}[f]\frac{(-\i2\pi\xi)^k}{k!}(z-\omega)^{k-j},
\end{align*}
whose coefficient of $(z-\omega)^{-1}$ is found by summing over all the coefficients with $j\in\ZZ_+$ and $k = j - 1$.
\end{proof}

\item\label{rems:DecayFT:Sing:Poles}
When $\omega\in\Omega$ is a pole of $f$ of order $\mu(\omega)$, equation \eqref{eq:FTDecayImplyHolomExt:Residues} becomes
\begin{align*}
\Res{z=\omega} f(z)\e^{- \i 2\pi\xi z}
&= \e^{-\i 2\pi\omega\xi} \sum_{j=1}^{\mu(\omega)} P_{\omega}^{j}[f]\frac{(-\i 2\pi\xi)^{j-1}}{(j-1)!} \\
&\sim C\xi^{\mu(\omega)-1}\e^{2\pi\IM(\omega)\xi} \text{ as $\xi\to-\sigma\infty$,}
\end{align*}
for some constant $C\neq 0$.  Therefore, if each member $\omega$ of the set
\[
\Omega' := \{\omega'\in\Omega : |\IM(\omega')| = \delta_{\sigma}^{H}(f)\}
\]
is a pole of order $\mu(\omega)$, and if we let $\mu = \max\{\mu(\omega') : \omega'\in\Omega'\}$, then Lemma \ref{lem:HolExtImplyFTDecay} shows that
\begin{equation}\label{eq:DecayFT:Sing:Poles}
\FT[f](\xi) \sim C\xi^{\mu-1}\e^{\sigma 2\pi\delta_{\sigma}^{H}(f)\xi} \quad\text{as $\xi\to-\sigma\infty$,}
\end{equation}
for some constant $C\neq 0$.  When $f\in\C^{\KK}(\RR)$, Theorem \ref{thm:CK:FTexpDecay} shows that
\begin{equation}\label{eq:DecayFT:Sing:Thm1}
\FT[f](\xi) = O(\e^{\sigma 2\pi a\xi}) \,\,\text{as $\xi\to-\sigma\infty$, for each $a < \delta_{\sigma}^{H}(f)$.}
\end{equation}
Observe how \eqref{eq:DecayFT:Sing:Poles} is a stronger assertion than \eqref{eq:DecayFT:Sing:Thm1}.
\end{enumerate}
\end{rems}

If $f\in\C^{\KK}(\RR)$ and the set $\Omega'$ contains an essential singularity, then to the best of our knowledge, it can be difficult to improve upon the asymptotic estimate of $\FT[f]$ given by Theorem \ref{thm:CK:FTexpDecay}.  This is because it seems to be difficult to precisely determine the asymptotic behavior of the residue \eqref{eq:FTDecayImplyHolomExt:Residues} as $\xi\to-\sigma\infty$ when $\omega$ is an essential singularity, as demonstrated in the next example.   To prepare for the example, we first prove a lemma which is also of independent interest.

\begin{lem}\label{lem:DecayFT:CKmeromFT}
Let $f\in\C^{\KK}(\RR)\cap\L^1(\RR)$ and $\sigma\in\{+,-\}$, assume that $f$ extends to an analytic function on $H_\sigma\setminus\Omega$ for some finite subset $\Omega$ of $H_\sigma$.  Then using the notation from the Remarks \ref{rems:DecayFT:Sing}, we have
\begin{equation}\label{eq:DecayFT:CKmeromFT}
\FT[f](\xi) = \sum_{\omega\in\Omega} \left(\e^{-\i2\pi\omega\xi} \sum_{j=1}^{\mu(\omega)} P_{\omega}^{j}[f]\frac{(-\i 2\pi\xi)^{j-1}}{(j-1)!}\right) \quad\text{for all $\xi\in\RR_{-\sigma}$,}
\end{equation}
where $\mu(\omega)$ is the order of the singularity $\omega$, with the understanding that $\mu(\omega) = +\infty$ when $\omega$ is an essential singularity.
\end{lem}

\begin{proof}
Let $\xi\in\RR_{-\sigma}$.  Theorem \ref{thm:HolomExtInfty}\ref{thm:HolomExtInfty:CK} supplies constants $C,R > 0$ and $r\in(1,2)$ such that $f$ is analytic and satisfies
\[
|f(t+\i y)| \leq C|t+\i y|^{-r} \leq g(t)
\]
on the set $\{(t+\i y\in H_\sigma : |t+\i y| > R\}$, where $g\in \L^1(\RR)$ is defined by
\[
g(t) = C(\max\{R,|t|\})^{-r} \quad\text{for all $t\in\RR$.}
\]
Lemma \ref{lem:HolExtImplyFTDecay} and the Remarks \ref{rems:DecayFT:Sing} therefore show that for each $\alpha > R$, we may write $\FT[f](\xi)$ as the sum of the function on the right side of equation \eqref{eq:DecayFT:CKmeromFT} and a remainder term $\R_{\alpha}[f](\xi)$ satisfying
\[
|\R_{\alpha}[f](\xi)| \leq \e^{\sigma 2\pi\alpha\xi}\int_{\RR} g(t) \d t.
\]
We have $\lim_{\alpha\to+\infty} \e^{\sigma 2\pi\alpha\xi} = 0$ because $\sigma\xi < 0$, so the lemma follows by letting $\alpha\to+\infty$.
\end{proof}

\begin{exam}\label{exam:DecayFT:EssentialSing}
Consider the analytic function $f:\RR\to\CC$ defined by
\[
f(t)=\frac{\e^{\frac{\i}{1+t^2}}}{1+t^2}, \quad\text{for all $t\in\RR$,}
\]
which extends to a holomorphic function on $\CC\setminus\{\pm\i\}$ with essential singularities at $\pm\i$.  Observe that $f\in\L^1(\RR) \cap \S_{\CC}(\RR)$ and that $\S_{\CC}(\RR) \subseteq \C^{\KK}(\RR)$ when $\KK\not\subseteq\RR$.  Therefore for each $\sigma\in\{+,-\}$, we have $\delta_{\sigma}^{\FTs}(f) = 1$ by Theorem \ref{thm:CK:FTexpDecay}.  This observation does not use the residues of $f$ at its two singular points, so let us now use those residues to compute $\FT[f]$ more explicitly.

Consider any $\sigma\in\{+,-\}$, $\xi\in\RR_{-\sigma}$, and $\alpha > 1$.  An elementary computation shows that the principal part of $f$ at $\sigma\i$ is
\[
\sum_{n=0}^{\infty}\left(\sum_{k=n}^{\infty} \frac{(-1)^{k-n} \i^k (2k-n)!}{(k-n)!(k!)^2(\sigma 2\i)^{2k-n+1}}\right) \frac{1}{(z-\sigma\i)^{n+1}},
\]
so by Remark \ref{rems:DecayFT:Sing}\ref{rems:DecayFT:Sing:Residues},
\[
\Res{z = \sigma\i} f(z)\e^{-\i 2\pi\xi z}
=
\e^{-\i 2\pi(\sigma\i)\xi}
\sum_{n=0}^{\infty} \frac{(-\i 2\pi\xi)^n}{n!}\left(\sum_{k=n}^{\infty} \frac{(-1)^{k-n} \i^k (2k-n)!}{(k-n)!(k!)^2(\sigma 2\i)^{2k-n+1}}\right).
\]
Thus, after some simplification,
\[
\sigma\i 2\pi\Res{z = \sigma\i} f(z)\e^{-\i 2\pi\xi z} = \e^{-2\pi|\xi|}\varphi_{\sigma}(\xi)
\]
for the function
\[
\varphi_{\sigma}(\xi) = \pi \sum_{n=0}^{\infty} \frac{(-\sigma\pi\xi)^n}{n!}\left(
\sum_{k=n}^{\infty} \frac{\i^k (2k-n)!}{(k-n)!(k!)^2 4^{k-n}}\right).
\]
Therefore by Lemma \ref{lem:DecayFT:CKmeromFT},
\[
\FT[f](\xi) = \e^{-2\pi|\xi|}\varphi_{\sigma}(\xi) \,\,\text{on $\RR_{\sigma}$.}
\]
Since $\delta_{\sigma}^{\FTs}(f) = 1$, we see that for each $r > 0$, the function $\varphi_\sigma(\xi)$ grows slower than $\e^{r|\xi|}$ as $\xi \to -\sigma\infty$, but there does not exist an $s > 0$ such that $\varphi_\sigma(\xi)$ grows slower than $\e^{-s|\xi|}$ as $\xi \to -\sigma\infty$.  Nevertheless, the precise asymptotic behavior of $\varphi_\sigma(\xi)$ as $\xi\to-\sigma\infty$ seems to be not easily accessible from its series representation.
\end{exam}

The next example demonstrates how a function $f\in\C^{\KK}(\RR)$ can extend analytically to horizontal strips of different widths in the upper and lower half planes when $\KK\not\subseteq\RR$, which then causes different asymptotic behavior of $\FT[f]$ at $+\infty$ and $-\infty$.

\begin{exam}\label{exam:DecayFT:HalfPlanes}
Consider the analytic function $f\in\L^1(\RR)$ defined by
\[
f(t) = \frac{1}{(t-\i)^2} \quad\text{for all $t\in\RR$,}
\]
and note that $f\in (\S+\i\S)(\RR)$ and that $(\S+\i\S)(\RR) \subseteq \S^{\KK}(\RR)$ when $\KK\not\subseteq\RR$.  Then $\delta_{+}^{H}(f) = 1$ and $\delta_{-}^{H}(f) = +\infty$, so $\delta_{+}^{\FTs}(f) = 1$ and $\delta_{-}^{\FTs}(f) = +\infty$ by Theorem \ref{thm:CK:FTexpDecay}.  More specifically, we have
\[
\FT[f](\xi) = \begin{cases}
0,                  & \text{if $\xi \geq 0$,} \\
4\pi^2\xi\e^{2\pi\xi}, & \text{if $\xi < 0$,}
\end{cases}
\]
with the formula for $\xi\geq 0$ following from Proposition \ref{prop:CK:deltaExtremeCases}\ref{prop:CK:deltaExtremeCases:infty} and the formula for $\xi < 0$ following from Lemma \ref{lem:DecayFT:CKmeromFT}.
\end{exam}

Proposition \ref{prop:CK:deltaExtremeCases}\ref{prop:CK:deltaExtremeCases:infty} shows that the only $f\in\C^{\KK}(\RR)\cap\L^1(\RR)$ with $\delta_{+}^{\FTs}(f) = \delta_{-}^{\FTs}(f) = +\infty$ is the function $f=0$.  The next example shows that this assertion fails in the larger category $\C^{\KK,\FTs}$.  We do not know whether this assertions holds for the intermediary category $\C^{\KK,\i\S}$.

\begin{exam}\label{exam:DecayFT:Schwartz}
The Schwartz function $f(t) = \e^{-t^2}$ is in $\C^{\QQ,\FTs}(\RR) \cap \L^1(\RR)$, and its Fourier transform $\FT[f](\xi) = \sqrt{\pi}\e^{-\pi^2\xi^2}$ decays faster than any exponential function as $\xi\to\pm\infty$.  So $\delta_{\sigma}^{\FTs}(f) = +\infty$ for each $\sigma\in\{+,-\}$, even though $f\neq 0$.
\end{exam}

Our final example gives a smooth function $f\in\C^{\QQ,\i\S}(\RR)\cap\L^1(\RR)$ for which $\FT[f]$ decays slowly at $\pm\infty$, thereby showing that Theorem \ref{thm:CK:FTexpDecay} fails in the larger category $\C^{\KK,\i\S}$.

\begin{exam}\label{exam:DecayFT:Thm1Fail}
Define the analytic function $f\in\C^{\QQ,\i\S}(\RR)\cap\L^1(\RR)$ by
\[
f(t) = \frac{\e^{\i t^2}}{1+t^2} \quad\text{for all $t\in\RR$,}
\]
and let $\sigma\in\{+,-\}$.  Then
\begin{equation}\label{eq:CKiS-SlowDecayFT}
\FT[f](\xi) \sim -\sigma\frac{\e^{\i\pi/4}}{\pi^{3/2}}
\frac{ \e^{-\i(\pi\xi)^2}}{\xi^2}
\quad\text{as $\xi\to \sigma\infty$.}
\end{equation}
This shows that $\FT[f]$ has a slow decay even though $f$ is analytic.  Thus, $\delta_{\sigma}^{\FTs}(f) = 0$, but clearly $\delta_{\sigma}^{H}(f) = 1$.
\end{exam}

\begin{proof}[Proof of \eqref{eq:CKiS-SlowDecayFT}]
Consider any $\xi\in\RR_\sigma$, and let
\[
I(\xi)=\FT[f]\left(\frac{\xi}{2\pi}\right)=
\int_{-\infty}^{+\infty} \frac{\e^{\i t^2}}{1+t^2}\e^{-\i t\xi } \d t.
\]
The change of variables $t=\xi u$ gives
\[
I(\xi) = \sigma\xi \int_{-\infty}^{+\infty} \frac{\e^{ \xi^2S(u)}}{1+\xi^2u^2} \d u,
\]
where $S(u)=\i(u^2-u)$. At the nondegenerate critical point $\frac{1}{2}$ of $S$, we can apply the saddle point method, and  we use the residue theorem to compute $I(\xi)$ by integrating on the oriented line $t\mapsto \frac{1}{2}-t\e^{\i \frac{\pi}{4}}$.  This gives
\[
I(\xi) =
\sigma\xi 2\i \pi \underset{u=-\frac{\i}{\sigma\xi}}{\mathrm{Res}}\left(\frac{\e^{\xi^2S(u)}}{1+\xi^2u^2} \right) +\sigma\xi\e^{\frac{\i\pi}{4}}\e^{-\frac{\i\xi^2}{4}} \int_{-\infty}^{+\infty} \frac{\e^{ -\xi^2 t^2}}{1+\xi^2(\frac{1}{2}-t\e^{\frac{\i\pi}{4}})^2} \d t,
\]
with $\displaystyle\underset{u=-\frac{\i}{\sigma\xi}}{\mathrm{Res}}\left(\frac{\e^{\xi^2S(u)}}{1+\xi^2u^2} \right)=
-\frac{\e^{-\i}\e^{-|\xi|}}{\sigma\xi 2\i}$.  Next write
\[
h(\xi,t) := \frac{h(t)}{\xi^2} :=  \frac{1}{  \xi^2(\frac{1}{2}-t\e^{\frac{\i\pi}{4}})^2}.
\]
With this notation, we have
\begin{align*}
\frac{ 1 }{1+\xi^2(\frac{1}{2}-t\e^{\frac{\i\pi}{4}})^2}
    &=\frac{h(t)}{\xi^2}\frac{1}{1+\frac{h(t)}{\xi^2}}\\
    &=\sum_{k=0}^N \frac{(-1)^{k}}{\xi^{2k+2}} h^{k+1}(t) + \frac{(-1)^{N+1}}{\xi^{2N+4}}\frac{h^{N+2}(t)}{1+h(\xi,t)}.
\end{align*}
We obtain
\begin{align}\label{eq:CKiS-SlowDecayFT:FiniteExpansion}
I(\xi)
    &= -\pi\e^{-\i}\e^{-|\xi|} +\sigma \e^{\frac{\i\pi}{4}}\e^{-\frac{\i\xi^2}{4}} \sum_{k=0}^N \frac{(-1)^{k}}{\xi^{2k+1}} \int_{-\infty}^{+\infty} h^{k+1}(t) \e^{ -\xi^2 t^2} \d t
    \\
    &\phantom{=}
    +\sigma\frac{(-1)^{N}}{\xi^{2N+3}} \e^{\frac{\i\pi}{4}}\e^{-\frac{\i\xi^2}{4}} \int_{-\infty}^{+\infty} \frac{h^{N+2}(t)}{1+h(\xi,t)} \e^{ -\xi^2 t^2} \d t. \nonumber
\end{align}
For all $\xi,t\in\RR$ with $|\xi|\geq 3$, we have $|\frac{1}{2}-t\e^{\i\pi/4}| \geq \dfrac{1}{2\sqrt{2}}$ by elementary geometry, so
\[
|1+h(\xi,t)| \geq 1-|h(\xi,t)| \geq 1 - \frac{8}{\xi^2} \geq \frac{1}{9}.
\]
Thus, for $|\xi|\geq 3$, the last integral in equation \eqref{eq:CKiS-SlowDecayFT:FiniteExpansion} is bounded by a constant only depending on $N$.  On the other hand,
    the other integrals in \eqref{eq:CKiS-SlowDecayFT:FiniteExpansion} are Laplace integrals that
    can be directly expanded in the monomial scale
    (see \cite[43.21 and 43.27]{SidFedSha}):
\[
\int_{-\infty}^{+\infty}
h^{k+1}(t) \e^{-\xi^2 t^2} \d t
\sim
\sum_{p=0}^{+\infty} \frac{c_{p,k}}{|\xi|^{2p+1}} \quad\text{as $\xi\to\sigma\infty$,}
\]
where $\displaystyle c_{p,k}=\Gamma\left(p+\frac{1}{2}\right)\frac{(h^{k+1})^{(2p)}(0)}{(2p)!}$.  In particular
\[
I(\xi) \sim
 \sigma 4\sqrt{\pi} \e^{\frac{\i\pi}{4}}
      \frac{\e^{-\frac{\i\xi^2}{4}}}{\xi^2} \quad\text{as $\xi\to \sigma\infty$,}
\]
from which \eqref{eq:CKiS-SlowDecayFT} follows.
\end{proof}

\section{Integrability of the Fourier Transform}\label{s:IntegFT}

Consider a function $f\in\L^1(\RR)$.  As discussed in the Introduction, the Fourier inversion theorem shows that a necessary condition for $\FT[f]$ to be in $\L^1(\RR)$ is that $[f]$ is continuous.  The next theorem and Corollary \ref{cor:aeClass} together show that the continuity of $[f]$ is also a sufficient condition for $\FT[f]$ to be in $\L^1(\RR)$ when $f\in\C^{\KK}(\RR)$.  The theorem extends Cluckers and Miller \cite[Theorem 1.3]{CluMil}, which considered the special case of $\KK = \QQ$.

\begin{thm}\label{thm:IntegFT:CK}
If $f\in\C^{\KK}(\RR)\cap\L^1(\RR)$ is continuous, then $\FT[f]\in \L^1(\RR)$.
\end{thm}

\begin{proof}
Assume that $f\in\C^{\KK}(\RR)$ is Lebesgue integrable and continuous.  Then by Proposition \ref{prop:CFK:anal}, there exist finitely many points
\[
-\infty = a_0 < a_1 < \cdots < a_{n-1} < a_n = +\infty
\]
such that $f$ is analytic on $(a_{j-1},a_j)$ for all $j\in\{1,\ldots,n\}$, and $f'\in \L^1(\RR)$ by Lemma \ref{lem:CK:L1}\ref{lem:CK:L1:Cn}.  So for each $\xi\in\RR$, integration by parts gives
\begin{align}\label{eq:IntegFT:CK}
\FT[f](\xi)
    &= \sum_{j=1}^{n}\int_{a_{j-1}}^{a_j} f(t)\e^{-\i 2\pi\xi t} \d t \\
    &= \frac{1}{-\i 2\pi\xi} \sum_{j=1}^{n}\left([f(t)\e^{-\i 2\pi\xi t}]_{a_{j-1}}^{a_j} - \int_{a_{j-1}}^{a_j} f'(t)\e^{-\i2\pi\xi t} \d t\right) \nonumber\\
    &= \frac{1}{\i2\pi\xi}\FT[f'](\xi). \nonumber
\end{align}
Since $\FT[f']\in\C^{\KK,\FTs}(\RR)$ by Remark \ref{rem:intro:CategFctsProps}\ref{rem:intro:CategFctsProps:CKF} and since $\lim_{\xi\to\pm\infty}\FT[f'](\xi) = 0$ by the Riemann-Lebesgue lemma, Lemma \ref{lem:LE1}\ref{lem:LE1:lim0} shows that $\LE_1(\FT[f']) < 0$ in the affine coordinates $y_{\AF} = \pm y$ at $\pm\infty$.  Therefore $\LE_1(\FT[f]) < -1$ at $\pm\infty$ by equation \eqref{eq:IntegFT:CK}, so $\FT[f]\in \L^1(\RR)$.
\end{proof}

\begin{que}
Does Theorem \ref{thm:IntegFT:CK} remain true if we replace $\C^{\KK}$ with either $\C^{\KK,\i\S}$ or $\C^{\KK,\FTs}$?
\end{que}

We do not know the answer to this question for $\C^{\KK,\i\S}$, but Example \ref{exam:NonintegFT:Cexp} below gives a negative answer for $\C^{\KK,\FTs}$.  We introduce the next concept to prepare for the example.

\begin{defn}\label{def:RiemFT}
Consider a function $f:\RR\to\CC$ for which there exist finitely many points
\[
-\infty = a_0 < a_1 < a_2 < \cdots < a_{n-1} < a_n = +\infty
\]
such that $f$ is Riemann integrable\footnote{More generally, one may assume that $f$ is Lebesgue integrable on every compact interval contained in $\RR\setminus\{a_1,\ldots,a_{n-1}\}$.} on every compact interval contained in $\RR\setminus\{a_1,\ldots,a_{n-1}\}$.  By choosing any points $b_1,\ldots,b_n$ with
\[
-\infty = a_0 < b_1 < a_1 < b_2 < a_2 < \cdots < a_{n-1} < b_n < a_n = +\infty,
\]
define the improper integral
\[
\int_{\RR}f(t)\d t := \sum_{j=1}^{n}\left(\left(\lim_{r\to a_{j-1}^{+}}\int_{r}^{b_j}f(t)\d t\right) + \left(\lim_{r\to a_{j}^{-}}\int_{b_j}^{r}f(t)\d t\right)\right),
\]
with the understanding that, by definition, this improper integral converges if and only if each of these individual limits exist.  When the improper integral
\[
\RFT[f](\xi) := \int_{\RR}f(t)\e^{-\i 2\pi\xi t}\d t
\]
converges for each $\xi\in\RR$, we call $\RFT[f]$ the \DEF{Riemann-Fourier transform of $f$}.
\end{defn}

\begin{lem}\label{lem:ImproperInt}
Let $f,\phi\in C^1(\RR)$ (both real-valued), and suppose that $\frac{f(t)}{\phi'(t)} \to 0$ monotonically as $t\to\pm\infty$ and that $\phi$ has finitely many stationary points.  Then the improper integral $\int_{\RR}f(t)\e^{\i\phi(t)}\d t$ converges.
\end{lem}

\begin{proof}
Let $f$ and $\phi$ be as hypothesized, and fix $R_0 > 0$ such that all of the stationary points of $\phi$ are contained in $(-R_0,R_0)$ and such that $\frac{f}{\phi'}$ is monotonic on both $(-\infty,-R_0]$ and $[R_0,+\infty)$.  It suffices to show that the improper integrals $\int_{R_0}^{+\infty}f(t)\e^{\i\phi(t)}\d t$ and $\int_{-\infty}^{-R_0}f(t)\e^{\i\phi(t)}\d t$ both converge.  Both proofs are similar, so we verify this for the first integral.  The functions $f$ and $\phi'$ both have constant positive or negative signs on $[R_0,+\infty)$, so up to possibly negating and taking the complex conjugate of the integrand, we may assume without loss of generality that $f$ and $\phi'$ are both positively-valued on $[R_0,+\infty)$.  The function $\phi$ is invertible on $[R_0,+\infty)$, so by considering any $S,R\in\RR$ with $S > R > R_0$ and using the substitution $u = \phi(t)$, we get
\begin{align*}
\int_{R}^{S} f(t)\e^{\i\phi(t)}\d t
    &= \int_{\phi(R)}^{\phi(S)} \left(\frac{f}{\phi'}\right)\circ\phi^{-1}(u) \e^{\i u}\d u \\
    &= \left(\frac{f}{\phi'}\right)\circ \phi^{-1}(\phi(R)) \int_{\phi(R)}^{u_0}\e^{\i u} \d u
\end{align*}
for some $u_0\in[\phi(R),\phi(S)]$ by the second mean value theorem for integrals, which may be applied because $(\frac{f}{\phi'})\circ\phi^{-1}$ is decreasing on $[\phi(R),\phi(S)]$.  Thus,
\[
\left|\int_{R}^{S} f(t)\e^{\i\phi(t)}\d t\right| \leq \frac{2 f(R)}{\phi'(R)} \to 0 \quad\text{as $R\to+\infty$,}
\]
which is a Cauchy condition that shows that $\int_{R_0}^{+\infty}f(t)\e^{i\phi(t)}\d t$ converges.
\end{proof}

\begin{exams}\label{exam:RFT:RL-Lemma-Fails}
We now give two examples demonstrating how the conclusion of the Riemann-Lebesgue lemma can fail for the Riemann-Fourier transform of a $\C^{\QQ,\i\S}$-function $f$ with $f\not\in\L^1(\RR)$.
\begin{enumerate}
\item\label{exam:RFT:RL-Lemma-Fails:Fresnel}
Define $f\in\C^{\QQ,\i\S}(\RR)$ by $f(t) = \e^{\i t^2}$ for all $t\in\RR$.  Then $f\not\in\L^1(\RR)$, but
\[
\RFT[f](\xi)=\sqrt{\pi}\e^{\i\frac{\pi}{4}}\e^{-\i(\pi\xi)^2}
\]
for all $\xi\in\RR$.  Thus, $\RFT[f]$ does not decay at $\pm\infty$.

\begin{proof}
Let $\xi\in\RR$.  Lemma \ref{lem:ImproperInt} shows that $\RFT[f](\xi)$ converges, we have
\[
\RFT[f](\xi) = \int_{\RR}\e^{\i t^2}\e^{-\i2\pi\xi t}\d t = \e^{-\i(\pi\xi)^2}\int_{\RR}\e^{\i(t-\pi\xi)^2}\d t,
\]
and the Fresnel integral $\int_{\RR}\e^{\i(t-\pi\xi)^2}\d t = \int_{\RR}\e^{\i u^2}\d u$ equals $\sqrt{\pi}\e^{i\pi/4}$.
\end{proof}

\item\label{exam:RFT:RL-Lemma-Fails:MainExample}
Define $f\in\C^{\QQ,\i\S}(\RR)$ by $f(t) = t\e^{\i t^4/4}$ for all $t\in\RR$.  Then $f\not\in\L^1(\RR)$, but $\RFT[f](\xi)$ converges for all $\xi\in\RR$ by Lemma \ref{lem:ImproperInt}, and for each $\sigma\in\{+,-\}$,
\begin{equation}\label{eq:RFT:RL-Lemma-Fails:MainExample}
\RFT[f]\left(\frac{\xi}{2\pi}\right) \sim \sigma\sqrt{\frac{2\pi}{3}} \e^{\i\frac{\pi}{4}}\e^{-\i\frac{3}{4}\xi^{\frac{4}{3}}} \quad\text{as $\xi\to\sigma\infty$.}
\end{equation}
Thus, $\RFT[f]$ does not decay at $\pm\infty$.

\begin{proof}[Proof of \eqref{eq:RFT:RL-Lemma-Fails:MainExample}]
We apply the method of the stationary phase.  For any $\xi\neq 0$, the change of variables $t=\xi^{\frac{1}{3}}u$ gives
\[
\RFT[f]\left(\frac{\xi}{2\pi}\right)
= \sigma\xi^{\frac{2}{3}}\int_\RR u\e^{\i\xi^{\frac{4}{3}}S(u)} \d u,
\]
where $S(u)=\frac{1}{4}u^4 - u$ and $\sigma$ is the sign of $\xi$.  Since $1$ is the only real stationary point of $S$, $S(1)=-\frac{3}{4}$, and $S''(1)=3\neq 0$, we can apply \cite[44.2 Theorem 2]{SidFedSha} to obtain \eqref{eq:RFT:RL-Lemma-Fails:MainExample}.
\end{proof}
\end{enumerate}
\end{exams}

We are now ready for the main example of the section.

\begin{exam}\label{exam:NonintegFT:Cexp}
Consider the function $f(t)=t\e^{\i t^4/4}$ from Example \ref{exam:RFT:RL-Lemma-Fails}\ref{exam:RFT:RL-Lemma-Fails:MainExample}, and define the function $g:\RR\to\CC$ by
\[
g(x) = -\int_{x}^{+\infty}f(t)\d t, \quad\text{for all $x\in\RR$,}
\]
where the integral is improper and convergent by Lemma \ref{lem:ImproperInt}.  Then $g$ is a continuous and Lebesgue integrable function in $\C^{\QQ,\FTs}(\RR)$ with $\FT[g]\not\in\L^1(\RR)$.
\end{exam}

\begin{proof}
The continuity of $g$ follows immediately from the fact that $g(y)-g(x) = \int_{x}^{y}f(t)\d t$ for all $x,y\in\RR$.  Letting $D = \{(x,y)\in\RR^2 : x<y\}$, the function $(x,y,t)\mapsto -\chi_{[x,y]}(t)f(t)$ is in $\C^{\QQ,\i\S}(D\times\RR)$ and is Lebesgue integrable in $t$ for each $(x,y)\in D$.  So $(x,y)\mapsto -\int_{x}^{y} f(t) \d t$ is in $\C^{\QQ,\FTs}(D)$ by stability of $\C^{\QQ,\FTs}$ under parametric integration (see Remark \ref{rem:intro:CategFctsProps}\ref{rem:intro:CategFctsProps:CKF}).  Since the improper integral defining $g$ converges, we see that $g\in\C^{\QQ,\FTs}(\RR)$ by letting $y\to+\infty$ and applying Proposition \ref{prop:Limits}.

We now show that $g\in \L^1(\RR)$.  Since $g$ is continuous on $\RR$, it suffices to study the asymptotic behavior of $g$ at $\pm\infty$.   Let us first assume that $x\geq 1$. Then using integration by parts for improper integrals, we have
\begin{align*}
g(x) &=\int_{x}^{+\infty} \frac{1}{t^2}  t^3\e^{\i \frac{t^4}{4}} \d t
=
\left[ \frac{1}{t^2} \e^{\i \frac{t^4}{4}}\right]_{x}^{+\infty} +
2\int_{x}^{+\infty}    \frac{\e^{\i \frac{t^4}{4}}}{t^3} \d t  \\
&
= {-\frac{\e^{\i\frac{x^4}{4}}}{x^2}} + \widetilde{g}(x),
\end{align*}
with
\[
|\widetilde{g}(x)|\leq \int_{x}^{+\infty} \frac{2}{t^3} \d t= \frac{1}{x^2},
\]
which shows that $g$ is Lebesgue integrable on $[1,+\infty)$.   On the other hand, since $\int_{\RR} f(t) \d t=0$ as an improper integral, we have that $g(x)=\int_{-\infty}^{x} f(t) \d t$, and we can deal with the case $x\to -\infty$ in the same way by assuming that $x\leq -1$ and performing the same integration by parts to see that $g$ is Lebesgue integrable on $(-\infty,-1]$.  Thus, $g\in\L^1(\RR)$.

It remains only to show that $\FT[g]\not\in\L^1(\RR)$.  Since $g' = f$, integration by parts shows that for any $R > 0$,
\[
\int_{-R}^{R}g(t)\e^{-\i2\pi\xi t}\d t = -\frac{1}{\i2\pi\xi}\left([g(t)\e^{-\i2\pi\xi t}]_{-R}^{R} - \int_{-R}^{R}f(t)\e^{-\i2\pi\xi t}\d t\right).
\]
Since $g\in\L^1(\RR)$, we may let $R\to+\infty$ and apply Lebesgue's dominated convergent theorem to the left side to obtain $\FT[g](\xi) = \dfrac{1}{\i2\pi\xi}\RFT[f](\xi)$, so $\FT[g]\not\in\L^1(\RR)$ by the asymptotic \eqref{eq:RFT:RL-Lemma-Fails:MainExample}.
\end{proof}

\bibliographystyle{siam}
\bibliography{Fourier-Biblio}
\end{document}